\documentclass[11pt]{article}
\usepackage[utf8]{inputenc}
\usepackage{geometry}

\usepackage[dvipsnames]{xcolor}

\usepackage{epsfig}
\usepackage{graphicx}
\usepackage{amsbsy}
\usepackage{amsmath}
\usepackage{amsfonts}
\usepackage{amssymb}
\usepackage{textcomp}
\usepackage{hyperref}
\usepackage{aliascnt}
\usepackage{bbm}
\usepackage{tikz-cd}
\usepackage{enumitem}\setlist[enumerate]{label=\textnormal{(\roman*)}}

\usepackage[section]{placeins}

\renewcommand{\leq}{\leqslant}
\renewcommand{\geq}{\geqslant}

\renewcommand{\ge}{\geqslant}
\newcommand{\pmo}{^{\pm1}}

\newcommand{\mcm}[3]{\newcommand{#1}[#2]{{\ensuremath{#3}}}} 
\mcm{\tuple}{1}{\langle #1 \rangle}
\mcm{\name}{1}{\ulcorner #1 \urcorner}
\mcm{\Nbb}{0}{\mathbb{N}}
\mcm{\N}{0}{\mathbb{N}}
\mcm{\Zbb}{0}{\mathbb{Z}}
\mcm{\Rbb}{0}{\mathbb{R}}
\mcm{\Cbb}{0}{\mathbb{C}}
\mcm{\Qbb}{0}{\mathbb{Q}}
\mcm{\Sbb}{0}{\mathbb{S}}
\mcm{\Hbb}{0}{\mathbb{H}}
\mcm{\Kbb}{0}{\mathbb{K}}
\mcm{\Acal}{0}{\cal A}
\mcm{\Bcal}{0}{\mathcal B}
\mcm{\Ccal}{0}{\cal C}
\mcm{\Dcal}{0}{\cal D}
\mcm{\Ecal}{0}{\cal E}
\mcm{\Fcal}{0}{\cal F}
\mcm{\Gcal}{0}{\cal G}
\mcm{\Hcal}{0}{\cal H}
\mcm{\Ical}{0}{\cal I}
\mcm{\Jcal}{0}{\cal J}
\mcm{\Kcal}{0}{\cal K}
\mcm{\Lcal}{0}{\cal L}
\mcm{\Mcal}{0}{\cal M}
\mcm{\Ncal}{0}{\cal N}
\mcm{\Ocal}{0}{{\cal O}}
\mcm{\Pcal}{0}{{\cal P}}
\mcm{\Qcal}{0}{{\cal Q}}
\mcm{\Rcal}{0}{{\cal R}}
\mcm{\Scal}{0}{{\cal S}}
\mcm{\Tcal}{0}{{\cal T}}
\mcm{\Ucal}{0}{{\cal U}}
\mcm{\Vcal}{0}{{\cal V}}
\mcm{\Wcal}{0}{{\cal W}}
\mcm{\Xcal}{0}{{\cal X}}
\mcm{\Ycal}{0}{{\cal Y}}
\mcm{\Zcal}{0}{{\cal Z}}
\mcm{\Mfrak}{0}{\mathfrak M}

\usepackage{amsthm}

\newcommand{\theoremize}[2]{\newaliascnt{#1}{thm} \newtheorem{#1}[#1]{#2} \aliascntresetthe{#1}}

\theoremstyle{plain}
\newtheorem{thm}{Theorem}[section]
\theoremize{lem}{Lemma}

\theoremize{skolem}{Skolem}
\theoremize{fact}{Fact}
\theoremize{setting}{Setting}

\newtheorem{claim}{Claim}[thm]
\theoremize{obs}{Observation}
\theoremize{prop}{Proposition}
\theoremize{cor}{Corollary}
\theoremize{que}{Question}
\theoremize{oque}{Open Question}
\theoremize{con}{Conjecture}

\theoremstyle{definition}
\theoremize{dfn}{Definition}
\theoremize{rem}{Remark}
\theoremize{eg}{Example}
\theoremize{exercise}{Exercise}
\theoremstyle{plain}

\newcommand{\myref}[1]{\autoref{#1}~(\nameref{#1})}

\hypersetup{
    draft = false,
    bookmarksopen=true,
    colorlinks,
    linkcolor={red!60!black},
    citecolor={green!60!black},
    urlcolor={blue!60!black}
}

\newcommand{\sm}{\setminus}

\makeatletter
\def\thm@space@setup{\thm@preskip=2pt
\thm@postskip=2pt}
\makeatother

\makeatletter
\newenvironment{pf}[1][\proofname]{\par
  \pushQED{\qed}%
  \normalfont \topsep0\p@\relax
  \trivlist
  \item[\hskip\labelsep\itshape
  #1\@addpunct{.}]\ignorespaces
}{%
  \popQED\endtrivlist\@endpefalse
}
\makeatother

\newenvironment{cproof}{\noindent\textit{Proof of Claim.}}{\phantom{!}\!\hfill$\diamondsuit$\vspace{2pt}}

\usepackage{titlesec}

\titleformat{\section}
  {\normalfont\scshape\filcenter}{\thesection}{1em}{}

\titleformat{\subsection}
 {\normalfont\scshape\filcenter}{\thesubsection}{1em}{}

\titlespacing*{\section}
{0pt}{5pt plus 1pt minus 1pt}{5pt plus 1pt minus 1pt}

\titlespacing*{\subsection}
{0pt}{5pt plus 1pt minus 1pt}{5pt plus 1pt minus 1pt}
\usepackage{tikzit}

\tikzstyle{black dot}=[fill=black, draw=none, shape=circle, scale=0.3]
\tikzstyle{backing}=[fill=white, draw=none, shape=circle, scale=0.85]
\tikzstyle{small back}=[fill=white, draw=none, shape=circle, scale=0.95]
\tikzstyle{medium back}=[fill=white, draw=none, shape=circle, scale=1.2]

\tikzstyle{c0}=[-, line width=0.08em]
\tikzstyle{c0 arrow}=[->, line width=0.08em]
\tikzstyle{c0 bold}=[-, line width=0.15em]
\tikzstyle{c1}=[-, draw={rgb,255: red,135; green,206; blue,235}, line width=0.08em]
\tikzstyle{c1 arrow}=[->, draw={rgb,255: red,135; green,206; blue,235}, line width=0.08em]
\tikzstyle{c1 bold}=[-, draw={rgb,255: red,135; green,206; blue,235}, line width=0.15em]
\tikzstyle{c2}=[-, draw={rgb,255: red,255; green,69; blue,0}, line width=0.08em]
\tikzstyle{c2 arrow}=[->, draw={rgb,255: red,255; green,69; blue,0}, line width=0.08em]
\tikzstyle{c3}=[-, draw={rgb,255: red,102; green,176; blue,50}, line width=0.08em]
\tikzstyle{c3 arrow}=[->, draw={rgb,255: red,102; green,176; blue,50}, line width=0.08em]
\tikzstyle{d1}=[-, dashed, line width=0.08em]
\tikzstyle{d1 arrow}=[->, line width=0.08em, dashed]
\tikzstyle{c1 dashed}=[-, draw={rgb,255: red,135; green,206; blue,235}, line width=0.08em, dashed, dash pattern=on 1em off 1em]
\tikzstyle{c2 dashed}=[-, draw={rgb,255: red,255; green,69; blue,0}, line width=0.08em, dashed, dash pattern=on 1em off 1em]

\usepackage{verbatim}

\newcommand{\triple}{{\scriptsize(\hspace{-0.23em}|\hspace{-0.23em})}}

\newcommand{\ct}{^{C}}
\mcm{\restric}{0}{\upharpoonright}

\newcommand{\tseq}[1]{(#1_1\minorm #1_2\minorp #1_3)}

\mcm{\minorp}{0}{\ \preccurlyeq \ }
\DeclareRobustCommand{\minorm}{\text{\reflectbox{$\ \preccurlyeq \ $}}}

\newcommand{\rcol}{\textcolor{RedOrange}{\textsc{red}}}
\newcommand{\bcol}{\textcolor{SkyBlue}{\textsc{blue}}}
\newcommand{\gcol}{\textcolor{Green}{\textsc{green}}}
\newcommand{\xcol}{\textsc{black}}
\newcommand{\se}{\subseteq}

\mcm{\Fbb}{0}{\mathbb{F}}

\title{\vspace{-1cm}A Graph Minors Approach to Temporal Sequences}
\author{Johannes Carmesin\thanks{TU Freiberg, funded by DFG, project number 546892829.} and Will J.\ Turner\footnotemark[1]}
\date{}

\begin{document}

\maketitle
\thispagestyle{empty}

\setcounter{page}{0}

\vspace{-0.75cm}
\begin{center}
\begin{minipage}{0.9\textwidth}
\small
\textbf{Abstract.} We develop a structural approach to simultaneous embeddability in temporal sequences of graphs, inspired by graph minor theory. Our main result is a classification theorem for \(2\)-connected temporal sequences: we identify five obstruction classes and show that every \(2\)-connected temporal sequence is either simultaneously embeddable or admits a sequence of improvements leading to an obstruction. This structural insight leads to a polynomial-time algorithm for deciding the simultaneous embeddability of \(2\)-connected temporal sequences.

\hspace{0.4cm}The restriction to \(2\)-connected sequences is necessary, as the problem is NP-hard for connected graphs, while trivial for \(3\)-connected graphs. As a consequence, our framework also resolves the rooted-tree SEFE problem, a natural extension of the well-studied Sunflower SEFE. 

\hspace{0.4cm}Our results uncover a rich structural theory of temporal planarity, laying the groundwork for a temporal graph minors theory.
\end{minipage}
\end{center}
\vspace{2cm}
\begin{center}
\begin{minipage}{\textwidth}
        \centering
        \Large
        \scalebox{0.65}{\input{figures/main_strategy.tikz}}
\end{minipage}
\end{center}
\vspace{2.5cm}

\pagebreak\section{Introduction}
Temporal graphs provide a natural framework for modelling networks that evolve over time. Such networks arise in diverse fields, including computer science (e.g., the internet), neuroscience (e.g., brain connectivity) and epidemiology (e.g., social networks) \cite{sporns2018graph,tang2013applications,hosseinzadeh2022temporal,
einarsson2014one, giusti2015clique}. The temporal dimension is widely believed to be essential for understanding the full behaviour of such networks \cite{thompson2017static}. As a result, temporal graphs have attracted increasing interest in recent years. 
Much of the existing work in temporal graph theory has focused on modelling, problem formulation, and heuristic approaches \cite{dall2024embedding,huang2023temporal,{GuptaBedathur2022}}, reflecting the field’s broad range of applications and recent emergence. While there is growing interest and a rich set of algorithmic questions, fundamental structural characterisations and general algorithmic frameworks offer exciting challenges \cite{{Fluschnik2020},{michail2016introduction}}. In this work, we study temporal planarity.

A fundamental landmark in graph theory is the theory of graph minors. Since Robertson and Seymour's proof of the Graph Minor Theorem and their Graph Minor Structure Theorem in 2004, the impact of graph minor theory has steadily expanded beyond the traditional graph-theoretic setting. Over the past two decades, tools from graph minor theory have proven effective in reasoning about other graph-like structures. Examples include matroid theory, where 
Geelen, Gerards and Whittle proved Rota’s Conjecture~\cite{GeelenGerardsWhittle2014}, width parameters such as rank-width by  Oum and Seymour \cite{OumSeymour2006}, directed graphs~\cite{kawarabayashi2015towards,kawarabayashi2015directed,amiri2016erdos}, and $2$-dimensional simplicial complexes \cite{carmesin2023embedding}.

Temporal graphs are typically defined as sequences of graphs on a fixed vertex set, where edges may be added or deleted over time. As we shall see, allowing contractions in addition to edge deletions enables new applications and unlocks stronger structural techniques, particularly in the planar setting. Moreover, allowing edge contractions can lead to cleaner structural results; for example, Kuratowski's Theorem provides a finite characterisation of graph planarity in terms of forbidden minors. With only the subgraph relation, one requires infinitely many obstructions.

A \emph{temporal sequence (for the minor relation)} is a sequence $(G_1, \ldots, G_n)$ of graphs such that for every $i \in [n-1]$, either $G_i$ is a minor of $G_{i+1}$, or $G_{i+1}$ is a minor of $G_i$. A temporal graph, in the classic sense, is a temporal sequence in which minor relations are restricted to edge deletions.

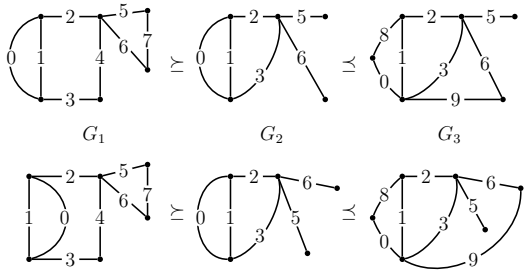
\begin{figure}
    \begin{minipage}{0.44\textwidth}
    \centering
       \large
       \centering
        \scalebox{0.63}{\begin{tikzpicture}
	\begin{pgfonlayer}{nodelayer}
		\node [style=black dot] (0) at (-7, 4.75) {};
		\node [style=black dot] (1) at (-7, 1.25) {};
		\node [style=black dot] (2) at (-2.5, 2.5) {};
		\node [style=black dot] (3) at (-4.5, 1.25) {};
		\node [style=black dot] (4) at (-4.5, 4.75) {};
		\node [style=black dot] (5) at (-2.5, 5) {};
		\node [style=black dot] (14) at (1, 4.75) {};
		\node [style=black dot] (15) at (1, 1.25) {};
		\node [style=black dot] (16) at (5, 1.25) {};
		\node [style=black dot] (18) at (3, 4.75) {};
		\node [style=black dot] (19) at (5, 4.75) {};
		\node [style=black dot] (27) at (8.25, 4.75) {};
		\node [style=black dot] (28) at (8.25, 1.25) {};
		\node [style=black dot] (29) at (12.5, 1.25) {};
		\node [style=black dot] (30) at (10.75, 4.75) {};
		\node [style=black dot] (31) at (13, 4.75) {};
		\node [style=black dot] (39) at (7, 3) {};
		\node [style=none] (41) at (6, 2.75) {$\preceq$};
		\node [style=none] (42) at (-1.25, 2.75) {$\succeq$};
		\node [style=none] (43) at (-4.75, -0.25) {$G_1$};
		\node [style=none] (44) at (2.75, -0.25) {$G_2$};
		\node [style=none] (45) at (10.25, -0.25) {$G_3$};
		\node [style=black dot] (46) at (-7.5, -2) {};
		\node [style=black dot] (47) at (-7.5, -5.5) {};
		\node [style=black dot] (48) at (-2.5, -3.75) {};
		\node [style=black dot] (49) at (-4.5, -5.5) {};
		\node [style=black dot] (50) at (-4.5, -2) {};
		\node [style=black dot] (51) at (-2.5, -1.5) {};
		\node [style=black dot] (60) at (1, -2) {};
		\node [style=black dot] (61) at (1, -5.5) {};
		\node [style=black dot] (62) at (4.25, -5.25) {};
		\node [style=black dot] (63) at (3, -2) {};
		\node [style=black dot] (64) at (5.5, -2.5) {};
		\node [style=black dot] (71) at (8.25, -2) {};
		\node [style=black dot] (72) at (8.25, -5.5) {};
		\node [style=black dot] (73) at (13.25, -2.5) {};
		\node [style=black dot] (74) at (10.5, -2) {};
		\node [style=black dot] (75) at (11.75, -4.25) {};
		\node [style=black dot] (83) at (7, -3.75) {};
		\node [style=none] (85) at (6, -3.5) {$\preceq$};
		\node [style=none] (86) at (-1.25, -3.5) {$\succeq$};
		\node [style=small back] (87) at (-8.25, 3) {};
		\node [style=none] (88) at (-8.25, 3) {$0$};
		\node [style=small back] (89) at (-7, 3) {};
		\node [style=none] (90) at (-7, 3) {$1$};
		\node [style=small back] (91) at (-5.75, 4.75) {};
		\node [style=none] (92) at (-5.75, 4.75) {$2$};
		\node [style=small back] (93) at (-5.75, 1.25) {};
		\node [style=none] (94) at (-5.75, 1.25) {$3$};
		\node [style=small back] (95) at (-4.5, 3) {};
		\node [style=none] (96) at (-4.5, 3) {$4$};
		\node [style=small back] (97) at (-3.5, 5) {};
		\node [style=none] (98) at (-3.5, 5) {$5$};
		\node [style=small back] (99) at (-3.5, 3.5) {};
		\node [style=none] (100) at (-3.5, 3.5) {$6$};
		\node [style=small back] (101) at (-2.5, 3.75) {};
		\node [style=none] (102) at (-2.5, 3.75) {$7$};
		\node [style=small back] (103) at (2, 4.75) {};
		\node [style=none] (104) at (2, 4.75) {$2$};
		\node [style=small back] (105) at (4, 4.75) {};
		\node [style=none] (106) at (4, 4.75) {$5$};
		\node [style=small back] (107) at (4, 3) {};
		\node [style=none] (108) at (4, 3) {$6$};
		\node [style=small back] (109) at (-0.25, 3) {};
		\node [style=none] (110) at (-0.25, 3) {$0$};
		\node [style=small back] (111) at (1, 3) {};
		\node [style=none] (112) at (1, 3) {$1$};
		\node [style=small back] (113) at (2.25, 2.25) {};
		\node [style=none] (114) at (2.25, 2.25) {$3$};
		\node [style=small back] (115) at (9.5, 4.75) {};
		\node [style=none] (116) at (9.5, 4.75) {$2$};
		\node [style=small back] (117) at (12, 4.75) {};
		\node [style=none] (118) at (12, 4.75) {$5$};
		\node [style=small back] (119) at (7.5, 4) {};
		\node [style=none] (120) at (7.5, 4) {$8$};
		\node [style=small back] (121) at (7.5, 2) {};
		\node [style=none] (122) at (7.5, 2) {$0$};
		\node [style=small back] (123) at (8.25, 3) {};
		\node [style=none] (124) at (8.25, 3) {$1$};
		\node [style=small back] (125) at (11.75, 2.75) {};
		\node [style=none] (126) at (11.75, 2.75) {$6$};
		\node [style=small back] (127) at (10, 2.5) {};
		\node [style=none] (128) at (10, 2.5) {$3$};
		\node [style=small back] (129) at (10.5, 1.25) {};
		\node [style=none] (130) at (10.5, 1.25) {$9$};
		\node [style=small back] (131) at (-6, -3.75) {};
		\node [style=none] (132) at (-6, -3.75) {$0$};
		\node [style=small back] (133) at (-7.5, -3.75) {};
		\node [style=none] (134) at (-7.5, -3.75) {$1$};
		\node [style=small back] (135) at (-5.75, -2) {};
		\node [style=none] (136) at (-5.75, -2) {$2$};
		\node [style=small back] (137) at (-5.75, -5.5) {};
		\node [style=none] (138) at (-5.75, -5.5) {$3$};
		\node [style=small back] (139) at (-4.5, -3.75) {};
		\node [style=none] (140) at (-4.5, -3.75) {$4$};
		\node [style=small back] (141) at (-3.5, -1.75) {};
		\node [style=none] (142) at (-3.5, -1.75) {$5$};
		\node [style=small back] (143) at (-3.5, -3) {};
		\node [style=none] (144) at (-3.5, -3) {$6$};
		\node [style=small back] (145) at (-2.5, -2.75) {};
		\node [style=none] (146) at (-2.5, -2.75) {$7$};
		\node [style=small back] (147) at (2, -2) {};
		\node [style=none] (148) at (2, -2) {$2$};
		\node [style=small back] (149) at (-0.25, -3.75) {};
		\node [style=none] (150) at (-0.25, -3.75) {$0$};
		\node [style=small back] (151) at (1, -3.75) {};
		\node [style=none] (152) at (1, -3.75) {$1$};
		\node [style=small back] (153) at (2.25, -4.5) {};
		\node [style=none] (154) at (2.25, -4.5) {$3$};
		\node [style=small back] (155) at (3.75, -3.75) {};
		\node [style=none] (156) at (3.75, -3.75) {$5$};
		\node [style=small back] (157) at (4.25, -2.25) {};
		\node [style=none] (158) at (4.25, -2.25) {$6$};
		\node [style=small back] (159) at (7.5, -2.75) {};
		\node [style=none] (160) at (7.5, -2.75) {$8$};
		\node [style=small back] (161) at (7.5, -4.75) {};
		\node [style=none] (162) at (7.5, -4.75) {$0$};
		\node [style=small back] (163) at (8.25, -3.75) {};
		\node [style=none] (164) at (8.25, -3.75) {$1$};
		\node [style=small back] (165) at (9.25, -2) {};
		\node [style=none] (166) at (9.25, -2) {$2$};
		\node [style=small back] (167) at (10, -4) {};
		\node [style=none] (168) at (10, -4) {$3$};
		\node [style=small back] (169) at (11.25, -3.25) {};
		\node [style=none] (170) at (11.25, -3.25) {$5$};
		\node [style=small back] (171) at (12, -2.25) {};
		\node [style=none] (172) at (12, -2.25) {$6$};
		\node [style=small back] (173) at (11.25, -5.5) {};
		\node [style=none] (174) at (11.25, -5.5) {$9$};
	\end{pgfonlayer}
	\begin{pgfonlayer}{edgelayer}
		\draw [style=c0] (0) to (4);
		\draw [style=c0] (4) to (3);
		\draw [style=c0] (3) to (1);
		\draw [style=c0] (1) to (0);
		\draw [style=c0, bend right=75, looseness=1.25] (0) to (1);
		\draw [style=c0] (4) to (2);
		\draw [style=c0] (2) to (5);
		\draw [style=c0] (5) to (4);
		\draw [style=c0] (14) to (18);
		\draw [style=c0] (15) to (14);
		\draw [style=c0, bend right=75, looseness=1.25] (14) to (15);
		\draw [style=c0] (18) to (16);
		\draw [style=c0] (19) to (18);
		\draw [style=c0, bend right=45, looseness=0.75] (15) to (18);
		\draw [style=c0] (27) to (30);
		\draw [style=c0] (28) to (27);
		\draw [style=c0] (30) to (29);
		\draw [style=c0] (31) to (30);
		\draw [style=c0, bend right=45, looseness=0.75] (28) to (30);
		\draw [style=c0] (28) to (29);
		\draw [style=c0, bend right=15, looseness=0.75] (27) to (39);
		\draw [style=c0, bend right=15] (39) to (28);
		\draw [style=c0] (46) to (50);
		\draw [style=c0] (50) to (49);
		\draw [style=c0] (49) to (47);
		\draw [style=c0] (47) to (46);
		\draw [style=c0, bend left=75, looseness=1.50] (46) to (47);
		\draw [style=c0] (50) to (48);
		\draw [style=c0] (48) to (51);
		\draw [style=c0] (51) to (50);
		\draw [style=c0] (60) to (63);
		\draw [style=c0] (61) to (60);
		\draw [style=c0, bend right=90, looseness=1.25] (60) to (61);
		\draw [style=c0] (63) to (62);
		\draw [style=c0] (64) to (63);
		\draw [style=c0, bend right=45, looseness=0.75] (61) to (63);
		\draw [style=c0] (71) to (74);
		\draw [style=c0] (72) to (71);
		\draw [style=c0] (74) to (73);
		\draw [style=c0] (75) to (74);
		\draw [style=c0, bend right=45, looseness=0.75] (72) to (74);
		\draw [style=c0, bend right=60] (72) to (73);
		\draw [style=c0, bend right=15, looseness=0.75] (71) to (83);
		\draw [style=c0, bend right=15, looseness=0.75] (83) to (72);
	\end{pgfonlayer}
\end{tikzpicture}}
    \end{minipage}
    \hfill%
    \begin{minipage}{0.54\textwidth}
        \caption{A temporal triple $(G_1,G_2,G_3)$. We obtain $G_2$ from $G_1$ by contracting the edge $4$ and deleting $7$ and from $G_3$ by contracting the edge $8$ and deleting~$9$. Depicted are two possible sequences of embeddings, one of which (top) is simultaneous whereas the other (bottom) fails to be consistent.}
    \label{fig:to_sim_or_not_to_sim}
    \end{minipage}
    \vspace{-.25cm}
\end{figure}

In this paper, we investigate an analogue of graph planarity for temporal sequences. We say that a temporal sequence $(G_1, \ldots, G_n)$ is \emph{simultaneously embeddable} if there is a respective sequence $(\iota_1, \ldots, \iota_n)$ of planar embeddings that are consistent with respect to the minor relation; that is, if $G_i$ is a minor of $G_j$ for $j\in \{i+1,i-1\}$, then the embedding $\iota_i$ is obtained from the embedding $\iota_j$ of $G_j$ by applying the respective minoring to $\iota_j$. In \autoref{fig:to_sim_or_not_to_sim} we depict one temporal sequence along with two possible sequences of (plane) embeddings. Later, in Examples~\ref{eg:parity_obs} and \ref{eg:amazingly_bad}, we will see temporal sequences which admit no simultaneous embedding.

There are three main motivations for studying a temporal analogue of planarity. First, graph minor theory originated with Kuratowski’s theorem, so it is not unreasonable to expect that a graph minor theory for temporal graphs might emerge from the planar case. Second, planarity serves as a test case for the challenges inherent in temporal sequences: even in the deletion-only setting, deciding simultaneous embeddability can quickly become non-trivial. 
Third, our notion of simultaneous embedding generalises \emph{SEFE} (simultaneous embedding with fixed edges), providing a unified framework for problems of this kind; see the survey \cite{rutter2020simultaneous} for an overview of the extensive literature on SEFE and related notions. As one application, our methods extend to problems previously studied in this setting, yielding new algorithmic consequences.

Working in the context of the SEFE problems,  Gassner, Jünger, Percan, Schaefer and Schulz showed that deciding simultaneous embeddability is NP-hard for temporal sequences of connected graphs \cite{gassner2006simultaneous}. 
Given that the $3$-connected case is a straightforward application of Kuratowski's theorem, NP-hardness of the connected case suggests that the natural problem to look at is simultaneous embeddability for temporal sequences of \(2\)-connected graphs; we refer to them as \emph{\(2\)-connected temporal sequences}.

The main result of this paper given below provides a list of obstructions to the simultaneous embeddability of \(2\)-connected temporal sequences that completely characterises the problem, perhaps after applying a finite number of simplifying operations. To begin, we present two examples illustrating distinct ways in which a temporal sequence may fail to be simultaneously embeddable.

\begin{eg}\label{eg:parity_obs}
    Consider the temporal triple $(G_1,G_2,G_3)$ illustrated in \autoref{fig:parity_obs}. Observe that in all embeddings induced on $G_2$ by embeddings of $G_1$, the edges $1$ and $3$ appear on the boundary of the same face. On the other hand, in all embeddings of $G_2$ induced by embeddings of $G_3$, they lie in distinct face boundaries. Hence no embedding of $G_2$ is induced by embeddings of both $G_1$ and $G_3$. Thus $(G_1,G_2,G_3)$ admits no simultaneous embedding.
\end{eg}

\begin{eg}\label{eg:amazingly_bad}
    Consider the temporal triple $(G_1,G_2,G_3)$ depicted in \autoref{fig:amazingly_bad_obs}. Observe that in any embedding of $G_1$, the edge $2$ has to appear directly between the edges $1$ and $3$ in the cyclic order at $x_1$. On the other hand, in any embedding of $G_3$, the edges $1$ and $3$ must be adjacent in the cyclic order at $x_3$. Hence there is no embedding of $G_2$ which is induced both by an embedding of $G_1$ and an embedding of $G_3$, since they disagree on the cyclic order at $x_2$. Hence $(G_1,G_2,G_3)$ admits no simultaneous embedding.
\end{eg}

\begin{figure}
    \begin{minipage}{0.49\textwidth}
        \centering
        \Large
        \scalebox{0.55}{\begin{tikzpicture}
	\begin{pgfonlayer}{nodelayer}
		\node [style=none] (42) at (-3.5, 0) {$\succeq$};
		\node [style=black dot] (87) at (-10, 0) {};
		\node [style=black dot] (88) at (-12.5, 0) {};
		\node [style=black dot] (89) at (-10, 1.5) {};
		\node [style=black dot] (90) at (-12.5, 2) {};
		\node [style=black dot] (91) at (-12.5, -2) {};
		\node [style=black dot] (92) at (-10, -1.5) {};
		\node [style=black dot] (93) at (-8, 4) {};
		\node [style=black dot] (94) at (-8, -4) {};
		\node [style=black dot] (95) at (-7, 0) {};
		\node [style=black dot] (96) at (-4.5, 0) {};
		\node [style=none] (101) at (-0.5, 0) {$x$};
		\node [style=black dot] (103) at (-1.25, 0) {};
		\node [style=black dot] (104) at (0, 1.5) {};
		\node [style=black dot] (105) at (-2.5, 2) {};
		\node [style=black dot] (106) at (-2.5, -2) {};
		\node [style=black dot] (107) at (0, -1.5) {};
		\node [style=black dot] (108) at (2, 4) {};
		\node [style=black dot] (109) at (2, -4) {};
		\node [style=black dot] (110) at (3, 0) {};
		\node [style=black dot] (111) at (5.5, 0) {};
		\node [style=black dot] (116) at (10, 0) {};
		\node [style=black dot] (117) at (7.5, 0) {};
		\node [style=black dot] (118) at (10, 1.5) {};
		\node [style=black dot] (119) at (7.5, 2) {};
		\node [style=black dot] (120) at (7.5, -2) {};
		\node [style=black dot] (121) at (10, -1.5) {};
		\node [style=black dot] (122) at (12, 4) {};
		\node [style=black dot] (123) at (12, -4) {};
		\node [style=black dot] (124) at (13, 0) {};
		\node [style=black dot] (125) at (15.5, 0) {};
		\node [style=none] (127) at (6.5, 0) {$\preceq$};
		\node [style=small back] (132) at (-10.25, 3) {};
		\node [style=none] (133) at (-10.25, 3) {$1$};
		\node [style=small back] (134) at (-9.25, 2.5) {};
		\node [style=none] (135) at (-9.25, 2.5) {$2$};
		\node [style=small back] (136) at (-0.25, 3) {};
		\node [style=none] (137) at (-0.25, 3) {$1$};
		\node [style=small back] (138) at (0.75, 2.5) {};
		\node [style=none] (139) at (0.75, 2.5) {$2$};
		\node [style=small back] (140) at (-10.25, -3) {};
		\node [style=none] (141) at (-10.25, -3) {$3$};
		\node [style=small back] (142) at (-9.25, -2.5) {};
		\node [style=none] (143) at (-9.25, -2.5) {$4$};
		\node [style=small back] (144) at (-0.25, -3) {};
		\node [style=none] (145) at (-0.25, -3) {$3$};
		\node [style=small back] (146) at (0.75, -2.5) {};
		\node [style=none] (147) at (0.75, -2.5) {$4$};
		\node [style=small back] (148) at (-11.25, 0) {};
		\node [style=none] (149) at (-11.25, 0) {$5$};
		\node [style=small back] (150) at (9.75, 3) {};
		\node [style=none] (151) at (9.75, 3) {$1$};
		\node [style=small back] (152) at (10.75, 2.5) {};
		\node [style=none] (153) at (10.75, 2.5) {$2$};
		\node [style=small back] (154) at (9.75, -3) {};
		\node [style=none] (155) at (9.75, -3) {$4$};
		\node [style=small back] (156) at (10.75, -2.5) {};
		\node [style=none] (157) at (10.75, -2.5) {$3$};
		\node [style=small back] (158) at (8.75, 0) {};
		\node [style=none] (159) at (8.75, 0) {$6$};
	\end{pgfonlayer}
	\begin{pgfonlayer}{edgelayer}
		\draw [style=c0] (90) to (89);
		\draw [style=c0] (92) to (91);
		\draw [style=c0] (88) to (87);
		\draw [style=c0] (89) to (92);
		\draw [style=c0] (91) to (90);
		\draw [style=c0] (90) to (93);
		\draw [style=c0] (93) to (89);
		\draw [style=c0] (92) to (94);
		\draw [style=c0] (94) to (91);
		\draw [style=c0] (93) to (95);
		\draw [style=c0] (96) to (93);
		\draw [style=c0] (96) to (95);
		\draw [style=c0] (95) to (94);
		\draw [style=c0] (96) to (94);
		\draw [style=c0] (105) to (104);
		\draw [style=c0] (107) to (106);
		\draw [style=c0] (105) to (108);
		\draw [style=c0] (108) to (104);
		\draw [style=c0] (107) to (109);
		\draw [style=c0] (109) to (106);
		\draw [style=c0] (108) to (110);
		\draw [style=c0] (111) to (108);
		\draw [style=c0] (111) to (110);
		\draw [style=c0] (110) to (109);
		\draw [style=c0] (111) to (109);
		\draw [style=c0] (105) to (103);
		\draw [style=c0] (103) to (106);
		\draw [style=c0] (103) to (107);
		\draw [style=c0] (103) to (104);
		\draw [style=c0] (119) to (118);
		\draw [style=c0] (121) to (120);
		\draw [style=c0] (117) to (116);
		\draw [style=c0] (118) to (121);
		\draw [style=c0] (120) to (119);
		\draw [style=c0] (119) to (122);
		\draw [style=c0] (122) to (118);
		\draw [style=c0] (121) to (123);
		\draw [style=c0] (123) to (120);
		\draw [style=c0] (122) to (124);
		\draw [style=c0] (125) to (122);
		\draw [style=c0] (125) to (124);
		\draw [style=c0] (124) to (123);
		\draw [style=c0] (125) to (123);
	\end{pgfonlayer}
\end{tikzpicture}}
        \caption{A temporal triple $(G_1,G_2,G_3)$ which has no simultaneous embedding. The graph $G_2$ is obtained from $G_1$ by contracting the edge $5$ and is obtained from $G_3$ by contracting the edge $6$. Although $G_1$ and $G_3$ are isomorphic, they minor down to $G_2$ in different ways, as can be seen by considering the edges labelled $1,2,3$ and~$4$.}
        \label{fig:parity_obs}
    \end{minipage}
    \hfill%
    \begin{minipage}{0.49\textwidth}
        \centering
        \Large
        \scalebox{0.55}{\begin{tikzpicture}
	\begin{pgfonlayer}{nodelayer}
		\node [style=none] (42) at (-3.5, 0) {$\succeq$};
		\node [style=black dot] (89) at (-7.5, 0) {};
		\node [style=black dot] (90) at (-10, 0) {};
		\node [style=black dot] (93) at (-7.5, 4) {};
		\node [style=black dot] (94) at (-7.5, -4) {};
		\node [style=black dot] (95) at (-5.75, 0) {};
		\node [style=black dot] (96) at (-4.5, 0) {};
		\node [style=none] (127) at (6.5, 0) {$\preceq$};
		\node [style=black dot] (132) at (-12.5, 0) {};
		\node [style=none] (133) at (-7.5, 4.75) {$x_1$};
		\node [style=none] (134) at (-7.5, -4.75) {$y_1$};
		\node [style=none] (135) at (-10.75, 3) {$A$};
		\node [style=black dot] (140) at (2.5, 0) {};
		\node [style=black dot] (141) at (0, 0) {};
		\node [style=black dot] (142) at (2.5, 4) {};
		\node [style=black dot] (143) at (2.5, -4) {};
		\node [style=black dot] (144) at (4.25, 0) {};
		\node [style=black dot] (145) at (5.5, 0) {};
		\node [style=black dot] (147) at (-2.5, 0) {};
		\node [style=none] (148) at (2.5, 4.75) {$x_2$};
		\node [style=none] (149) at (2.5, -4.75) {$y_2$};
		\node [style=black dot] (153) at (12.5, 0) {};
		\node [style=black dot] (154) at (10, 0) {};
		\node [style=black dot] (155) at (12.5, 4) {};
		\node [style=black dot] (156) at (12.5, -4) {};
		\node [style=black dot] (157) at (14.25, 0) {};
		\node [style=black dot] (158) at (15.5, 0) {};
		\node [style=black dot] (159) at (7.5, 0) {};
		\node [style=none] (160) at (12.5, 4.75) {$x_3$};
		\node [style=none] (161) at (12.5, -4.75) {$y_3$};
		\node [style=small back] (163) at (-10.25, 1.75) {};
		\node [style=none] (164) at (-10.25, 1.75) {$1$};
		\node [style=small back] (165) at (-9, 1.75) {};
		\node [style=none] (166) at (-9, 1.75) {$2$};
		\node [style=small back] (167) at (-7.5, 1.75) {};
		\node [style=none] (168) at (-7.5, 1.75) {$3$};
		\node [style=small back] (169) at (-11, 0) {};
		\node [style=none] (170) at (-11, 0) {$e$};
		\node [style=small back] (171) at (-8.75, 0) {};
		\node [style=none] (172) at (-8.75, 0) {$f$};
		\node [style=small back] (173) at (-0.25, 1.75) {};
		\node [style=none] (174) at (-0.25, 1.75) {$1$};
		\node [style=small back] (175) at (1, 1.75) {};
		\node [style=none] (176) at (1, 1.75) {$2$};
		\node [style=small back] (177) at (2.5, 1.75) {};
		\node [style=none] (178) at (2.5, 1.75) {$3$};
		\node [style=small back] (179) at (9.75, 1.75) {};
		\node [style=none] (180) at (9.75, 1.75) {$1$};
		\node [style=small back] (181) at (11, 1.75) {};
		\node [style=none] (182) at (11, 1.75) {$2$};
		\node [style=small back] (183) at (12.5, 1.75) {};
		\node [style=none] (184) at (12.5, 1.75) {$3$};
		\node [style=small back] (185) at (11.25, -0.75) {};
		\node [style=none] (186) at (11.25, -0.75) {$g$};
	\end{pgfonlayer}
	\begin{pgfonlayer}{edgelayer}
		\draw [style=c1] (90) to (89);
		\draw [style=c1] (90) to (93);
		\draw [style=c1] (93) to (89);
		\draw [style=c0] (93) to (95);
		\draw [style=c0] (96) to (93);
		\draw [style=c0] (96) to (95);
		\draw [style=c0] (95) to (94);
		\draw [style=c0] (96) to (94);
		\draw [style=c1] (89) to (94);
		\draw [style=c1] (94) to (90);
		\draw [style=c1] (90) to (132);
		\draw [style=c1] (132) to (94);
		\draw [style=c1] (132) to (93);
		\draw [style=c0] (141) to (142);
		\draw [style=c0] (142) to (140);
		\draw [style=c0] (142) to (144);
		\draw [style=c0] (145) to (142);
		\draw [style=c0] (145) to (144);
		\draw [style=c0] (144) to (143);
		\draw [style=c0] (145) to (143);
		\draw [style=c0] (140) to (143);
		\draw [style=c0] (143) to (141);
		\draw [style=c0] (147) to (143);
		\draw [style=c0] (147) to (142);
		\draw [style=c0] (154) to (155);
		\draw [style=c2] (155) to (153);
		\draw [style=c0] (155) to (157);
		\draw [style=c0] (158) to (155);
		\draw [style=c0] (158) to (157);
		\draw [style=c0] (157) to (156);
		\draw [style=c0] (158) to (156);
		\draw [style=c2] (153) to (156);
		\draw [style=c0] (156) to (154);
		\draw [style=c2] (159) to (156);
		\draw [style=c2] (159) to (155);
		\draw [style=c2, in=-135, out=-15] (159) to (153);
	\end{pgfonlayer}
\end{tikzpicture}}
        \caption{A temporal triple $(G_1,G_2,G_3)$ which has no simultaneous embedding. The graph $G_2$ is obtained from $G_1$ by deleting the edges $e$ and $f$ and is obtained from $G_3$ by deleting the edge~$g$.}
        \label{fig:amazingly_bad_obs}
    \end{minipage}
    \vspace{-.2cm}
\end{figure}

Later, we will introduce two broader classes of obstructions, \emph{disagreeable obstructions} and \emph{amazingly bad obstructions}, which generalise the examples seen earlier: \autoref{eg:parity_obs} belongs to the former, while \autoref{eg:amazingly_bad} belongs to the latter. Together with \emph{gaudy obstructions}, these will complete our list of obstructions that are themselves temporal sequence---each consisting of three graphs. However, as often occurs when proving an inductive result that does not quite hold, a slight generalisation is required---much like how the study of colourings naturally leads to list-colourings. The same necessity arises here, leading us to extend our framework further as follows.

Roughly speaking, a \emph{list} \(L\) for a temporal sequence \( (G_1, \dots, G_n) \) encodes the relative orientations of the induced embeddings of maximal \(3\)-connected minors throughout the sequence. A \emph{simultaneous list-embedding} of \( (G_1, \dots, G_n) \) with \(L\) is a simultaneous embedding that respects these relative orientations. We introduce \emph{inconsistent} and \emph{uncombinable obstructions} as concrete examples where no such embedding exists. A temporal sequence with a family of lists is \emph{obstructed} if it contains one of five types of obstructions: \emph{disagreeable}, \emph{amazingly-bad}, \emph{gaudy}, \emph{inconsistent}, or \emph{uncombinable}.
An \emph{improvement} of a temporal sequence \(\mathcal{G}\) with lists is a temporal sequence \(\mathcal{G'}\) with lists so that: $\mathcal{G'}$ is \(2\)-connected, simultaneously list-embeddable if and only if \(\mathcal{G}\) is, and \lq strictly smaller\rq\ than $\mathcal{G}$---for this we introduce an nonnegative integer-valued measure \lq $\mathbb{C}$\rq\ of complexity for temporal sequences. The main result of this paper is the following.
\begin{thm}\label{main-intro37}
    Every nonempty $2$-connected temporal sequence with two lists is obstructed or admits an improvement.
\end{thm}
For brevity in the definition of improvement which appears below, we omit any conditions which relate the structure of a temporal sequence with that of its improvement. It is important to note however, that such a structure can be recovered by inspecting the key lemmas involved in the proof. In particular, improvements are always simple---in the sense that they act on graphs individually without considering their neighbours---or take the highly-structured form depicted in \autoref{fig:twist}.

\begin{figure}
    \centering
    {\large
    \scalebox{0.9}{\begin{tikzpicture}
	\begin{pgfonlayer}{nodelayer}
		\node [style=none] (491) at (0.75, 0.25) {};
		\node [style=none] (492) at (0.75, -0.25) {};
		\node [style=none] (493) at (0.75, 0) {};
		\node [style=none] (494) at (2.5, 0) {};
		\node [style=none] (495) at (2.25, 0.25) {};
		\node [style=none] (496) at (2.25, -0.25) {};
		\node [style=none] (603) at (1.5, 0.75) {improvement};
		\node [style=none] (611) at (-12, 1.25) {$G_{i-2}$};
		\node [style=none] (612) at (-10.5, -1.25) {$G_{i-1}$};
		\node [style=none] (613) at (-9, 1.25) {$G_i$};
		\node [style=none] (614) at (-13.75, -0.5) {$\ldots$};
		\node [style=none] (615) at (-12, 0.75) {};
		\node [style=none] (616) at (-10.5, -0.75) {};
		\node [style=none] (617) at (-12, 0.75) {};
		\node [style=none] (618) at (-13, -0.25) {};
		\node [style=none] (619) at (-10.5, -0.75) {};
		\node [style=none] (620) at (-9, 0.75) {};
		\node [style=none] (621) at (-9, 0.75) {};
		\node [style=none] (622) at (-7.5, -0.75) {};
		\node [style=none] (623) at (-7.5, -0.75) {};
		\node [style=none] (624) at (-6, 0.75) {};
		\node [style=none] (625) at (-7.5, -1.25) {$G_{i+1}$};
		\node [style=none] (626) at (-6, 1.25) {$G_{i+2}$};
		\node [style=none] (627) at (-6, 0.75) {};
		\node [style=none] (628) at (-4.5, -0.75) {};
		\node [style=none] (629) at (-4.5, -0.75) {};
		\node [style=none] (630) at (-3, 0.75) {};
		\node [style=none] (631) at (-4.5, -1.25) {$G_{i+3}$};
		\node [style=none] (632) at (-3, 1.25) {$G_{i+4}$};
		\node [style=none] (633) at (-1.75, -0.5) {$\ldots$};
		\node [style=none] (634) at (-3, 0.75) {};
		\node [style=none] (635) at (-2, -0.25) {};
		\node [style=none] (636) at (6, 1.25) {$G_{i-2}$};
		\node [style=none] (637) at (7.5, -1.25) {$G_{i-1}$};
		\node [style=none] (638) at (9.25, -0.75) {$H_i$};
		\node [style=none] (639) at (4.75, -0.5) {$\ldots$};
		\node [style=none] (640) at (6, 0.75) {};
		\node [style=none] (641) at (7.5, -0.75) {};
		\node [style=none] (642) at (6, 0.75) {};
		\node [style=none] (643) at (5, -0.25) {};
		\node [style=none] (644) at (7.5, -0.75) {};
		\node [style=none] (645) at (9.25, 0) {};
		\node [style=none] (646) at (9.25, 0) {};
		\node [style=none] (647) at (10.75, -1.25) {};
		\node [style=none] (648) at (10.75, -0.75) {};
		\node [style=none] (649) at (12.25, 1.5) {};
		\node [style=none] (650) at (10.75, -1.75) {$H_{i+1}$};
		\node [style=none] (651) at (12.25, 2) {$H_{i+2}$};
		\node [style=none] (652) at (12.25, 1.5) {};
		\node [style=none] (653) at (13.75, -0.75) {};
		\node [style=none] (654) at (13.75, -0.75) {};
		\node [style=none] (655) at (15.25, 0.75) {};
		\node [style=none] (656) at (13.75, -1.25) {$G_{i+3}$};
		\node [style=none] (657) at (15.25, 1.25) {$G_{i+4}$};
		\node [style=none] (658) at (17, -0.5) {$\ldots$};
		\node [style=none] (659) at (15.25, 0.75) {};
		\node [style=none] (660) at (16.25, -0.25) {};
		\node [style=none] (661) at (9.25, 0.75) {};
		\node [style=none] (662) at (12.25, 0.75) {};
		\node [style=none] (663) at (9.25, 1) {};
		\node [style=none] (664) at (9.25, 1.5) {};
		\node [style=none] (665) at (12.25, 0.25) {};
		\node [style=none] (666) at (12.25, -0.25) {};
	\end{pgfonlayer}
	\begin{pgfonlayer}{edgelayer}
		\draw [style=c0] (491.center) to (492.center);
		\draw [style=c0] (493.center) to (494.center);
		\draw [style=c0] (494.center) to (495.center);
		\draw [style=c0] (494.center) to (496.center);
		\draw [style=c0] (615.center) to (616.center);
		\draw [style=c0] (619.center) to (620.center);
		\draw [style=c0] (618.center) to (617.center);
		\draw [style=c0] (621.center) to (622.center);
		\draw [style=c0] (623.center) to (624.center);
		\draw [style=c0] (627.center) to (628.center);
		\draw [style=c0] (629.center) to (630.center);
		\draw [style=c0] (635.center) to (634.center);
		\draw [style=c0] (640.center) to (641.center);
		\draw [style=c0] (644.center) to (645.center);
		\draw [style=c0] (643.center) to (642.center);
		\draw [style=c0] (646.center) to (647.center);
		\draw [style=c0] (652.center) to (653.center);
		\draw [style=c0] (654.center) to (655.center);
		\draw [style=c0] (660.center) to (659.center);
		\draw [style=c0] (652.center) to (647.center);
		\draw [style=c0 arrow] (666.center) to (665.center);
		\draw [style=c0 arrow] (664.center) to (663.center);
		\draw [style=d1] (654.center)
			 to (662.center)
			 to (648.center)
			 to (661.center)
			 to (644.center);
	\end{pgfonlayer}
\end{tikzpicture}}}\vspace{-.2cm}
    \caption{A temporal sequence $\mathcal{G}:=(G_1,\ldots,G_n)$ and its improvement $\mathcal{G}'$. The improvement has replaced the graphs $G_i,G_{i+1},G_{i+2}$ by $H_i,H_{i+1},H_{i+2}$, respectively. These graphs satisfy the structural conditions $H_i\preceq G_i$, $G_{i+2}\preceq H_{i+2}$ and $H_i\succeq H_{i+1}\preceq H_{i+2}$. Furthermore, $\Cbb(\mathcal{G})<\Cbb(\mathcal{G}')$. Not all of the graphs get `smaller', but the overall complexity decreases.}
    \label{fig:twist}
    \vspace{-.3cm}
\end{figure}

Beyond its direct implications for temporal graphs, this result provides a new setting where graph minor techniques yield a complete structural characterisation and efficient algorithms, demonstrating the broader applicability of these methods. Applying \autoref{main-intro37} iteratively, we obtain that every nonempty \(2\)-connected temporal sequence is either simultaneously embeddable or admits a sequence of improvements leading to an obstruction in one of our five classes. This leads to the following algorithmic consequence.

\begin{cor}\label{algo-intro37}
    There exists a polynomial-time algorithm that decides the simultaneous embeddability of temporal sequences of \(2\)-connected graphs.
\end{cor}

Perhaps surprisingly, \autoref{algo-intro37} also implies (\autoref{rooted_tree_sefe}) that the $2$-connected rooted-tree SEFE problem is in~P, a natural extension of the well-studied $2$-connected Sunflower SEFE problem.

The remainder of the paper is organised as follows. In \autoref{sec:new_prelims}, we formalise the notions of temporal sequences, simultaneous embeddings, and lists. Our main results are stated in \autoref{sec:main}, while in \autoref{sec:obs} we define the five obstruction classes that underpin our classification theorem. The core of the proof is carried out in Sections~\ref{sec:nice}, \ref{sec:amazing_analogues}, and~\ref{sec:delightful_analogues}, where we establish three key lemmas. Finally, in \autoref{sec:conc} we give algorithmic consequences,  discuss follow-up works such as \cite{contractionsequences2025}, \cite{fpttreesofgraphs2025} and \cite{weaktemporal2025}, as well as new research directions opened by \autoref{thm:main}.

\section{Temporal sequences}\label{sec:new_prelims}
For general graph-theoretic and graph-embedding terminology, we follow the book of Mohar and Thomassen \cite{graphsonsurfaces2001}. There are some subtle issues concerning the definition of simultaneous embedding below (see \autoref{rem:weak_remark} for an explanation). 
These issues will be circumvented by working with the following notation. 

A \emph{graph} consists of a \emph{vertex set} $V$ and an \emph{edge set} $E$, together with an incidence relation between the vertices $v\in V$ and the edges $e\in E$ such that each edge is incident with at least one and at most two vertices. We shall throughout assume that $V$ and $E$ are disjoint sets and note that graphs are allowed to have parallel edges and loops. A \emph{minor} of a graph $G$ is a graph $H$ and a choice of disjoint edge sets $C,D\subseteq E(G)$ so that $H$ is obtained from $G$ by deleting the edges in $D$, contracting the edges in $C$, deleting some isolated vertices, and relabelling the vertices (but not the edges). If $H$ is a minor of $G$, we write $H\minorp G$ (and assume sets $C$ and $D$ are fixed).
Let $(G_1,\ldots,G_n)$ be a sequence of graphs, whose edge sets possibly overlap. We say that $(G_1,\ldots,G_n)$ is a \emph{temporal sequence (for the minor relation)} if for every $i\in [n-1]$ we have that either $G_i\minorp G_{i+1}$ or the other way round; that is, $G_{i+1}\minorp G_{i}$.

A \emph{rotation system} is a combinatorial way to describe an embedding of a connected graph $G$ in a $2$-dimensional manifold (such the oriented sphere $\Sbb^2$): it consists of a family $(\sigma(v)| v\in V(G))$ of cyclic orderings $\sigma(v)$ of the edges incident with the vertex $v$. For example, given a graph embedded in $\Sbb^2$, we obtain a rotation system by \lq walking clockwise around each vertex $v$\rq\ and storing the order in which we traverse the edges incident to $v$. 
We refer to this rotation system as the rotation system \emph{induced} by that embedding of $G$. 
Conversely, given a rotation system $\Sigma$ of a connected graph $G$, up to continuous deformation, there is at most one embedding of $G$ in $\Sbb^2$ that induces $\Sigma$. And if there is one, then we call $\Sigma$ planar. 

Let $G$ be a graph and $\sigma$ a planar rotation system (embedding) thereof.
The minor relation on graphs lifts to a relation on graphs with rotation systems as follows. Given a pair $(G,\sigma)$ consisting of a graph $G$ and a rotation system $\sigma$, and an edge $e$ of $G$, the graph with a rotation system obtained from $(G,\sigma)$ by \emph{deleting} $e$ is the pair $(G-e,\sigma')$ where for each vertex $v$, we let $\sigma'(v)=\sigma(v)-e$  be the rotator obtained from $\sigma(v)$ by removing $e$ from the cyclic order (if it is present). And the graph with a rotation system obtained from $(G,\sigma)$ by \emph{contracting} an edge $e=vw$ is the pair $(G/e,\sigma')$ where $\sigma'(u)=\sigma(u)$ for each vertex $u$ in $G$ distinct from both $v$ and $w$ and $\sigma'(e)$ is defined by concatenating\footnote{Let $S_v$ and $S_w$ denote the segments from $\sigma(v)$ and $\sigma(w)$, respectively. We take the unique cyclic order containing both $S_v$ and $S_w$ as segments so that the first element of $S_w$ follows the last element of $S_v$ and the first element of $S_v$ follows the last element of $S_w$.} the segment excluding $e$ in $\sigma(v)$ and the segment excluding $e$ in $\sigma(w)$. Let $(G,\sigma)$ be a graph with a rotation system and $H$ be a minor of $G$. Let $\sigma'$ be the rotation system for $H$ obtained from $\sigma$ by repeated application of the above two operations in accordance with the relation $H\minorp G$. Observe that $\sigma'$ does not depend on the order in which we perform these operations. We say that $\sigma'$ is the rotation system \emph{induced} by $\sigma$ on $H$. We sometimes write $\sigma\restric_H$ for $\sigma'$.

\begin{lem}[folklore \cite{graphsonsurfaces2001}]
    Let $(G,\iota)$ be an embedded graph. Then for any minor $H$ of $G$, the graph $H$ with the induced rotation system is also an embedded graph. \qed
\end{lem}

Similar to temporal sequences of graphs, we define temporal sequences of embedded graphs referring to the \lq embedded minor relation\rq\ in place of the \lq minor relation\rq. 
Each temporal sequence of embedded graphs induces a temporal sequence of graphs by forgetting the embedding.
In this paper we are interested in the question which temporal sequences of graphs are induced by temporal sequences of embedded graphs in this way. We call such sequences of graphs \emph{simultaneously embeddable}; that is, a temporal sequence $(G_1,\ldots,G_n)$ of graphs is simultaneously embeddable if and only if for each $i\in [n]$ and each graph $G_i$, there exist plane embeddings $\iota_i$ for $G_i$ so that if $G_i$ is a minor of $G_{i+1}$, then $\iota_i$ is an embedded minor of $\iota_{i+1}$; and otherwise $\iota_{i+1}$ is an embedded minor of $\iota_{i}$. We call the sequence $(\iota_1,\ldots,\iota_n)$ a \emph{simultaneous embedding} of $(G_1,\ldots,G_n)$.
Let $(G,\sigma)$ be a graph with a rotation system. The \emph{reorientation} of $(G,\sigma)$ is the pair $(G,\iota^-)$ where $\iota^-(v)$ is the reverse of $\iota(v)$ for each vertex $v$ in $G$.
\begin{lem}[Whitney 1933 \cite{whitney2iso}]\label{lem:unique_3-con}
Every $3$-connected planar graph admits a unique embedding up to reorientation. \qed
\end{lem}

\begin{cor}
    Let $(G_1,\ldots,G_n)$ be a temporal sequence of planar graphs. If each graph $G_i$ is $3$-connected, then $(G_1,\ldots,G_n)$ admits a simultaneous embedding.
\end{cor}
\begin{pf}[Proof sketch:]
    Inductively pick embeddings of $G_{i+1}$, compatible with the embeddings for $G_i$. Here the key observation is that reorientation commutes with deletion and contraction for embedded graphs.
\end{pf}

The above corollary motivates the study of our problem via connectivity. In \autoref{sec:conc} we give a polynomial-time reduction from the problem $k$-Sunflower SEFE to the problem of deciding simultaneous embeddability. Since $k$-Sunflower SEFE is in general NP-hard \cite{gassner2006simultaneous}, and since we do not expect good characterisations for NP-hard problems, this suggests that we should restrict our attention to temporal sequences of $2$-connected graphs. 

Let $\Gcal=(G_1,\ldots,G_n)$ be a temporal sequence and let $\Hcal=(H_1,\ldots,H_m)$ be obtained from $\Gcal$ by (1) recursively deleting entries $G_i$ where either $G_{i-1}\minorp G_i\minorp G_{i+1}$ or $G_{i-1}\minorm G_i\minorm G_{i+1}$; and (2)~recursively deleting graphs at the start or end of the sequence which are a minor of their neighbour. The minor relations in $\Hcal$ start with $\minorm$ and end with $\minorp$ and alternate between $\minorm$ and $\minorp$. For the remainder of this paper, we will assume that all temporal sequences have this form and it is not a hard exercise to see that the simultaneous embeddings of $\Hcal$ are in bijection with those of $\Gcal$. We will often write $(G_1\minorm\ldots\minorp G_n)$ for a temporal sequence $(G_1,\ldots,G_n)$.

\begin{dfn}[$2$-component]
    Let $G$ be a graph and $t$ a $2$-separator of $G$. Consider the graph $G-t$. For each connected component $c$ of $G-t$, we call the subgraph obtained from $G$ by deleting all of the vertices not in $c$ or $t$, and all edges between the two vertices of $t$ a \emph{$2$-component} at $t$. For each edge $e$ spanning $t$ in $G$, we also call the subgraph of $G$ consisting of $e$ and its endvertices $t$ a \emph{$2$-component}.
\end{dfn}
Observe that the $2$-components at a given $2$-separator pairwise share exactly the separator vertices and the family of their edge sets forms a partition of the edge set of $G$. Given a $2$-connected graph $G$ and a $2$-separator $t$ of $G$, we say that $t$ is a \emph{tuplet} if it has at least three $2$-components.

A \emph{$2$-separation} of a $2$-connected graph $G$ is a bipartition $(V_1,V_2)$ of the vertex set of $G$ so that both $G[V_1]$ and $G[V_2]$ have at least two edges and $V_1\cap V_2$ consists of exactly two vertices $x$ and $y$. We call $\{x,y\}$ the \emph{separator} of $(V_1,V_2)$.

Given a graph $G$ with an edge set $A$, the graph \emph{spanned} by $A$ is the subgraph of $G$ consisting of the endvertices of edges in $A$ and the edges in $A$; we denote this graph by $G[[A]]$. 

\begin{figure}
    \begin{minipage}{0.39\textwidth}
        \hspace{10pt}
        {\Large\scalebox{0.55}{\begin{tikzpicture}
	\begin{pgfonlayer}{nodelayer}
		\node [style=none] (491) at (1, 0.25) {};
		\node [style=none] (492) at (1, -0.25) {};
		\node [style=none] (493) at (1, 0) {};
		\node [style=none] (494) at (2.5, 0) {};
		\node [style=none] (495) at (2.25, 0.25) {};
		\node [style=none] (496) at (2.25, -0.25) {};
		\node [style=black dot] (579) at (-9.25, -1.25) {};
		\node [style=black dot] (580) at (-9.25, 1.25) {};
		\node [style=black dot] (581) at (-6.75, 1.25) {};
		\node [style=black dot] (582) at (-6.75, -1.25) {};
		\node [style=black dot] (583) at (-4.75, 1.25) {};
		\node [style=black dot] (584) at (-4.75, -1.25) {};
		\node [style=black dot] (585) at (-3, 2.25) {};
		\node [style=black dot] (586) at (-1.25, 1.25) {};
		\node [style=black dot] (587) at (-3, -2.25) {};
		\node [style=black dot] (588) at (-1.25, -1.25) {};
		\node [style=small back] (589) at (-6.75, 0) {};
		\node [style=small back] (590) at (-4.75, 0) {};
		\node [style=none] (591) at (-6.75, 0) {$e$};
		\node [style=none] (592) at (-4.75, 0) {$e$};
		\node [style=black dot] (593) at (4.75, -1.25) {};
		\node [style=black dot] (594) at (4.75, 1.25) {};
		\node [style=black dot] (595) at (7.25, 1.25) {};
		\node [style=black dot] (596) at (7.25, -1.25) {};
		\node [style=black dot] (597) at (7.25, 1.25) {};
		\node [style=black dot] (598) at (7.25, -1.25) {};
		\node [style=black dot] (599) at (9, 2.25) {};
		\node [style=black dot] (600) at (10.75, 1.25) {};
		\node [style=black dot] (601) at (9, -2.25) {};
		\node [style=black dot] (602) at (10.75, -1.25) {};
	\end{pgfonlayer}
	\begin{pgfonlayer}{edgelayer}
		\draw [style=c0] (491.center) to (492.center);
		\draw [style=c0] (493.center) to (494.center);
		\draw [style=c0] (494.center) to (495.center);
		\draw [style=c0] (494.center) to (496.center);
		\draw [style=c0] (580) to (581);
		\draw [style=c0] (579) to (580);
		\draw [style=c0] (580) to (582);
		\draw [style=c0] (582) to (579);
		\draw [style=c0] (579) to (581);
		\draw [style=c0] (583) to (585);
		\draw [style=c0] (585) to (586);
		\draw [style=c0] (586) to (588);
		\draw [style=c0] (588) to (587);
		\draw [style=c0] (587) to (584);
		\draw [style=c1] (583) to (584);
		\draw [style=c1] (581) to (582);
		\draw [style=c0] (594) to (595);
		\draw [style=c0] (593) to (594);
		\draw [style=c0] (594) to (596);
		\draw [style=c0] (596) to (593);
		\draw [style=c0] (593) to (595);
		\draw [style=c0] (597) to (599);
		\draw [style=c0] (599) to (600);
		\draw [style=c0] (600) to (602);
		\draw [style=c0] (602) to (601);
		\draw [style=c0] (601) to (598);
	\end{pgfonlayer}
\end{tikzpicture}}}
    \end{minipage}
    \hfill%
    \begin{minipage}{0.59\textwidth}
        \caption{The $2$-sum of the two graphs on the left via the \bcol\ edge $e$ is the graph on the right.}
    \label{fig:2-sum}
    \end{minipage} 
\end{figure}
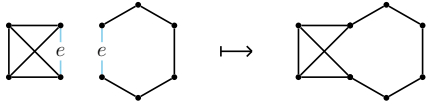

\begin{lem}[folklore]\label{2-sum-inverse}
Let $G$ be a $2$-connected graph with an edge-set bipartition $(A,B)$ that induces a $2$-separation. 
Let $e$ be an arbitrary edge of $B$. 
The graph $G$ has a minor $H$ with edge set $A+e$ that is equal to the graph obtained from $G[[A]]$ by adding an edge with label $e$ between the two vertices in the intersection of $G[[A]]$ and $G[[B]]$. 
\end{lem}
\begin{pf}
Pick an edge set $F$ that is a spanning forest of $G[[B]]$ so that the following holds: $F$ contains $e$; and in $F$ the edge $e$ separates the vertices of the separator of $(A,B)$. We let $H=G/(F-e)\sm (B\sm F)$. 
\end{pf}
We denote the graph $H$ of the above lemma by $G[[A:e]]$.
In a slight abuse of notation we also do that when we did not pick an edge $e$ in advance, in this case this simply means that one can pick $e\in B$ arbitrarily.
\begin{eg}
    Let $G$ be a $2$-connected graph with an edge-set bipartition $(A,B)$ that induces a $2$-separation. Then $G$ is the $2$-sum (see \autoref{fig:2-sum}) of the graphs $G[[A:e]]$ and $G[[B:f]]$ along their shared edge, after first identifying the labels $e$ and $f$.
\end{eg}

Basic examples of $2$-connected graphs are $3$-connected graphs and cycles. 
We will make extensive use of Tutte's Theorem (\ref{2SeparatorTheorem}), which canonically decomposes $2$-connected graphs into these basic pieces. Before stating it, we make the following definitions.

A~\emph{tree-decomposition} of a graph $G$ is a pair $(T,\Vcal)$ of a tree $T$ and a family 
$\Vcal=(V_t)_{t\in T}$ of vertex sets $V_t\subseteq V(G)$ indexed by the nodes $t$ of~$T$ such that the 
following conditions are satisfied:\vspace{-.25cm}
\begin{enumerate}
    \item[(T1)] $G=\bigcup_{t\in T}G[V_t]$; and \vspace{-.25cm}
    \item[(T2)] for every $v\in V(G)$, the vertex set $\{\,t\in T\mid v\in V_t\,\}$ is connected in~$T$.\vspace{-.25cm}
\end{enumerate}
The vertex sets $V_t$ and the subgraphs $G[V_t]$ they induce are the \emph{bags} of this
decomposition.
The \emph{torso} at a node $t\in T$ is obtained from the bag $G[V_t]$ by adding between any two vertices $x$ and $y$
of $V_t$ that are both contained in some bag  $V_{t'}$ with $t'\neq t$ an edge joining them, and deleting all other edges between $x$ and $y$.
The label of that edge is the label of an arbitrary edge in one of the subgraphs $H=\bigcup_{s\in S} G[V_s]$ of $G$, where $S$ is the subtree of $T-t$ that contain such a node $t'$ (whenever $H$ has at least one edge, and this will be satisfied in all torsos appearing in this paper). By repeated application of \autoref{2-sum-inverse}, we see that a torso is always a minor of $G$, provided that $G$ is $2$-connected.
Given a tree decomposition $(T,\Vcal)$ and an edge $e$ with endvertices $u$ and $v$, denote the component of $T-e$ containing $u$ by $R$ and the other by $S$. The pair $(\bigcup_{r\in R} V_r, \bigcup_{s\in S} V_s)$ is a separation of $G$; and we refer to it as the separation \emph{induced} by $(T,\Vcal)$ via the pair $(e,u)$. 

We say that a pair $(V_1,V_2),(U_1,U_2)$ of $2$-separations of a $2$-connected graph $G$ are \emph{crossed} if their separators are crossed. That is, one separator separates the vertices of the other. When $(V_1,V_2)$ and $(U_1,U_2)$ are not crossed, we say that they are \emph{nested}. Let us say that a $2$-separation of a graph is \emph{totally nested} if it is nested with every $2$-separation of the graph. Note that tuplets, $2$-separators with at least three $2$-components, are always nested with all other $2$-separators.

\begin{thm}[Tutte's $2$-Separation Theorem \cite{tutte1961theory} and \cite{triseparations}, simplified version]\label{2SeparatorTheorem}
For every $2$-connected graph~$G$, the totally-nested $2$-separa\-tions of~$G$ induce a tree-decomposition 
$(T,\Vcal)$ of $G$ all whose torsos are minors of $G$ and are $3$-connected, cycles, or have exactly two vertices.
Moreover, if a torso is $3$-connected, it is a maximal $3$-connected minor of $G$. Conversely, every maximal $3$-connected minor of~$G$ is a torso of  $(T,\Vcal)$. 
\end{thm}

We refer to the above canonical tree decomposition as the \emph{Tutte decomposition} of a given $2$-connected graph $G$ and the $3$-connected torsos of the Tutte decomposition as  the \emph{$3$-blocks} of $G$. 

Now in preparation for defining lists for temporal sequences, we make the following definitions. We say that two $3$-blocks of a graph $G$ are \emph{Tutte-equivalent} if they correspond to the same bag of the Tutte decomposition of $G$. We note that Tutte-equivalent $3$-blocks are the same up to a choice of labels for their torso edges. Now let $G$ be a planar graph and $b$ and $b'$ two Tutte-equivalent $3$-blocks. Since $b$ and $b'$ are $3$-connected, they admit unique embeddings up to reorientation. Let $\iota$ and $\iota'$ be embeddings of $b$ and $b'$, respectively, induced by the same embedding of $G$. Observe that any embedding of $G$ induces both $\iota$ and $\iota'$ or induces both of their reorientations. Hence we say that $\iota$ and $\iota'$ are \emph{Tutte-equivalent} embeddings of $b$ and $b'$.

Given a $2$-connected temporal sequence $\Gcal=(G_1\minorm\ldots\minorp G_n)$ of planar graphs, we refer to the (set of) \emph{$3$-blocks of} $\Gcal$ 
as the disjoint union of the $3$-blocks of all graphs of odd index $G_1,G_3,\ldots,G_n$.
A \emph{list} for $\Gcal$ is a partition of the set of $3$-blocks, together with a choice of an embedding for each $3$-block, so that each part is closed under Tutte-equivalence and Tutte-equivalent $3$-blocks receive Tutte-equivalent embeddings. For simplicity, we will often talk about lists as though they contain one representative from each Tutte-equivalence class of $3$-blocks and we note that storing one representative also suffices for computational implementations. This is justified by the notion of Tutte-equivalence of embeddings (see \autoref{rem:computing_with_lists}).
A \emph{simultaneous list-embedding} of $\Gcal$ for a list $L$ is a simultaneous embedding of $\Gcal$
so that if two $3$-blocks are in the same partition class of $L$, then their embeddings are both equal to the embeddings chosen in $L$ or both are the reorientations of their respective embeddings chosen in $L$. We define a simultaneous list-embedding of $\Gcal$ for a set of lists $\Lcal$ to be a simultaneous embedding of $\Gcal$ which satisfies all of the lists in $\Lcal$ individually. Two lists $L$ and $M$ are \emph{equivalent} if they both specify the same partition on the set of $3$-blocks and $M$ can be obtained from $L$ by reorienting all of the embeddings in some of the parts of $L$. It is immediate that a simultaneous embedding satisfies $L$ if and only if it satisfies an equivalent list~$M$.

\begin{dfn}[Trivial list]
    Let $\Gcal$ be a temporal sequence of planar graphs. Let $I$ be a list whose parts are exactly the Tutte-equivalence classes of the $3$-blocks of $\Gcal$ and whose specified embeddings are arbitrary. Then we call $I$ a \emph{trivial list} for $\Gcal$.
\end{dfn}

Let $\Gcal$ be a temporal sequence of planar graphs and let $I$ be a trivial list thereof. Then the simultaneous embeddings of $\Gcal$ are exactly the simultaneous list-embeddings of $(\Gcal,I)$. Also note that all the trivial lists of the same temporal sequence are equivalent. 

The following two definitions are also necessary for our main result.
\begin{dfn}[Tutte leaf]
    Let $G$ be a $2$-connected graph and $A$ a set of edges in $G$. Then we say that $A$ is a \emph{Tutte leaf} of $G$ if it is the edge set of a leaf bag of the Tutte decomposition of $G$ which corresponds to a $3$-block or a cyclic torso of $G$. We also say that $A$ is a \emph{Tutte leaf} of $G$ if it is a singleton consisting of an edge of $G$ which spans one of its $2$-separators.
\end{dfn}
A Tutte leaf is always the edge-set of a $2$-component at some $2$-separator of $G$.

\begin{dfn}[Shadow]
    Let $G$ be a graph and $A$ a set of edges in $G$. Let $H$ be a minor of~$G$. In the context of $G$ and $H$, we write $\pi(A)$ (or simply the juxtaposition $\pi A$) for the intersection $A\cap E(H)$. We refer to $\pi A$ as the \emph{shadow} of $A$ in $H$.
    In the context of $G$, we write $A^C$ for the complement of the edge set $A$ in $E(G)$. Since the sets $\pi(A)^C$ and $\pi(A^C)$ are the same (considering $\pi(A)^C$ as a complement in $E(H)$), we denote them simply by~$\pi A^C$.
\end{dfn}

\section{Main results}\label{sec:main}
We introduce the following measure of the complexity of a $2$-connected graph called `flexibility'. The general idea is that a less flexible graph admits fewer embeddings.

\begin{dfn}[Flexibility]\label{dfn:flex}
    Let $G$ be a $2$-connected graph with no vertices of degree $2$. Denote by $\alpha(G)$ the number of totally nested $2$-separators of $G$ with exactly two $2$-components and denote by $\beta(G)$ the number of tuplets $t$ of $G$, weighted by the number of $2$-components at $t$, take away $2$: in formulas: $\beta(G)=\Sigma_{\text{tuplets $t$}}(\Ccal(t)-2)$. Denote by $\gamma(G)$ the number of $2$-components at tuplets of $G$ which are single edges.
    Then the \emph{flexibility} of $G$ is the sum
$
        \Cbb(G):=1 + 2\cdot(\alpha(G)+3\beta(G))-\gamma(G).
$
    Now let $G$ be a $2$-connected graph, perhaps having some vertices of degree two. If $G$ is a cycle then the \emph{flexibility} $\Cbb(G)$ of $G$ is 1, otherwise it is defined as the flexibility of the graph obtained from $G$ after suppressing\footnote{That is, we replace all maximal paths whose internal vertices have degree two by single edges.} all of the vertices of degree two.
Given a temporal sequence $\Gcal=(G_1\minorm\ldots\minorp G_n)$, we define the \emph{flexibility} of $\Gcal$ to be the weighted sum $  \Cbb(\Gcal) = \Sigma_{i=1,3,5,\ldots,n}(n-i+1)\cdot\Cbb(G_i)$.

That is, we count the flexibility of the `big' (odd-index) graphs in the sequence and give more weight to those which are further from the end (the intuition behind this is that we will work our way through the sequence from left to right, so \lq prioritise\rq\ improvements further to the left).
\end{dfn}

The \emph{size} of the temporal sequence $\Gcal=(G_1\minorm\ldots\minorp G_n)$ is another measure of complexity and it is defined simply as\vspace{-.45cm}
\begin{equation*}
    \text{size}(\Gcal) = \Sigma_{i=1,3,5,\ldots,n} |E(G_i)|.\vspace{-.1cm}
\end{equation*}

\begin{dfn}[Analogue]
    Let $\Gcal=(G_1\minorm\ldots,\minorp G_n)$ and $\Hcal=(H_1\minorm\ldots\minorp H_m)$ be temporal sequences along with sets of lists $\Lcal$ and $\Mcal$ of $\Gcal$ and $\Hcal$, respectively. We say that the pair $(\Hcal,\Mcal)$ is an \emph{analogue} of the pair $(\Gcal,\Lcal)$ if $(\Hcal,\Mcal)$ has a simultaneous list-embedding if and only if $(\Gcal,\Lcal)$ has a simultaneous list-embedding. In this context, if $\Lcal=\{L\}$, we shall write $(G,L)$ as an abbreviation for $(G,\{L\})$ in an attempt to increase readability. 
\end{dfn}

\begin{dfn}[Nice]
    We say that a graph $G$ is \emph{nice} if it is $2$-connected and both:
    \begin{enumerate}
        \item the distance between every pair of tuplets of $G$ is at least two in $G$; and
        \vspace{-.25cm}
        \item each vertex $x$ in a tuplet $t$ of $G$ is incident with exactly one edge from each $2$-component at~$t$.
    \end{enumerate}
    \vspace{-.28cm}We say that a temporal sequence is \emph{nice} if all of its graphs of odd index are nice.
\end{dfn}

\begin{dfn}[Improvement]\label{dfn:improv}
    Let $\Gcal=(G_1\minorm\ldots,\minorp G_n)$ be a $2$-connected temporal sequence and $\{L,M\}$ a pair\footnote{We believe that the formulation with two lists (instead of some other number) gives us the cleanest result. Our improvement procedure will first determine if these two lists $L$ and $M$ can be combined into some list $L+M$ (if not, we get an `uncombinable obstruction'). The main part of the procedure generates one new list $N$, which may or may not be combinable with $L+M$. Hence, we are left with two lists in the output.} of lists thereof. An \emph{improvement} of $(\Gcal,\{L,M\})$ is an analogue $(\Gcal',\{L',M'\})$ so that $\Gcal'$ is nice and\vspace{-.4cm}
    \[
        \text{size}(\Gcal')\leq 100\cdot(\text{size}(\Gcal))^2. \vspace{-.1cm}
    \]
    Furthermore if $\Gcal$ is nice, then additionally $\Cbb(\Gcal')\leq\Cbb(\Gcal)-1$ and 
$\text{size}(\Gcal')\leq \text{size}(\Gcal) + 13\cdot(\Cbb(\Gcal))^2.$
\end{dfn}

The main result of this paper is the following.
\begin{thm}\label{thm:main}
    Every nonempty $2$-connected temporal sequence with two lists is obstructed or admits an improvement.
\end{thm}

\begin{thm}[Algorithmic strengthening of \autoref{thm:main}]\label{thm:algorithmic_main}
    Let $\Gcal$ be a nonempty temporal sequence of $2$-connected graphs and $\{L,M\}$ a pair of lists thereof. There is a polynomial-time algorithm which decides whether $(\Gcal,\{L,M\})$ is obstructed or returns an improvement of $(\Gcal,\{L,M\})$.
\end{thm}

\begin{cor}\label{cor:poly}
    There exists a polynomial-time algorithm that decides simultaneous embeddability of $2$-connected temporal sequences.
\end{cor}
\begin{pf}
Given a nonempty $2$-connected temporal sequence $\Gcal_0=(G_1,..,G_n)$, check whether all graphs $G_i$ are planar, which takes linear time using Hopcroft and Tarjan \cite{HopcroftTarjan1974}. Since we are done otherwise, assume that all $G_i$ are planar. 
Hence pick two trivial lists for $\Gcal_0$ and denote them by $L_0$ and $M_0$. 
Clearly, $\Gcal_0$ admits a simultaneous embedding if and only if $(\Gcal_0,\{L_0,M_0\})$ admits a simultaneous list-embedding. If $(\Gcal_0,\{L_0,M_0\})$ is obstructed, we are done, so assume that it is not. 
    By applying \autoref{thm:algorithmic_main} to $(\Gcal_0,\{L_0,M_0\})$ and replacing $(\Gcal_0,\{L_0,M_0\})$ by the output of \autoref{thm:algorithmic_main} if necessary, we assume that $\Gcal_0$ is nice. 
    We apply \autoref{thm:algorithmic_main} iteratively to obtain a sequence $(\Gcal_0,\{M_0,L_0\}),\ldots,(\Gcal_n,\{M_n,L_n\})$ of analogues of $(\Gcal_0,\{L_0,M_0\})$ terminating when either $(\Gcal_n,\{L_n,M_n\})$ is obstructed (in which case $\Gcal_0$ is not simultaneously embeddable by the definition of analogue) or $\Gcal_n$ is an empty sequence (in which case $\Gcal_0$ admits a simultaneous embedding). In each step $\Cbb(\Gcal_{i+1})\leq\Cbb(\Gcal_i)-1$ and so the process terminates after at most $\Cbb(\Gcal_0)$ steps; note that this term is polynomial in the size of $\Gcal_0$.
    
    In each step our input $\Gcal_i$ is nice and
$
        \text{size}(\Gcal_{i+1})\leq \text{size}(\Gcal_0) + (i+1)\cdot13\cdot(\Cbb(\Gcal_0))^2
    $. Thus each of the polynomially many steps runs in polynomial time, leading to polynomial running time overall. 
\end{pf}

\begin{figure}
    \begin{minipage}{0.54\textwidth}
        \centering
        \Large
        \scalebox{0.52}{\input{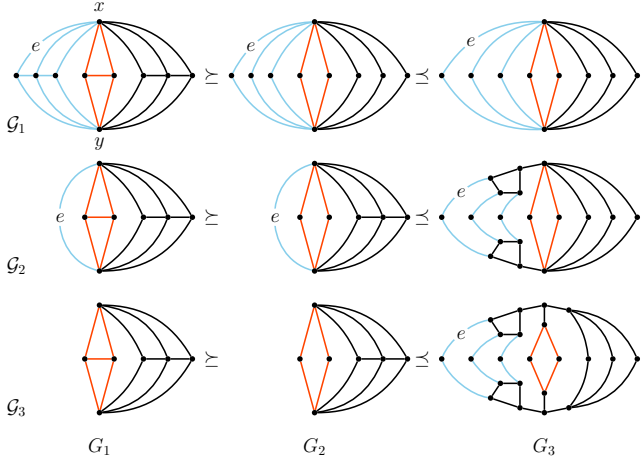}}
    \end{minipage}
    \hfill%
    \begin{minipage}{0.44\textwidth}
        \caption{Three temporal sequences $\Gcal_1$, $\Gcal_2$ and $\Gcal_3$. The temporal sequence $\Gcal_2$ is obtained from $\Gcal_1$ by reducing the \bcol\ $2$-component at $\{x,y\}$ in $G_1$ to a single edge and `uncontracting' some edges in $G_3$. The temporal sequence $\Gcal_3$ is obtained from $\Gcal_2$ by deleting the edge $e$ from $G_1$ and $G_2$ and `uncontracting' some edges in $G_3$. Each temporal sequence has the same number of simultaneous embeddings.}
    \label{fig:main-strategy}
    \end{minipage}
    \vspace{-.5cm}
\end{figure}

{\bf Overall structure of the proof of \autoref{thm:main}.} If a given $2$-connected temporal sequence is $3$-connected, then we are done. Otherwise we find a \lq suitable\rq\ Tutte leaf (an edge set of a leaf of the Tutte decomposition) of some graph in our temporal sequence, and then we try to contract it to a single edge (\nameref{lem:A}) or delete it if it is already a single edge (\nameref{lem:B_2}). This way the Tutte leaf is eliminated in a sequence of two steps, see \autoref{fig:main-strategy} for an example with a concrete sequence. To do this type of elimination in the most general context, we will sometimes need to `add information' to the graphs further on in the sequence in order to \lq pay\rq\ for this elimination. This makes the later graphs slightly larger and adds lists of additional information; to analyse how this works in detail is a major technical part of our proof. Examples \ref{eg:expansion_impossible} and \ref{eg:model_impossible} below demonstrate that the additional information cannot be displayed using graphs alone but one needs more abstract information---that is, lists. In the following we shall introduce notation to make all of this precise.

Given a temporal sequence $\Gcal=(G_1\minorm\ldots\minorp G_n)$, an \emph{upleaf} of $\Gcal$ is a Tutte leaf $A\subseteq E(G_i)$ of some graph $G_i$ with $i$ odd, so that either $i=1$ or the graph $G_{i-1}$ does not contain any edge of $A$ (that is, $A$ has empty shadow in $G_{i-1}$). In this case, we say that $A$ is \emph{rooted} at $i\in[n]$. We say an upleaf $A$ in $G_i$ is a \emph{last upleaf} if the index $i$ is maximal amongst all upleaves across the whole temporal sequence $\Gcal$.

\begin{lem}\label{unamazing_lemma}
 Let $\Gcal=(G_1\minorm\ldots\minorp G_n)$ be a nice $2$-connected temporal sequence with a list $L$, and suppose that $A$ is a last upleaf of $\Gcal$ rooted at $i$ for some $i\in[n]$ with $|A|\geq 2$ and $A\neq E(G_i)$. If $i=n$ or $A$ has empty shadow in $G_{i+1}$, then $(\Gcal,L)$ admits an improvement.
\end{lem}

\begin{pf}
Let $\Gcal':=(G_1\minorm\ldots\minorm G_{i-1}\minorp G_i[[A\ct:e]]\minorm G_{i+1}\minorp\ldots\minorp G_n)$. We obtain $L'$ from $L$ by deleting the torso corresponding to $A$ (so $L'=L$ if the torso is not a $3$-block). 
Since $\Cbb(G_i)>\Cbb(G_i[[A^C:e]])$, $\text{size}(\Gcal')\leq \text{size}(\Gcal)$ and $\Gcal'$ is nice, it suffices to show that $(\Gcal',L')$ is an analogue of $(\Gcal,L)$.
Since $\Gcal'$ is obtained by taking minors in $\Gcal$, every simultaneous list-embedding of $\Gcal$ for $L$ induces a
simultaneous list-embedding of $\Gcal'$ for $L'$.
For the converse, let a simultaneous list-embedding $(\phi_1,\ldots,\phi_n)$ of $(\Gcal',L')$ be given. 
Let $\iota$ be an embedding of the torso for the Tutte leaf $A$ that respects the list $L$. 
We combine the embedding $\phi_i$ of $G_i[[A:e]]$ with $\iota$ to obtain an embedding $\psi_i$ of their 2-sum along the edge $e$ via \autoref{lem:2-sum_embeddings} below.
By \autoref{2-sum-inverse}  this 2-sum is the graph $G_i$. Then $(\phi_1,\ldots,\phi_{i-1}, \psi_i, \phi_{i+1},\ldots, \phi_n)$ is a simultaneous list-embedding for $(\Gcal,L)$.
Hence $(\Gcal',L')$ is an analogue of $(\Gcal,L)$, and thus an improvement. 
\end{pf}

Next we present a version of \autoref{unamazing_lemma} for when $A$ has nonempty shadow in $G_{i+1}$.

\begin{dfn}[Amazing analogue]
Let $\Gcal=(G_1\minorm\ldots\minorp G_n)$ be a nice $2$-connected temporal sequence with a list $L$. Let $A$
 be a last upleaf of $\Gcal$ rooted at $i$ for some $i\in[n]$ so that some edge of $A$ appears in $G_{i+1}$. An \emph{amazing analogue} of $(\Gcal,L)$ for the pair $(A,i)$ is a temporal sequence \vspace{-.25cm}
    \[
    \Gcal':=(G_1\minorm\ldots\minorm G_{i-1}\minorp G_i[[A\ct:e]]\minorm G_{i+1}[[\pi A\ct:e]]\minorp H_{i+2}\minorm G_{i+3}\minorp\ldots\minorp G_n) \vspace{-.25cm}
    \]
    along with a list $L'$ so that $H_{i+2}\minorm G_{i+2}$ is nice; and $(\Gcal',L')$ is an analogue for $(\Gcal,L)$ and furthermore \vspace{-.6cm}
    \[
    \Cbb(H_{i+2})\leq \Cbb(G_{i+2})\quad\text{and}\quad |E(H_{i+2})|\leq |E(G_{i+2})| + 4\;\Cbb(G_{i+2}).
    \vspace{-.05cm}
    \]
\end{dfn}

\begin{lem}[Amazing Analogue Lemma]\label{lem:A}
    Let $\Gcal=(G_1\minorm\ldots\minorp G_n)$ be a nice $2$-connected temporal sequence with a list $L$, and suppose that $A$ is a last upleaf of $\Gcal$ rooted at $i$ for some $i\in[n]$. Suppose further that some edge of $A$ appears in $G_{i+1}$. If $\Gcal$ is not a disagreeable obstruction or an amazingly-bad obstruction and the pair $(\Gcal,L)$ is not an inconsistent obstruction, then $(\Gcal,L)$ admits an amazing analogue for the pair $(A,i)$.
\end{lem}

\begin{dfn}[Delightful analogue]
Given a temporal sequence $\Gcal=(G_1\minorm\ldots\minorp G_n)$ with a list $L$ and an edge $d$ of a graph $G_i$, a \emph{delightful analogue} of $(\Gcal,L)$ for the pair $(d,i)$ is a temporal sequence \vspace{-.55cm}
    \[
    \Gcal':=(G_1\minorm\ldots\minorm G_{i-1}\minorp G_i-d\minorm H_{i+1}\minorp H_{i+2}\minorm G_{i+3}\minorp\ldots\minorp G_n) \vspace{-.15cm}
    \]
    along with a pair of lists $M$ and $N$ so that $H_{i+1}\minorm G_{i+1}-d$ and $H_{i+2}\minorm G_{i+2}$ are nice, we have $(\Gcal',\{M,N\})$ is an analogue of $(\Gcal,L)$ and furthermore  $
        |E(H_{i+1})|\leq |E(G_{i+1})|+1,$ \vspace{-.25cm}
    \[
    \Cbb(H_{i+2})\leq \Cbb(G_{i+2})+5\quad\text{and}\quad |E(H_{i+2})|\leq |E(G_{i+2})| + 12\,(\Cbb(G_{i+2}))^2+1.
    \vspace{-.25cm}
    \]
\end{dfn}

\begin{lem}[Delightful Analogue Lemma]\label{lem:B_2}
Let $\Gcal=(G_1\minorm\ldots\minorp G_n)$ be a nice $2$-connected temporal sequence  with a list $L$,
and $A$ a last upleaf of $\Gcal$ consisting of a single edge and rooted at $i\in [n]$. If $\Gcal$ contains no gaudy obstruction, then $(\Gcal,L)$ admits a delightful analogue for $(d,i)$.
\end{lem}

Finally, we require the following lemmas, proved in \autoref{subsec:uncombinable} and \autoref{sec:nice}, respectively.
\begin{lem}\label{lem:list-combining_machine}
    Let $L$ and $M$ be a pair of lists for a temporal sequence $\Gcal$ and suppose that together they do not form an uncombinable obstruction. There is a list $N$ for $\Gcal$ so that $(\Gcal,N)$ is an analogue of $(\Gcal,\{L,M\})$.
\end{lem}
Computing with lists can be done efficiently, as outlined in \autoref{rem:computing_with_lists}.

\begin{lem}\label{lem:niceness_main}
For every $2$-connected temporal sequence $\Gcal$ with list $L$, there exists an analogue $(\Hcal,L')$ of $(\Gcal,L)$ so that $\Hcal$ is nice.
    Furthermore, $\text{size}(\Hcal)\leq 100\cdot(\text{size}(\Gcal))^2$.
\end{lem}

\begin{pf}[Proof of \autoref{thm:main}, given \ref{lem:A}, \ref{lem:B_2}, \ref{lem:list-combining_machine} and \ref{lem:niceness_main}.]
Let a nonempty temporal sequence\break $(\Gcal,\{L,M\})$ of $2$-connected graphs with lists be given and suppose that it is not obstructed. We are to construct an improvement. 
    Since $(\Gcal,\{L,M\})$ contains no uncombinable obstruction, by \autoref{lem:list-combining_machine} there is a list $N$ for $\Gcal$ so that $(\Gcal,N)$ is an analogue of $(\Gcal,\{L,M\})$. For the remainder of this proof, we will work with $(\Gcal,N)$. Now if $\Gcal$ is not nice, we obtain an improvement via  \autoref{lem:niceness_main}. So assume that $\Gcal=(G_1\minorm\ldots\minorp G_n)$ is nice. Recall that for nice $\Gcal$ we are to construct an improvement that is not only bounded in size but also in flexibility.
Since all of the Tutte leaves of $G_1$, of which there is at least one, are upleaves rooted at $1$, the temporal sequence $\Gcal$  has an upleaf. Now let $A$ be a last upleaf, say rooted at $i$ for some odd $i\in[n]$.

\begin{claim}\label{claimX}
If $E(G_i)=A$, then $(\Gcal,N)$ has an improvement.
\end{claim}

\begin{cproof}
Since $E(G_i)=A$, it has empty shadow in $G_{i-1}$ (assuming $i\geq 2$), so $G_{i-1}$ would be empty and contradict the assumption of $2$-connectivity. Thus $i=1$. We distinguish two cases. 

{\bf Case 1:} $G_1$ is a cycle. Since $G_1$ is a cycle, its 2-connected minor $G_2$ is also a cycle.
Hence $((G_3\minorm\ldots\minorp G_n), N)$ is an analogue of $(\Gcal,N)$.
Since it is nice, has reduced size and flexibility, it is an improvement. 

{\bf Case 2:} not Case 1. Since $A=E(G_1)$ is a Tutte leaf, it must be a 3-block. 
If $G_2$ is a 3-block or a cycle, then we argue as in Case 1, so assume that $\Cbb(G_2)>1$.
By 2-connectedness of $G_2$, it contains an edge $e$. 
Since $G_2\minorp G_1$, the edge $e$ is an edge of $G_1$. 
We obtain $G_1'$ from $G_1$ by subdividing the edge $e$; denote the subdivision edges by $e'$ and $e''$. 
Thus $G_1'$ is the 2-sum of $G_1$ with the triangle $\{e,e',e''\}$ along the edge $e$. 
Let $\Gcal'=(G_1'\minorm G_2 \minorp\ldots\minorp G_n)$.
Now $(\Gcal',N)$ is a nice analogue of $(\Gcal,N)$, whose size increased only by one. Applying \myref{lem:A} to $(\Gcal',N)$ yields an amazing analogue of $(\Gcal',N)$ whose first and second graphs are cycles (we contract the whole of $G_1$ onto a single edge $e$), so have flexibility one; denote it by $(\Gcal'',N'')$. 
Since $\Cbb(G_2)>1$, the temporal sequence $\Gcal''$ has smaller flexibility than $\Gcal$. Since $\Gcal''$ is nice and its size is appropriately bounded, $(\Gcal'',N'')$ is an improvement of $(\Gcal,N)$.
\end{cproof}

If $A$ contains exactly one edge, then by \myref{lem:B_2} applied to that edge as $d$
yields an improvement for $(\Gcal,N)$, since $\Cbb(G_i-d)\leq\Cbb(G_i)-2\cdot3+1=\Cbb(G_i)-5$, $G_i-d$ is nice and the size is bounded appropriately. 
So assume that $A$ corresponds to a $3$-block or a cyclic torso of $G_i$. 
By \autoref{claimX} assume that $E(G_i)\neq A$. 
Hence $\Cbb(G_i[A\ct:e])\leq\Cbb(G_i)-1$.
If $i=n$ or $A$ has empty shadow in $G_{i+1}$, apply \autoref{unamazing_lemma} and otherwise apply \myref{lem:A}, yielding the desired improvement of $(\Gcal,N)$ in either case. 
\end{pf}

\section{Obstructions}\label{sec:obs}
In this section, we detail the five obstructions to simultaneous embeddability arising in \autoref{thm:main}. Given a temporal sequence $\Gcal=(G_1\minorm\ldots\minorp G_n)$ one of whose graphs $G_i$ is non-planar, we refer to $G_i$ as a \emph{nonplanar obstruction}. 
Recall that a list of $\Gcal$ includes plane embeddings of all 3-blocks. Thus in the presence of a list, all torsos of the Tutte decomposition of the graphs $G_i$ are planar. And since 2-sums of planar graphs are planar, all graphs $G_i$ are planar. Hence the assumption in \autoref{thm:main} that we are given a temporal sequences with a list excludes these trivial obstructions. However, nonplanar obstructions arose implicitly already in the proof of \autoref{cor:poly}.
To summarise, for technical reasons, nonplanar obstructions do not arise in the context of \autoref{thm:main}.

\subsection{Disagreeable obstructions}
Given two temporal sequences $(G_1,\ldots,G_n)$ and $(H_1,\ldots,H_n)$ of the same length $n$,
we say that $(H_1,\ldots,H_n)$ is a \emph{(temporal) minor} of $(G_1,\ldots,G_n)$ if $H_i\minorp G_i$ for all $i\in [n]$. This definition naturally generalises to sequences of not necessarily the same length, yet the definition is a bit of a mouthful (this is why we have given the definition in the special case first). 
Given two temporal sequences $(G_1,\ldots,G_m)$ and $(H_1,\ldots,H_n)$, we say that $(H_1,\ldots,H_n)$ is a \emph{(temporal) minor} of $(G_1,\ldots,G_m)$ if $m\geq n$, and there is a sequence $a_1< a_2<  \ldots< a_{n} < a_{n+1}$ so that for all $i\in[n]$ and all $j\in [a_i,a_{i+1}-1]$ we have that $H_i\minorp G_j$. 

Given a temporal triple $\tseq{G}$ and embeddings $\phi_i$ of $G_i$ for $i=1,3$, we write
{$\phi_1\rightarrow\nolinebreak G_2\leftarrow \phi_3$} as a shorthand for the 
statement that both $\phi_1$ and $\phi_3$ induce the same embedding on $G_2$. We denote a graph consisting of three edges in parallel by $\triple$ 
(we refer to this graph as a \emph{theta-graph}). 
Let $H_1$ be a $3$-connected minor of $G_1$ and let $\alpha$ be one of the embeddings of $H_1$. We say that an embedding $\psi_3$ of $G_3$ is \emph{$\alpha$-greeable}\footnote{Pronounced `agreeable'.} if for every $3$-block $x$ of $G_3$ and every temporal minor {$(H_1\minorm\triple\minorp x)$} of $(G_1 \minorm G_2 \minorp G_3)$ we have that $\alpha \rightarrow \triple \leftarrow \psi_3\restric_x$; that is, the embeddings $\alpha$ and $\psi_3\restric x$ induce the same embedding on~$\triple$.

\begin{dfn}[Disagreeable obstruction]\label{obs:disagreeable}
    A \emph{disagreeable obstruction} is a temporal triple\\ $\tseq{G}$ along with a $3$-connected minor $H_1$ of $G_1$ and a $3$-connected minor $H_3$ of $G_3$ so that the following holds. There exists a pair of temporal minors $(H_1\minorm\triple_1\minorp H_3)$ and $(H_1\minorm\triple_2\minorp H_3)$ of $\tseq{G}$, with $\triple_1$ and $\triple_2$ theta graphs, so that no choice of embeddings $\alpha_1$ for $H_1$ and $\alpha_3$ for $H_3$ satisfy both $\alpha_1\rightarrow\triple_1\leftarrow\alpha_3$ and $\alpha_1\rightarrow\triple_2\leftarrow\alpha_3$.
\end{dfn}

\begin{lem}\label{lem:polynomial_disagreeable}
There is a polynomial algorithm that checks whether a given temporal triple is a disagreeable obstruction.
\end{lem}
We will see a proof of \autoref{lem:polynomial_disagreeable} in \autoref{appendix:algorithms}.

Given a temporal triple $\tseq{G}$, a minor $H_1$ of $G_1$ is \emph{linked to} a minor $H_3$ of $G_3$ if there is a temporal triple of the form $(H_1\minorm\triple\minorp H_3)$ that is a minor of $\tseq{G}$. 
Given a temporal triple $\tseq{G}$, a minor $H_1$ of $G_1$ is \emph{coupled with} a minor $H_3$ of $G_3$ if for every embedding $\alpha_1$ of $H_1$ there is a unique embedding $\alpha_3$ of $H_3$ so that 
 for every  temporal triple of the form $(H_1\minorm\triple\minorp H_3)$ that is a minor of $\tseq{G}$ we have that $\iota=\alpha_3$ is the only embedding of $H_3$ that satisfies
$\alpha_1\rightarrow\triple\leftarrow\iota$.

\begin{lem}\label{not-disagree}
Assume that $\tseq{G}$ admits a simultaneous embedding. Let $H_1$ and $H_3$ be $3$-connected minors of $G_1$ and $G_3$, respectively. Then $H_1$ is linked to $H_3$ if and only if $H_1$ is coupled with $H_3$. 
\end{lem}
\begin{pf}
Suppose first that $H_1$ is coupled to $H_3$. Suppose for a contradiction that there is no minor of $\tseq{G}$ of the form $(H_1\minorm\triple\minorp H_3)$. Then the last condition of the definition of \lq coupled with\rq\ is satisfied by every pair $(\alpha_1,\alpha_3)$ of embeddings $\alpha_1$ and $\alpha_3$ of $H_1$ and $H_3$, respectively. For each $i=1,3$, since $H_i$ is $3$-connected it has two embeddings which are non-identical. This gives us our contradiction and we conclude that $H_1$ is linked to $H_3$. For the converse implication, assume that $H_1$ is linked to $H_3$.

By assumption there is a simultaneous embedding of $\tseq{G}$; denote it by $(\phi_1,\phi_2,\phi_3)$. Let $\alpha_1$ be one of the two embeddings of $H_1$. By reorienting all of the embeddings of $(\phi_1,\phi_2,\phi_3)$ if necessary, we assume that $\alpha_1$ is the embedding induced on $H_1$ by $\phi_1$. Let $\alpha_3$ denote the embedding induced on $H_3$ by $\phi_3$.
Now since $H_1$ is linked with $H_3$ there is a temporal triple of the form $(H_1\minorm\triple\minorp H_3)$ that is a minor of $\tseq{G}$.
The embeddings $\alpha_1$ and $\alpha_3$ satisfy $\alpha_1\rightarrow\triple\leftarrow\alpha_3$ since they are each induced by the same simultaneous embedding of $\tseq{G}$.
Now since only one of $\alpha_1\rightarrow\triple\leftarrow\alpha_3$ and $\alpha_1\rightarrow\triple\leftarrow(\alpha_3)^-$ is true, and the only embeddings of $H_3$ are $\alpha_3$ and $(\alpha_3)^-$, we deduce that $\iota=\alpha_3$ is the only embedding of $H_3$ that satisfies
$\alpha_1\rightarrow\triple\leftarrow\iota$. In fact by the same analysis, the embedding $\iota=\alpha_3$ satisfies $\alpha_1\rightarrow\triple\leftarrow\iota$ for any choice of theta graph $\triple$ satisfying $H_1\minorm\triple\minorp H_3$. This completes our proof.
\end{pf}

\begin{cor}\label{cor:disagreeable_bad}
    No disagreeable obstruction admits a simultaneous embedding.
\end{cor}
\begin{pf}
Assume that $\tseq{G}$ along with a $3$-connected minor $H_1$ of $G_1$ and a $3$-connected minor $H_3$ of $G_3$ forms a disagreeable obstruction. Then $H_1$ is linked to $H_3$. Yet $H_1$ is not coupled with $H_3$. Hence the conclusion of \autoref{not-disagree} is violated. Thus one of the assumptions of that lemma does not hold. And the only possibility is that $\tseq{G}$ does not admit a simultaneous embedding. 
 \end{pf}

\subsection{Amazingly-bad obstructions}

We will need the following machinery.

\begin{lem}[\cite{whitney2iso}]\label{lem:comp_segs}
    Let $G$ be a $2$-connected graph, $t$ a $2$-separator of $G$ and $x$ a vertex thereof. Let $\iota$ be an embedding of $G$. Then the edges incident with $x$ in each $2$-component at $t$ form a segment in the rotator $\iota(x)$ at $x$. 

    Furthermore, the cyclic order of these segments at $x$ is the reverse of the cyclic order of these segments at the rotator from $\iota$ at the other vertex $y\in t\sm\{x\}$.
\end{lem}
\begin{pf}
    Suppose for a contradiction that the edges from some $2$-component $c$ do not form a segment in the rotator $\iota(x)$. Then there are two edges $a,b\in\iota(x)$, not in $c$, so that both intervals between them in $\iota(x)$ contain edges from $c$. We can find two internally disjoint paths from $x$ to the vertex $y\in t\sm\{x\}$, starting from $a$ and $b$, respectively. These form a cycle in~$G$, not using edges of $c$. The cycle topologically separates edges of $c$ in $\iota$ and meets $c$ exactly in $x$ and $y$. Hence $t$ is a separator of $c$, contradicting the fact that it is a single $2$-component. Hence, the edges incident with $x$ in each $2$-component at $t$ form a segment in the rotator $\iota(x)$ at $x$. 

    Let $H$ be the graph obtained by contracting each $2$-component $c$ at $t$ to a single edge labeled by $c$ that spans $x$ and $y$. Then clearly the rotator $\iota\restric_H(y)$ is the reversal of $\iota\restric_H(x)$ and this fact lifts to the cyclic order on the components at $t$ in $G$.
\end{pf}

Let $G$ be a $2$-connected graph and $t$ a $2$-separator thereof, and $x$ a vertex of $t$.
Let $\iota$ be an embedding of $G$. Consider the rotator $\iota(x)$ of $\iota$ at $x$.
By \autoref{lem:comp_segs}, for each 2-component at $t$, the edges of that 2-component incident with $x$ form a segment in the rotator $\iota(x)$. These segments partition the set of edges incident with $\iota(x)$. 
So $\iota(x)$ induces a cyclic order on the set of these segments (simply starting with $\iota(x)$ and contracting each segment to a point). We refer to this cyclic order as the \emph{cyclic order at the pair $(t,x)$}.

Suppose we are given a $2$-connected graph $G$ and an embedding $\iota$ of $G$, from $\iota$ we can read off an orientation of each $3$-block (indeed, it is a minor of $G$ so $\iota$ induces an embedding on it) and a  cyclic order of $2$-components at all pairs $(t,x)$, where $t$ is a tuplet and $x\in t$. We refer to this as the \emph{embedding data} of $\iota$.

\begin{lem}[\cite{whitney2iso}]\label{determine_embed}
Embeddings of $2$-connected graphs with the same embedding data are equal. 
\end{lem}

The following completes our toolset for deconstructing embeddings of $2$-connected graphs.
\begin{lem}\label{lem:2-sum_embeddings}
    Let $G$ and $H$ be $2$-connected graphs that share exactly one edge, which we denote by~$e$. Then for every pair of embeddings $\iota$ for $G$ and $\kappa$ for $H$, there exists a unique embedding of the $2$-sum $G\oplus_e H$ that induces $\iota$ and $\kappa$ on the minors $G$ and $H$, respectively.
\end{lem}
\begin{pf}
    Let $\iota$ and $\kappa$ be embeddings of $G$ and $H$, respectively. Define $\phi$ to be the embedding of $G\oplus_e H$ obtained in the following way. In both 2-spheres of the embeddings $\iota$ and $\kappa$, cut out a small disc around the edge $e$ that intersects $e$ only in its endvertices and only has $e$ on its inside. Now glue the two remaining 2-spheres with boundary at the boundary in the orientation-preserving way so that the endvertices of $e$ are identified from both copies in the correct way round. This gives a 2-sphere together with an embedding of the graph $G\oplus_e H$ into it.

    Now let $\phi'$ be an embedding of $G\oplus_e H$. The graphs $G$ and $H$ are minors of $G\oplus_e H$ (either $G$ or $H$ is a cycle or we apply \autoref{2-sum-inverse}). Let $\iota'$ be the embedding induced on $G$ by $\phi'$ and $\kappa'$ be the embedding induced on $H$ by $\phi'$. Observe that if $\phi'$ was constructed from embeddings $\iota$ and $\kappa$ as in the first part of the proof, that is $\phi'=\phi$, then $\iota'=\iota$ and $\kappa'=\kappa$. For the uniqueness of our construction, observe conversely that if $\iota'=\iota$ and $\kappa'=\kappa$, then we have $\phi=\phi'$.
\end{pf}

Let $G$ be a $2$-connected graph and $x$ a vertex in a tuplet $t$ of $G$. Let $\iota$ be an embedding of $G$ and recall the definition of `cyclic order at the pair $(t,x)$' from above.
By the \lq Furthermore\rq-part of \autoref{lem:comp_segs}, if $y$ is the vertex of $t$ other than $x$, then the cyclic order at $(t,y)$ is the reverse of the cyclic order at $(t,x)$. 
We call the unordered pair, consisting of this cyclic order and its reverse, the \emph{cyclic structure\footnote{In this paper, we use the word `cyclic structure' to refer to the unordered pair of a cyclic order and its reverse. This object describes a cyclic order up to reversal.} induced on the $2$-components at $t$}.

Let $G$ be a graph and $t$ a $2$-separator thereof. Let $A$ be a set of edges, perhaps not entirely contained in the edge set of $G$. We say that a $2$-component at $t$ is \emph{$A$-ffected} if it contains an edge of $A$.

\begin{lem}\label{lem:affected_G_2}
    Let $(G_1\minorm G_2)$ be a $2$-connected temporal pair. Let $A$ be a $3$-connected Tutte leaf of $G_1$ and let $t$ be a tuplet in $G_2$. Suppose that not every $2$-component at $t$ is $A$-ffected. Then in every simultaneous embedding $(\phi_1,\phi_2)$ of $(G_1\minorm G_2)$, the $A$-ffected $2$-components at $t$ form a segment in the cyclic structure induced on the $2$-components at $t$ by $\phi_2$, and furthermore the linear structure of this segment is independent of our choice of simultaneous embedding.
\end{lem}
\begin{pf}
    Since otherwise there would be nothing to prove, we assume that the shadow $\pi A$  of $A$ in $G_2$ contains at least one edge. By the assumption that not every $2$-component at $t$ is $A$-ffected, we also get that the shadow $\pi A\ct$ of $A\ct$ in $G_2$ contains at least one edge. Hence $(V(\pi A),V(\pi A\ct))$ is a separation of $G_2$ of adhesion two (but perhaps not a proper $2$-separation). We consider the minor $G_2[[\pi A: e]]$ of $G_2$ given by \autoref{2-sum-inverse}.

    If one side of $(V(\pi A),V(\pi A\ct))$ completely covers the vertex set of $G_2$, then both vertices of $t$ lie either in $G_2[[\pi A]]$ or both in $G_2[[\pi A\ct]]$. If neither side of $(V(\pi A),V(\pi A\ct))$ completely covers the vertex set of $G_2$, then $(V(\pi A),V(\pi A\ct))$ is a proper $2$-separation.
    Since $t$ is a tuplet, it is totally-nested and hence does not cross the separator $s$ of $(V(\pi A)),V(\pi A\ct))$. Hence we similarly conclude that both vertices of $t$ lie either in $G_2[[\pi A]]$ or both in $G_2[[\pi A\ct]]$. If a vertex of $t$ lies in $G_2[[\pi A]]$ but is not one of the vertices in the separator $s$, then all $2$-components at $t$ are $A$-ffected, but this contradicts our hypothesis. If a vertex of $t$ lies in $G_2[[\pi A\ct]]$ but is not one of the vertices in the separator $s$, then we claim that all but one of the $2$-components at $t$ are not $A$-ffected, in which case we are trivially done. Indeed, in this case at least one of the vertices of $s$ is internal to one of the $2$-components at $t$.

    So we assume that $s$ and $t$ coincide. Let $\phi_2'$ be the embedding induced on $G_2[[\pi A:e]]$ by $\phi_2$. Let $\phi_2''$ be the embedding induced on $G_2[[\pi A\ct: e]]$ by $\phi_2$. The embedding $\phi_2$ can be obtained via \autoref{lem:2-sum_embeddings} from $\phi_2'$ and $\phi_2''$ since $G_2$ is the $2$-sum of $G_2[[\pi A: e]]$ and $G_2[[\pi A\ct: e]]$ along $e$. The graph  $G_2[[\pi A: e]]$ consists exactly of the $A$-ffected $2$-components at $t$ along with the edge $e$ and the graph $G_2[[\pi A\ct: e]]$ consists of the $2$-components at $t$ which are not $A$-ffected along with the edge $e$. We conclude that the $A$-ffected $2$-components at $t$ form a segment in the cyclic structure induced on them by $\phi_2$. Furthermore, we observe that the embedding $\phi_2'$ is the the unique embedding (up to reorientation) induced on $G_2[[\pi A: e]]$ by simultaneous embeddings of $(G_1\minorm G_2)$, since $G_1[[A: e]]$ admits a unique embedding (up to reorientation) and contains $G_2[[\pi A: e]]$ as a minor. This completes our proof.
\end{pf}

\begin{lem}\label{lem:affected_G_2_cycle}
    Let $(G_1\minorm G_2)$ be a $2$-connected temporal pair. Let $A$ be a $3$-connected Tutte leaf of $G_1$ and let $t$ be a tuplet in $G_2$. Suppose that every $2$-component at $t$ is $A$-ffected. Then there is a unique cyclic structure on the $2$-components at $t$ that is induced by any simultaneous embedding of $(G_1\minorm G_2)$.
\end{lem}
\begin{pf}
    We proceed similarly to the proof of \autoref{lem:affected_G_2} and assume that $(V(\pi A),V(\pi A^C))$ is a $2$-separation of $G_2$ and that its separator $s$ is nested with $t$ and $t$ is contained in $G_2[[\pi A]]$. Then either $t=s$, in which case we get our required result analogously to in the proof of \autoref{lem:affected_G_2}, or there is a vertex in $s$ internal to a $2$-component at $t$. In the latter case, we observe that the $2$-components at $t$ in $G_2$ are in bijection with the $2$-components at the shadow $\pi t$ of $t$ in $G_2[[\pi A: e]]$. In particular, all of the $2$-components are the same except for one in which we minor down to the edge $e$. This bijection is order-preserving between cyclic orders of the $2$-components in embeddings of $G_2$ and the embeddings that they induce on $G_2[[\pi A: e]]$. Since $G_2[[\pi A: e]]$ is a minor of $G_1[[A: e]]$, which admits a unique embedding up to reorientation, we have our required result.
\end{pf}

In the context of \autoref{lem:affected_G_2}, we say that the unique linear structure induced on the $A$-ffected $2$-components at $t$ is the \emph{linear structure suggested by $A$}. 
In the context of \autoref{lem:affected_G_2_cycle}, we say that the unique cyclic structure induced on the $A$-ffected $2$-components at $t$ is the \emph{cyclic structure suggested by $A$}.

\begin{figure}
    \centering
    \includegraphics[width=0.5\linewidth]{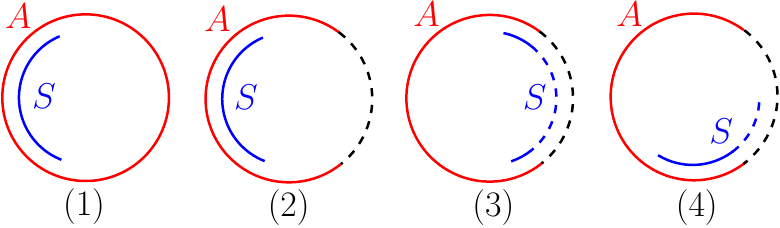}
    \caption{An illustration of the four ways in which a segment $S$ (in blue) can be affectively connected. The (cyclic or linear) structure suggested by $A$ is drawn as a red line.}
    \label{fig:affectively_connected}
\end{figure}

These definitions allow us to reason outside of the context of a given simultaneous embedding. Given a 2-separator $t$ of $G_2$ and a set $S$ of 2-components at $t$, we say that $S$ is \emph{affectively connected} if one of the following holds (see \autoref{fig:affectively_connected}):
\begin{enumerate}
\item[(1)] every $2$-component at $t$ in $S$ is $A$-ffected, and in the cyclic structure suggested by $A$, the set $S$ forms a segment; or
\item[(2)] there are $2$-components at $t$ that are not $A$-ffected, and in the linear structure suggested by $A$, the set $S$ forms a segment; or
\item[(3)] there are $2$-components at $t$ that are not $A$-ffected, and all such $2$-components are in $S$, and in the linear structure suggested by $A$, the $A$-ffected $2$-components in $S$ form the complement of a segment; or
\item[(4)] there are 2-components at $t$ that are not $A$-ffected, and some of them are in $S$ and some are not, and in the linear structure suggested by $A$, the $A$-ffected components in $S$ form a segment that covers one of the linear order's endpoints.
\end{enumerate}

Given a $2$-component $c$ at a tuplet $s$ of $G_3$, we denote by $x(c)$ the set of 2-components at $\pi s$ in the shadow of $c$. 

\begin{lem}\label{pre-obstruc}
Let $\tseq{G}$ be a $2$-connected temporal triple of planar graphs along with a Tutte leaf $A$ of $G_1$ and a $2$-component $c$ at a tuplet $s$ of $G_3$. If $\tseq{G}$ is simultaneously embeddable, then $x(c)$ is affectively connected. 
\end{lem}
\begin{pf}
    Let $(\phi_1,\phi_2,\phi_3)$ be a simultaneous embedding of $\tseq{G}$.
    By \autoref{lem:comp_segs}, the edges of $c$ incident with a vertex $y\in s$ form a segment in the rotator $\phi_3(y)$. By another application of \autoref{lem:comp_segs}, we conclude that the components in $x(s)$ form a segment in the cyclic structure induced on the $2$-components at $\pi s$ by $\phi_2$. By \autoref{lem:affected_G_2} and \autoref{lem:affected_G_2_cycle}, the $A$-ffected $2$-components form a segment in the cyclic structure induced on the $2$-components at $\pi s$ by $\phi_2$. The (linear or cyclic) structure of this segment is the (linear or cyclic) structure suggested by $A$. From this, we see that $x(c)$ is affectively connected.
\end{pf}

\begin{dfn}[amazingly-bad obstruction]\label{obstruction:amazing}
    An \emph{amazingly-bad obstruction} is a $2$-connected temporal triple $\tseq{G}$ of planar graphs along with a Tutte leaf $A$ of $G_1$ and a $2$-component $c$ at a tuplet $s$ of $G_3$ so that $x(c)$ is not affectively connected.
\end{dfn}

\begin{lem}\label{lem:amazingly_bad}
    No amazingly-bad obstruction is simultaneously embeddable.
\end{lem}
\begin{pf}[Proof:] immediate consequence of \autoref{pre-obstruc} and \autoref{obstruction:amazing}. 
\end{pf}

\subsection{Gaudy obstructions}

\begin{figure}
    \begin{minipage}{0.49\textwidth}
    \Large
    \scalebox{0.6}{\begin{tikzpicture}
	\begin{pgfonlayer}{nodelayer}
		\node [style=none] (42) at (-1.5, 1) {$\succeq$};
		\node [style=none] (127) at (7.5, 1) {$\preceq$};
		\node [style=black dot] (317) at (-8, 2.25) {};
		\node [style=black dot] (319) at (-6, 3.75) {};
		\node [style=black dot] (320) at (-6, -1.5) {};
		\node [style=black dot] (321) at (-4, 2.25) {};
		\node [style=none] (324) at (-6, 4.5) {$x_1$};
		\node [style=none] (325) at (-6, -2.25) {$y_1$};
		\node [style=black dot] (337) at (3, 3.25) {};
		\node [style=black dot] (338) at (3, 3.25) {};
		\node [style=black dot] (339) at (3, -1.25) {};
		\node [style=black dot] (340) at (3, 3.25) {};
		\node [style=none] (341) at (3, 4) {$x_2$};
		\node [style=none] (342) at (3, -2) {$y_2$};
		\node [style=black dot] (352) at (12, 4) {};
		\node [style=black dot] (353) at (12, -1.5) {};
		\node [style=none] (355) at (12, 4.75) {$x_3$};
		\node [style=none] (356) at (12, -2.25) {$y_3$};
		\node [style=black dot] (365) at (9.25, 1.75) {};
		\node [style=black dot] (366) at (11.5, 1.5) {};
		\node [style=black dot] (367) at (13, 1.25) {};
		\node [style=black dot] (368) at (15, 1.25) {};
		\node [style=none] (373) at (-1.5, -8) {$\succeq$};
		\node [style=none] (374) at (7.5, -8) {$\preceq$};
		\node [style=black dot] (375) at (-8, -6.75) {};
		\node [style=black dot] (376) at (-6, -5.25) {};
		\node [style=black dot] (377) at (-6, -10.5) {};
		\node [style=black dot] (378) at (-4, -6.75) {};
		\node [style=black dot] (389) at (3, -5.75) {};
		\node [style=black dot] (390) at (3, -5.75) {};
		\node [style=black dot] (391) at (3, -10.25) {};
		\node [style=black dot] (392) at (3, -5.75) {};
		\node [style=none] (393) at (3, -5) {$u$};
		\node [style=none] (394) at (3, -11) {$v$};
		\node [style=black dot] (402) at (12, -10.5) {};
		\node [style=black dot] (412) at (10.25, -5) {};
		\node [style=black dot] (413) at (12, -7.25) {};
		\node [style=black dot] (414) at (13.75, -5) {};
		\node [style=none] (420) at (-6, -3.25) {$P_1$};
		\node [style=none] (421) at (3, -3.25) {$P_2$};
		\node [style=none] (422) at (12, -3.25) {$P_3$};
		\node [style=none] (423) at (-6, -12) {$Q_1$};
		\node [style=none] (424) at (3, -12) {$Q_2$};
		\node [style=none] (425) at (12, -12) {$Q_3$};
		\node [style=small back] (438) at (-8.5, -0.75) {};
		\node [style=none] (439) at (-8.5, -0.75) {$1$};
		\node [style=small back] (440) at (-7.5, 0) {};
		\node [style=none] (441) at (-7.5, 0) {$2$};
		\node [style=small back] (442) at (-6.5, 0.25) {};
		\node [style=none] (443) at (-6.5, 0.25) {$3$};
		\node [style=small back] (444) at (-5.5, 0.25) {};
		\node [style=none] (445) at (-5.5, 0.25) {$4$};
		\node [style=small back] (446) at (-4.5, 0) {};
		\node [style=none] (447) at (-4.5, 0) {$5$};
		\node [style=small back] (448) at (-3.5, -0.75) {};
		\node [style=none] (449) at (-3.5, -0.75) {$6$};
		\node [style=small back] (450) at (0.25, 1) {};
		\node [style=none] (451) at (0.25, 1) {$1$};
		\node [style=small back] (452) at (1.25, 1) {};
		\node [style=none] (453) at (1.25, 1) {$2$};
		\node [style=small back] (454) at (2.25, 1) {};
		\node [style=none] (455) at (2.25, 1) {$3$};
		\node [style=small back] (456) at (3.75, 1) {};
		\node [style=none] (457) at (3.75, 1) {$4$};
		\node [style=small back] (458) at (4.75, 1) {};
		\node [style=none] (459) at (4.75, 1) {$5$};
		\node [style=small back] (460) at (5.75, 1) {};
		\node [style=none] (461) at (5.75, 1) {$6$};
		\node [style=small back] (462) at (9.5, -0.75) {};
		\node [style=none] (463) at (9.5, -0.75) {$1$};
		\node [style=small back] (464) at (11.25, 0.25) {};
		\node [style=none] (465) at (11.25, 0.25) {$2$};
		\node [style=small back] (466) at (13, 0) {};
		\node [style=none] (467) at (13, 0) {$3$};
		\node [style=small back] (468) at (10.25, -0.25) {};
		\node [style=none] (469) at (10.25, -0.25) {$4$};
		\node [style=small back] (470) at (12, 0.25) {};
		\node [style=none] (471) at (12, 0.25) {$5$};
		\node [style=small back] (472) at (14.5, -0.25) {};
		\node [style=none] (473) at (14.5, -0.25) {$6$};
		\node [style=small back] (474) at (-7, 3) {};
		\node [style=none] (475) at (-7, 3) {$a$};
		\node [style=small back] (476) at (-5, 3) {};
		\node [style=none] (477) at (-5, 3) {$b$};
		\node [style=small back] (478) at (10.5, 2.75) {};
		\node [style=none] (479) at (10.5, 2.75) {$p$};
		\node [style=small back] (480) at (11.75, 2.75) {};
		\node [style=none] (481) at (11.75, 2.75) {$q$};
		\node [style=small back] (482) at (12.75, 2.25) {};
		\node [style=none] (483) at (12.75, 2.25) {$r$};
		\node [style=small back] (484) at (14, 2.25) {};
		\node [style=none] (485) at (14, 2.25) {$s$};
		\node [style=small back] (486) at (14, 1.25) {};
		\node [style=none] (487) at (14, 1.25) {$t$};
		\node [style=small back] (488) at (0.25, -8) {};
		\node [style=none] (489) at (0.25, -8) {$1$};
		\node [style=small back] (490) at (1.25, -8) {};
		\node [style=none] (491) at (1.25, -8) {$2$};
		\node [style=small back] (492) at (2.25, -8) {};
		\node [style=none] (493) at (2.25, -8) {$3$};
		\node [style=small back] (494) at (3.75, -8) {};
		\node [style=none] (495) at (3.75, -8) {$4$};
		\node [style=small back] (496) at (4.75, -8) {};
		\node [style=none] (497) at (4.75, -8) {$5$};
		\node [style=small back] (498) at (5.75, -8) {};
		\node [style=none] (499) at (5.75, -8) {$6$};
		\node [style=small back] (501) at (-8.5, -9.75) {};
		\node [style=none] (502) at (-8.5, -9.75) {$1$};
		\node [style=small back] (503) at (-7.5, -9) {};
		\node [style=none] (504) at (-7.5, -9) {$2$};
		\node [style=small back] (505) at (-6.5, -8.75) {};
		\node [style=none] (506) at (-6.5, -8.75) {$3$};
		\node [style=small back] (507) at (-5.5, -8.75) {};
		\node [style=none] (508) at (-5.5, -8.75) {$4$};
		\node [style=small back] (509) at (-4.5, -9) {};
		\node [style=none] (510) at (-4.5, -9) {$5$};
		\node [style=small back] (511) at (-3.5, -9.75) {};
		\node [style=none] (512) at (-3.5, -9.75) {$6$};
		\node [style=small back] (513) at (-7, -6) {};
		\node [style=none] (514) at (-7, -6) {$a$};
		\node [style=small back] (515) at (-5, -6) {};
		\node [style=none] (516) at (-5, -6) {$b$};
		\node [style=small back] (517) at (9.5, -8.75) {};
		\node [style=none] (518) at (9.5, -8.75) {$1$};
		\node [style=small back] (519) at (11.5, -8.5) {};
		\node [style=none] (520) at (11.5, -8.5) {$2$};
		\node [style=small back] (521) at (13.5, -8.5) {};
		\node [style=none] (522) at (13.5, -8.5) {$3$};
		\node [style=small back] (523) at (10.5, -8.5) {};
		\node [style=none] (524) at (10.5, -8.5) {$4$};
		\node [style=small back] (525) at (12.5, -8.5) {};
		\node [style=none] (526) at (12.5, -8.5) {$5$};
		\node [style=small back] (527) at (14.5, -8.75) {};
		\node [style=none] (528) at (14.5, -8.75) {$6$};
		\node [style=small back] (529) at (11.25, -6.25) {};
		\node [style=none] (530) at (11.25, -6.25) {$e$};
		\node [style=small back] (531) at (12, -5) {};
		\node [style=none] (532) at (12, -5) {$f$};
		\node [style=small back] (533) at (12.75, -6.25) {};
		\node [style=none] (534) at (12.75, -6.25) {$g$};
	\end{pgfonlayer}
	\begin{pgfonlayer}{edgelayer}
		\draw [style=c1] (319) to (317);
		\draw [style=c1, bend right] (317) to (320);
		\draw [style=c1, bend left=75, looseness=1.50] (320) to (317);
		\draw [style=c2] (319) to (321);
		\draw [style=c2, bend right=15] (321) to (320);
		\draw [style=c2, bend right] (320) to (321);
		\draw [style=c2, bend left=75, looseness=1.50] (321) to (320);
		\draw [style=c1, bend left=15] (317) to (320);
		\draw [style=c1] (338) to (337);
		\draw [style=c1, bend right=60, looseness=1.50] (337) to (339);
		\draw [style=c1, bend right=270, looseness=2.00] (339) to (337);
		\draw [style=c2] (338) to (340);
		\draw [style=c1, bend right] (340) to (339);
		\draw [style=c2, bend right=60, looseness=1.50] (339) to (340);
		\draw [style=c2, bend left=90, looseness=2.00] (340) to (339);
		\draw [style=c2, bend left] (337) to (339);
		\draw [style=c0] (352) to (365);
		\draw [style=c0] (352) to (366);
		\draw [style=c0] (367) to (352);
		\draw [style=c0] (368) to (352);
		\draw [style=c0] (367) to (368);
		\draw [style=c2, bend right, looseness=0.75] (365) to (353);
		\draw [style=c2, bend left=15] (366) to (353);
		\draw [style=c1, bend left, looseness=0.75] (367) to (353);
		\draw [style=c1, bend right=60, looseness=1.25] (365) to (353);
		\draw [style=c1, bend right=45, looseness=0.75] (366) to (353);
		\draw [style=c2, bend left=45] (368) to (353);
		\draw [style=c1] (376) to (375);
		\draw [style=c1, bend right] (375) to (377);
		\draw [style=c1, bend left=75, looseness=1.50] (377) to (375);
		\draw [style=c2] (376) to (378);
		\draw [style=c2, bend right=15] (378) to (377);
		\draw [style=c2, bend right] (377) to (378);
		\draw [style=c2, bend left=75, looseness=1.50] (378) to (377);
		\draw [style=c1, bend left=15] (375) to (377);
		\draw [style=c1] (390) to (389);
		\draw [style=c1, bend right=60, looseness=1.50] (389) to (391);
		\draw [style=c1, bend left=90, looseness=2.00] (391) to (389);
		\draw [style=c2] (390) to (392);
		\draw [style=c1, bend right] (392) to (391);
		\draw [style=c2, bend right=60, looseness=1.50] (391) to (392);
		\draw [style=c2, bend left=90, looseness=2.00] (392) to (391);
		\draw [style=c2, bend left] (389) to (391);
		\draw [style=c2, bend right] (412) to (402);
		\draw [style=c2, bend left] (413) to (402);
		\draw [style=c1, bend left] (414) to (402);
		\draw [style=c1, bend right=60, looseness=1.25] (412) to (402);
		\draw [style=c1, bend right] (413) to (402);
		\draw [style=c2, bend left=60, looseness=1.25] (414) to (402);
		\draw [style=c0] (412) to (414);
		\draw [style=c0] (414) to (413);
		\draw [style=c0] (413) to (412);
	\end{pgfonlayer}
\end{tikzpicture}}
    \end{minipage}
    \hfill%
    \begin{minipage}{0.49\textwidth}
        \caption{Temporal triples $\tseq{P}$ and $\tseq{Q}$, neither of which has a simultaneous embedding. The graph $P_2$ is obtained from $P_1$ by contracting the edges $a$ and $b$ and obtained from the graph $P_3$ by contracting the edges $p,q,r,s,t$. The graph $Q_2$ is obtained from $Q_1$ by contracting the edges $a$ and $b$ and obtained from $Q_3$ by contracting the edges $e,f,g$.}
    \label{fig:gaudy}
    \end{minipage}
\end{figure}
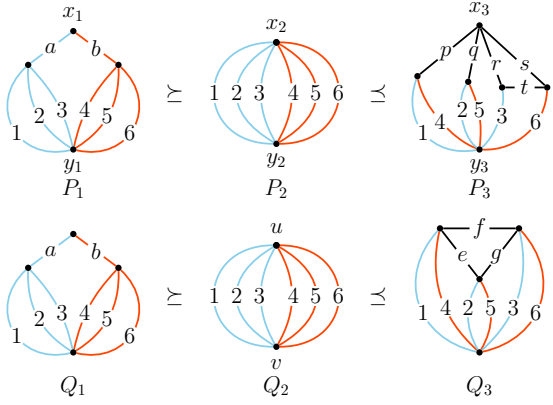

\begin{eg}\label{eg:gaudy}
    In \autoref{fig:gaudy}, we illustrate two temporal triples, $\tseq{P}$ and $\tseq{Q}$, neither of which have a simultaneous embedding. For the triple $\tseq{P}$, we observe the following: in any embedding of $P_2$ induced by an embedding of $P_1$, the edges from each $2$-component at $\{x_1,y_1\}$ in $P_1$ must form a segment in the rotator at $x_2$. In \autoref{fig:gaudy} this is illustrated by colouring the edges from each of the two $2$-components at $\{x_1,y_1\}$ by \rcol\ or \bcol, respectively. Similarly, in any embedding of $P_2$ induced by an embedding of $P_3$, the edges from each $2$-component at $\{x_3,y_3\}$ in $P_3$ must form a segment in the rotator at $x_2$. We observe that there is no cyclic order for the edges incident with $x_2$ that satisfies both of these constraints. Hence $\tseq{P}$ admits no simultaneous embedding.
    
    Now observe that the triple $\tseq{Q}$ exhibits analogous behaviour, except that this time instead of a tuplet in $Q_3$, there are three (crucially more than two) $2$-separators which project to $\{u,v\}$ in $Q_2$. Their branches induce segments at the rotator at the vertex $u$ in any embedding induced on $Q_2$ by an embedding of $Q_3$.
\end{eg}

Let $G$ be a $2$-connected graph and $b$ be a bond in $G$. For each embedding $\iota$ of $G$, the edges of $b$ form a cycle in the planar dual $G^*$ and hence come equipped with a cyclic order induced by $\iota$. By contracting one side of the bond to a single vertex, we may recover this cyclic order as the rotator at this vertex in the induced embedding.

\begin{dfn}[Bond-respectful]
    A \emph{bond graph} is pair $(G,b)$ where $G$ is a graph and $b$ is a bond thereof. A bond graph $(G,b)$ is \emph{partially-coloured} if some of the edges of $b$ (but maybe not all of them) are each assigned at least one colour. The edges outside of $b$ are not coloured. We say that a partially-coloured bond graph $(G,b)$ is \emph{$3$-overcoloured} if there is a colour, say \rcol, and three edges of $b$ each receiving colour \rcol\ as well as each at least one other colour. We say that a partially-coloured bond graph $(G,b)$ is \emph{$4$-overcoloured} with respect to a $4$-tuple of edges $(e_1,e_2,e_3,e_4)$ of $b$ if there is a colour, say \rcol, so that both $e_1$ and $e_3$ receive the colour \rcol\ and $e_2$ and $e_4$ both receive a colour other that \rcol\ (perhaps as well as \rcol).
    A \emph{bond-respectful} embedding of a partially-coloured bond graph $(G,b)$ is an embedding $\iota$ of $G$ so that $(G,b)$ is not $3$-overcoloured nor is $4$-overcoloured with respect to a $4$-tuplet $(e_1,e_2,e_3,e_4)$ that appears in that order in the cyclic order that $\iota$ induces on $b$. A bond graph $(G,b)$ itself is said to be \emph{bond-respectful} if all of its embeddings are bond-respectful.
\end{dfn}

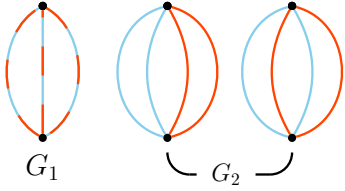
\begin{figure}
    \begin{minipage}{0.29\textwidth}
    \Large
        \scalebox{0.78}{\begin{tikzpicture}
	\begin{pgfonlayer}{nodelayer}
		\node [style=black dot] (337) at (-2, 4.5) {};
		\node [style=black dot] (338) at (-2, 4.5) {};
		\node [style=black dot] (339) at (-2, 0) {};
		\node [style=black dot] (340) at (-2, 4.5) {};
		\node [style=backing] (421) at (0, -1.25) {$G_2$};
		\node [style=black dot] (422) at (2.25, 4.5) {};
		\node [style=black dot] (423) at (2.25, 4.5) {};
		\node [style=black dot] (424) at (2.25, 0) {};
		\node [style=black dot] (425) at (2.25, 4.5) {};
		\node [style=black dot] (426) at (-6.25, 4.5) {};
		\node [style=black dot] (427) at (-6.25, 4.5) {};
		\node [style=black dot] (428) at (-6.25, 0) {};
		\node [style=black dot] (429) at (-6.25, 4.5) {};
		\node [style=none] (430) at (-6.25, -1) {$G_1$};
		\node [style=none] (431) at (-2, -0.5) {};
		\node [style=none] (432) at (2.25, -0.5) {};
		\node [style=none] (433) at (-1.25, -1.25) {};
		\node [style=none] (434) at (1.5, -1.25) {};
		\node [style=none] (436) at (1.5, -1.25) {};
		\node [style=none] (437) at (-1.25, -1.25) {};
	\end{pgfonlayer}
	\begin{pgfonlayer}{edgelayer}
		\draw [style=c1] (338) to (337);
		\draw [style=c1, bend right=75, looseness=1.25] (337) to (339);
		\draw [style=c2] (338) to (340);
		\draw [style=c1, bend right] (340) to (339);
		\draw [style=c2, bend right=75, looseness=1.25] (339) to (340);
		\draw [style=c2, bend left] (337) to (339);
		\draw [style=c1] (423) to (422);
		\draw [style=c1, bend right=75, looseness=1.25] (422) to (424);
		\draw [style=c2] (423) to (425);
		\draw [style=c2, bend right] (425) to (424);
		\draw [style=c2, bend right=75, looseness=1.25] (424) to (425);
		\draw [style=c1, bend left] (422) to (424);
		\draw [style=c1] (427) to (426);
		\draw [style=c1, bend right=60] (426) to (428);
		\draw [style=c2] (427) to (429);
		\draw [style=c1] (429) to (428);
		\draw [style=c1, bend right=60] (428) to (429);
		\draw [style=c2 dashed] (428) to (429);
		\draw [style=c2 dashed, bend left=60] (429) to (428);
		\draw [style=c2 dashed, bend left=60] (428) to (429);
		\draw [style=c0, bend right=45] (431.center) to (437.center);
		\draw [style=c0] (437.center) to (436.center);
		\draw [style=c0, bend right=45] (434.center) to (432.center);
	\end{pgfonlayer}
\end{tikzpicture}}
    \end{minipage}
    \hfill%
    \begin{minipage}{0.69\textwidth}
        \caption{Two bond graph $(G_1,b_1)$ and $(G_2,b_2)$, both of whose bonds are the full edge set of their respective graph. The graph $G_1$ does not admit a bond-respectful embedding because $b_1$ contains three edges having both colour \rcol\ and colour \bcol. The graph $G_2$ is depicted with two embeddings, one of which (left) is bond-respectful and one of which (right) is not. The graph bond graph $(G_2,b_2)$ itself is not bond-respectful because not all of its embeddings are bond-respectful.}
    \label{fig:bond-respectful}
    \end{minipage}
\end{figure}

See \autoref{fig:bond-respectful} for an illustration of the definition of bond-respectful embeddings.

Given a graph $G$ and $H\minorp G$ a minor thereof and $x$ some vertex of $H$, the \emph{branch set} of $x$ in $G$ is the maximal subset of vertices of $G$ which are all identified to $x$ in $H$.
\begin{lem}\label{lem:branch_set_bond}
    Let $G$ be a graph and $H\minorp G$ a minor of thereof. Let $x$ be a vertex of $H$ and denote by $Bx$ its branch set in $G$. If $G$ is $2$-connected, then the boundary $\partial Bx$ of $Bx$ is a bond.
\end{lem}
The above is a well-known fact, but we include a proof for convenience.
\begin{pf}[Proof of \autoref{lem:branch_set_bond}]
    It is clear that $\partial Bx$ is a cut in $G$. Now suppose that $G-\partial Bx$ has more than two components. One of these components must be $Bx$, since $Bx$ is connected. Let $G'$ be the minor of $G$ obtained by contracting $Bx$ back to a single vertex $x$. Now observe that $x$ forms a cut vertex between the shadows in $G'$ of the other two components of $G-\partial Bx$. But also notice that $H$ is a minor of $G'$ and that the shadow of each component in $H$ contains at least one edge. This contradicts the fact that $H$ is $2$-connected.
\end{pf}

\begin{dfn}\label{dfn:t-centred}
    Let $\tseq{G}$ be a $2$-connected temporal triple and $t$ be a tuplet in $G_1$ and choose a vertex $x$ in $t$. Assign a distinct colour to each $2$-component at $t$. Let $B\pi x$ be the branch set of $\pi x$ in $G_3$. By \autoref{lem:branch_set_bond}, the boundary $b:=\partial B\pi x$ is a bond of $G_3$; call it the \emph{$(t,x)$-centred bond of $G_3$}. We define the \emph{$(t,x)$-centred partial colouring of $(G_3,b)$} to be the partial colouring of the edges of $b$ which appear in $G_1$ by the colour of the $2$-component at $t$ in which they appear. All other edges remain colourless.
\end{dfn}

\begin{eg}
    Let $\tseq{G}$ be a $2$-connected temporal triple and $t$ be a tuplet in $G_1$ and choose a vertex $x$ in $t$. In the $(t,x)$-centred partial colouring of $G_3$, no edge is coloured using more than one colour.
\end{eg}

\begin{lem}\label{lem:b-respect_from_simultaneous}
    Let $\tseq{G}$ be a $2$-connected temporal triple. Let $t$ be a $2$-separator of $G_1$ and $x$ a vertex of $t$. Let $(\phi_1,\phi_2,\phi_3)$ be a simultaneous embedding of $\tseq{G}$. Let $b$ be the $(t,x)$-centred bond of $G_3$ along with the $(t,x)$-centred partial colouring. Then $\phi_3$ is bond-respectful for the bond graph $(G_3,b)$.
\end{lem}
\begin{pf}
    Suppose for a contradiction that $\phi_3$ is not bond-respectful.
    Since in the $(t,x)$-centred partial colouring of $(G_3,b)$ no edge receives more than one colour, there are edges $(e_1,e_2,e_3,e_4)$ in $b$, appearing in the cyclic order $(e_1,e_2,e_3,e_4)$ in the cyclic order induced on $b$ by $\phi_3$, so that $e_1$ and $e_3$ share a colour, say \rcol, and $e_2$ and $e_4$ both receive some colour that is not \rcol.
    By the construction of the $(t,x)$-centred bond $b$, we observe that the coloured edges of $b$ in $G_3$ are exactly the edges incident with $\pi x$ in $G_2$, where $\pi x$ denotes the shadow of $x$ in $G_2$. Hence the \rcol\ edges do not form a segment in the rotator $\phi_2(\pi x)$. However, the colouring of these edges was determined by the inclusion into $2$-components at $t$ in $G_1$. That is, the \rcol\ edges incident with $\pi x$ in $G_2$ are exactly the edges from some $2$-component $c$ at $t$ in $G_1$. But by \autoref{lem:comp_segs}, we know that the edges from $c$ must form a segment in the rotator $\phi_1(x)$. But then $\phi_1$ does not induce $\phi_2$ and this contradicts the fact that $(\phi_1,\phi_2,\phi_3)$ is a simultaneous embedding.
\end{pf}

In the remainder of this subsection, we will prove a characterisation of when bond-respectful embeddings exist, for partial colourings using two colours. We will then use this characterisation to translate \autoref{lem:b-respect_from_simultaneous} into another obstruction to simultaneous embeddability.

\begin{dfn}[Bond minors]
    Let $(G,b)$ be a bond graph along with a partial colouring. Let $e_1,\ldots,e_k$ be a class of parallel edges of $G$. We say that $e_1,\ldots,e_k$ is \emph{mergeable} if it is a maximal parallel class and the $2$-separator spanned by the class has at most one $2$-component which is not a single edge. The graph obtained by \emph{merging} $e_1,\ldots,e_k$ is $J:=G-\{e_2,\ldots,e_k\}$. The edge $e_1$ in $J$ is coloured using the union of the colours of each edge from $e_1,\ldots,e_k$. A \emph{bond-minor} of $G$ is a minor of $H$ of $G$ along with the restriction $b\restric H$ of $b$ to $H$ so that $H$ is obtained from $G$ by successive applications of the following two operations: firstly by contracting edges outside of $b$, or secondly by merging mergeable classes.
\end{dfn}

\begin{lem}\label{bondy}
    Let $(G,b)$ be a bond graph and $(H,b')$ a bond-minor of $(G,b)$. Then $(H,b')$ is a bond graph.
\end{lem}
The above is a well-known fact, but we include a proof for convenience.
\begin{pf}[Proof of \autoref{bondy}]
    We are required to show that $b'$ is a bond of $H$. Recall that a bond of $G$ is a circuit in the dual matroid of $G$. Deleting of edges from $b$ corresponds to contraction of edges in the dual. Thus $b$ remains a circuit in the dual setting and hence a bond in $H$.
\end{pf}

\begin{lem}\label{lem:bond_segment}
    Let $G$ be a $2$-connected graph and $b$ a bond in $G$. Let $t$ be a $2$-separator of $G$ and $c$ a $2$-component at $t$. Then in any embedding $\iota$ of $G$, the edges of $b$ which are in $c$ form a segment in the cyclic order induced on $b$ by $\iota$.
\end{lem}
\begin{pf}
    Let $\iota$ be an embedding of $G$. Consider the dual graph $G^*$ given by this embedding. Then the bond $b$ corresponds to a cycle in $G^*$. Furthermore, there are a pair of vertices $f,g$ in this cycle representing the faces that bound the component $c$ at $t$. That is, the two faces in $\iota$ which meet $c$ and both the vertices in $t$. In the cycle formed by $b$ in $G^*$, the vertices $f$ and $g$ bound a segment containing exactly the edges from $c$.
\end{pf}

The next lemma says that the class of partially-coloured bond graphs which admit bond-respectful embeddings is closed under taking bond-minors.

\begin{lem}\label{lem:b-respect_from_above}
    Let $(G,b)$ be partially-coloured bond graph. Let $(H,b')$ be a bond-minor of $(G,b)$. Any embedding of $H$ that is induced by a bond-respectful embedding for $(G,b)$ is itself bond-respectful for $(H,b')$.
\end{lem}
\begin{pf}
    By induction, we assume that a $H$ is obtained from $G$ by a single operation. If this operation is contraction of an edge outside of $b$, then we are done trivially. So assume that $H$ is obtained from $G$ by merging to an edge $e$ one mergeable class $e_1,\ldots,e_k$ in $(G,b)$ and let $b'$ be the restriction of $b$ to $H$. Let an embedding $\iota$ of $G$ be given. We first make the following claim.
    \begin{claim}\label{claim:bondy_lift_seg}
        The edges $e_1,\ldots,e_k$ form a segment in the cyclic order induced on $b$ by $\iota$.
    \end{claim}
    \begin{cproof}
        Since otherwise we are done trivially, suppose that $k\geq 2$ and that there is some other $2$-component $c$ at the $2$-separator $t$ spanned by $e_1,\ldots,e_k$. Since the class $e_1,\ldots,e_k$ is mergeable, it accounts for all of the $2$-components at $t$ except for $c$. By \autoref{lem:bond_segment}, the edges of $b$ that are in $c$ form a segment in the cyclic order induced on $b$ by $\iota$. Hence so does the complement of this set of edges in $b$, which is exactly the set $\{e_1,\ldots,e_k\}$.
    \end{cproof}
    
    Assume that the embedding $\iota\restric H$ that $\iota$ induces on $H$ is not a bond-respectful embedding. Suppose first that $\iota\restric H$ is $3$-overcoloured. Then there are three edges $a,b$ and $c$ in $b'$ coloured using a shared colour \rcol\ and each coloured some other colour in addition to \rcol. If none of these three edges are $e$, then we get that $\iota$ is $3$-overcoloured and hence not bond-respectful. So suppose without loss of generality that $a=e$.
    
    Now if one of the edges from $e_1,\ldots,e_k$, say $e_i$, is coloured both \rcol\ and some other colour in addition to \rcol, then $\iota$ is $3$-overcoloured, witnessed by $e_i,b,c$. So suppose that each \rcol\ edge in $e_1,\ldots,e_k$ receives no other colour. But since $e$ receives more than one colour there is one edge, say $e_i$, which is \rcol\ and one edge, say $e_j$ coloured using a colour distinct from \rcol. Now $e_i,e_j,a,b$ witness that $\iota$ is $4$-overcoloured.

    Finally, assume that $\iota\restric H$ is not $3$-overcoloured. Then $b'$ has at most two edges coloured both \rcol\ and each some other colour, but $b'$ contains four edges $r_1,b_1,r_2,b_2$ which witness the fact that $\iota\restric H$ is $4$-overcoloured. If $e$ is not amongst the edges $r_1,b_1,r_2,b_2$, then this quadruple witnesses that $\iota$ is $4$-overcoloured. So assume that $e$ is one of the edges $r_1,b_1,r_2,b_2$. Then we find $r_1',b_1',r_2',b_2'$ in $b$, with one edge from the segment $e_1,\ldots,e_k$, so that $r_1'$ and $r_2'$ receive colour \rcol\ and $b_1'$ and $b_2'$ receive some colour aside from \rcol. The tuplet $(r_1',b_1',r_2',b_2')$ witnesses that $\iota$ is $4$-overcoloured and hence is not bond-respectful.
\end{pf}

\begin{lem}\label{lem:seps_seps}
    Let $G$ be a $2$-connected graph and $b$ a bond in $G$. Let $\{x,y\}$ be $2$-separator of $G$ which separates some edges of $b$. Then the vertices $x$ and $y$ are separated by the bond~$b$.
\end{lem}
The above is a well-known fact, but we include a proof for convenience.
\begin{pf}[Proof of \autoref{lem:seps_seps}]
    Let $e$ and $f$ denote a pair of edges in $G$ separated by $\{x,y\}$. Since $b$ is a bond, there is a cycle $o$ in $G$ which meets $b$ exactly in $e$ and $f$. Then observe that $o-e-f$ is a pair of disjoint paths, each on a different side of $b$, both of which meet the separator $\{x,y\}$. This completes our proof.
\end{pf}

Let $(G,b)$ be a partially-coloured bond graph. Let $c$ be a $2$-component at some totally-nested $2$-separator of a $2$-connected graph $G$. A \emph{branch corresponding to $c$} is a minor $G[[E(c):e]]$ from \autoref{2-sum-inverse}, obtained from $G$ by minoring all of the $2$-components at $t$ except for $c$ down to a single edge $e$. A \emph{bond-branch} is a bond-minor of $(G,b)$ that is also a branch of~$G$.
\begin{lem}\label{lem:b-branch_construction}
    Let $(G,b)$ be a partially-coloured bond graph. Then for every $2$-component $c$ at a totally-nested $2$-separator $t$ of $G$, there is a bond-branch of $G$ which is a branch corresponding to $c$.
\end{lem}
\begin{pf}
    Let $c$ be a $2$-component at some totally-nested $2$-separator $t=\{x,y\}$ of $G$ and first suppose that more than one of the $2$-components at $t$ contains an edge of $b$. Observe by \autoref{lem:seps_seps} that the vertices $x$ and $y$ lie on different sides of the bond $b$. Now denote by $X$ and $Y$ be the subgraphs corresponding to the sides of $b$ in $G$, so that $X$ contains $x$ and $Y$ contains $y$. Let $H$ be the graph obtained from $G$ by contracting all of the edges of $X$ and $Y$ except those in $c$. Observe that our new graph is obtained from $G$ be replacing the $2$-components at $t$ except $c$ by a parallel class containing exactly the edges in these $2$-components that lie in $b$. This class is a mergeable parallel class and so we merge it into a single torso edge. This completes our construction.

    Now instead suppose that no more than one of the $2$-components at $t$ contains an edge of $b$. So either $c$ contains no edges of $b$ or $c$ contains all the edges of $b$. In both cases, we obtain our bond-branch by replacing all of the edges outside of $c$ be a single torso edge which doesn't come from $b$.
\end{pf}

Let $(G,b)$ be a partially-coloured bond graph. A \emph{bond-torso} of $(G,b)$ is a bond-minor of $(G,b)$ that is also a torso of a bag of the Tutte decomposition of $G$.

\begin{lem}\label{lem:b-torso_construction}
    Let $(G,b)$ be a partially-coloured bond graph. Then for every bag $w$ of the Tutte decomposition of $G$, there is a bond-torso of $G$ which is a torso of $w$.
\end{lem}
\begin{pf}
    For each $2$-separator $t$ of $G$ in $w$, apply \autoref{lem:b-branch_construction} to the $2$-component at $t$ which contains~$w$.
\end{pf}
We remark that the above construction is unique, up to the choice of edges which replace branches at $w$. However it is easy to see that the partial colouring is unique. For this reason, we often say `the' bond-torso instead of `a' bond-torso.

Let $(G,b)$ be a $2$-connected bond graph with a partial \rcol/\bcol-colouring of its bond. A tuplet $t$ of $G$ is called \emph{gaudy} if it has at least three $2$-components which each contain some \rcol\ edge and some \bcol\ edge.
\begin{dfn}[Gaudy]
    Let $(G,b)$ be a $2$-connected bond graph with a partial \rcol/\bcol-colouring of its bond. We say that this colouring is \emph{gaudy} if either of the following two conditions are met:
    \begin{enumerate}
        \item the graph $G$ contains a gaudy tuplet; or
        \item there is a $3$-connected bond-minor of $(G,b)$ which is not bond-respectful.
    \end{enumerate}
\end{dfn}

\begin{rem}
    We note that to test condition (ii) of the above definition of gaudy, it suffices to check only the maximal $3$-connected bond-minors of $G$. That is, the bond-torsos of the $3$-blocks from the Tutte decomposition of $G$, of which there are only linearly many (up to Tutte-equivalence). Furthermore, testing whether a $3$-connected bond graph is bond-respectful is easy since $3$-connected graphs admit unique embeddings up to reorientation, and hence induce unique cyclic structures on their bonds.
\end{rem}

\begin{cor}\label{cor:gaudy_disrespect}
    Let $(G,b)$ be a bond graph along with a partial \rcol/\bcol-colouring. If this partial colouring is gaudy, then $G$ admits no $b$-respectful embedding.
\end{cor}
\begin{pf}
    Suppose we are in case (i) of the definition of gaudy. Let $c_1,c_2,c_3$ be three $2$-components at a tuplet $t$ of $G$ so that all of $c_1,c_2,c_3$ contain both \rcol\ and \bcol\ edges. Applying \autoref{lem:bond_segment} to each of the $2$-separations $(V(c_i),V(c_i^C))$, we obtain that there are three distinct subsets $a_i=c_i\cap b$ of the bond $b$ so that in every embedding $\iota$ of $G$, the sets $a_i$ each from a segment in the cyclic order induced on $b$. Since each $a_i$ contains both \rcol\ and \bcol\ edges, no embedding of $(G,b)$ is bond-respectful.
    Now suppose we are in case (ii) of the definition of gaudy. Then by \autoref{lem:b-respect_from_above}, no embedding of $(G,b)$ is bond-respectful and this completes our proof.
\end{pf}

\begin{lem}\label{lem:glue_colourfully}
    Let $(G,b)$ be a bond graph along with a partial \rcol/\bcol-colouring. Suppose that this partial colouring is not gaudy. Let $t$ be a totally-nested $2$-separator of $G$. Suppose that the bond-branch corresponding to each $2$-component $c$ at $t$ admits a bond-respectful embedding. Then $(G,b)$ admits a bond-respectful embedding.
\end{lem}
\begin{pf}
    Since otherwise we would be trivially done, assume that $t$ has more than one $2$-component containing edges from $b$. Let $c_1,\ldots,c_m$ denote the $2$-components at $t$ and $C_i$ their respective branches for $i\in[m]$, given by \autoref{lem:b-branch_construction}. By \autoref{lem:seps_seps}, for each $i\in[m]$, the edge $e_i$ spanning $t$ in $C_i$ is part of the restriction of the bond $b$ to $C_i$ and is coloured using all of the colours present in the other components at $t$ aside from $c_i$. By \autoref{lem:2-sum_embeddings}, a choice of embedding $\iota_i$ for each of the bond-branches $C_i$, respectively, corresponds to an embedding $\iota$ of $G$, up to the cyclic order in which the $b$-branches attach. For each $i\in[m]$, denote by $s_i$ the segment of the cyclic order of $b\cap E(C_i)$ induced by $\iota_i$ obtained by removing $e_i$. By \autoref{lem:bond_segment}, we observe that the cyclic structure on $b$ induced by $\iota$ can be obtained by joining the segments $s_i$ in the chosen cyclic order.
    Since the colouring of $G$ is not gaudy, there are at most two $2$-components at $t$ that contain more than one colour. There are three cases.
    
    \textbf{Case 1.} There are two $2$-components at $t$, say $c_1$ and $c_2$, that contain both \rcol\ and \bcol\ edges. By reorienting $\iota_1$ if necessary, we assume that in the cyclic order induced on the bond $b\cap C_1$, the edge $e_1$ appears between a \rcol\ and a \bcol\ edge in the order $(\rcol,e_1,\bcol)$. We choose $\iota_2$ so that, in the cyclic order induced on the bond $b\cap C_2$, the edge $e_2$ appears between a \rcol\ and a \bcol\ edge in the order $(\bcol,e_2,\rcol)$.
    We then choose the cyclic order at $t$ in the following way. First comes the component $c_1$, then all of the purely \bcol\ components, then the component $c_2$, then all of the purely \rcol\ components. Hence $b$ is composed of four segments in the following way: first the segment $s_1$, then the segments $s_i$ having only \bcol\ edges, then the segment $s_2$, then the segments $s_i$ having only \rcol\ edges. We observe that our embedding $\iota$ of $(G,b)$ is bond-respectful.

    \textbf{Case 2.} There is one $2$-component, say $c_1$, that has both colours. We proceed similarly to in the first case, choosing the cyclic ordering of the $2$-components at $t$ so that they form a \rcol\ and a \bcol\ segment either side of $c_1$ and that these appear in an order compatible with the embedding of $\iota_1$ of $C_1$.

    \textbf{Case 3.} There are no $2$-components containing more than one colour. Then we arrange the $2$-components $c_i$ into a \rcol\ segment and a \bcol\ segment, given by the unique colour that they contain. When a $2$-component contains no coloured edge, then we place it arbitrarily. The resulting embedding $\iota$ of $(G,b)$ is bond-respectful.
\end{pf}

\begin{lem}\label{lem:existence_b-respect}
    Let $(G,b)$ be a bond graph along with a partial \rcol/\bcol-colouring. If this colouring is not gaudy, then there exists a bond-respectful embedding of $(G,b)$.
\end{lem}
\begin{pf}
    First suppose that $G$ has no totally-nested $2$-separator. Then either $G$ is a cycle or a $3$-connected graph. In the first case $b$ consists of at most two edges and we conclude that every embedding of $G$ is trivially bond-respectful. In the latter case, the graph is its own torso and we are done because our colouring is not gaudy. So assume that $G$ has some totally-nested $2$-separator $t$. We apply \autoref{lem:glue_colourfully} to $G$ at $t$ and, by induction on the number of totally-nested $2$-separators, we are done.
\end{pf}

\begin{lem}\label{lem:respect_characterisation}
    Let $(G,b)$ be a planar bond graph along with a partial \rcol/\bcol-colouring. Then $(G,b)$ admits a bond-respectful embedding if and only if the partial \rcol/\bcol-colouring is not gaudy.
\end{lem}
\begin{pf}
    We combine \autoref{cor:gaudy_disrespect} and \autoref{lem:existence_b-respect}.
\end{pf}

\begin{dfn}[Gaudy obstruction]\label{obstruction:gaudy}
    A \emph{gaudy obstruction} is a temporal triple $\tseq{G}$ along with a $2$-separator $t$ of $G_1$ so that the following holds. There is a vertex $x$ in $t$ along with a colour, say \rcol, in the $(t,x)$-centred colouring of $G_3$ so that the colouring obtained from the $(t,x)$-centred colouring by recolouring every non-\rcol\ edge to \bcol\ is a gaudy colouring.
\end{dfn}

\begin{cor}\label{cor:gaudy_bad}
    No gaudy obstruction admits a simultaneous embedding.
\end{cor}
\begin{pf}[Proof:]
    follows from \autoref{lem:b-respect_from_simultaneous} and \autoref{lem:respect_characterisation}.
\end{pf}

\subsection{Inconsistent obstructions}
Recall that a list for a temporal sequence $(G_1 \minorm\ldots \minorp G_n)$ is a partition of the set of $3$-blocks, together with a choice of an embedding for each $3$-block, so that each part is closed under Tutte-equivalence and Tutte-equivalent $3$-blocks receive Tutte-equivalent embeddings. For simplicity, we will often talk about lists as though they contain one representative from each Tutte-equivalence class of $3$-blocks.
Recall that a simultaneous list-embedding of $(G_1 \minorm\ldots \minorp G_n)$ for a list $L$ is a simultaneous embedding of $(G_1 \minorm\ldots \minorp G_n)$
so that if two $3$-blocks are in the same partition class of $L$, then their embeddings are both equal to the embeddings chosen in $L$ or both are the reorientations of their respective embeddings chosen in $L$. Finally, recall that two lists $L$ and $M$ are equivalent if they both specify the same partition on the set of $3$-blocks and $M$ can be obtained from $L$ by reorienting all of the embeddings in some of the parts of $L$. It is immediate that a simultaneous embedding satisfies $L$ if and only if it satisfies an equivalent list $M$.

When presented with a list on a temporal sequence, we might like to test whether this list actually respects the structure of the graphs on short scale. The following obstruction is a list-based obstruction that roughly tests whether the list is `locally bad' at a particular temporal triple.
\begin{dfn}[Inconsistent obstruction]\label{obstruction:consistent}
    An \emph{inconsistent obstruction} is a temporal triple $\tseq{G}$ along with a list $L$ so that either:
    \begin{enumerate}
        \item there are $3$-blocks $w_1$ of $G_1$ and $w_2$ of $G_3$, both in the same part of $L$, along with a minor of $\tseq{G}$ of the form $(w_1\minorm\triple\minorp w_3)$ so that the embeddings $\alpha_1$ and $\alpha_3$ of $w_1$ and $w_3$, respectively, given by $L$ do not satisfy $\alpha_1\rightarrow\triple\leftarrow\alpha_3$; or
        \item there are $3$-blocks $w_1$ of $G_1$ and $w_3$ and $w_3'$ of $G_3$, all in the same part of $L$, along with minors of $\tseq{G}$ of the form $(w_1\minorm\triple\minorp w_3)$ and $(w_1\minorm\triple'\minorp w_3')$ so that the embeddings $\alpha_3$ and $\alpha_3'$ of $w_3$ and $w_3'$, respectively, given by $L$ do not satisfy both $\alpha\rightarrow\triple\leftarrow\alpha_3$ and $\alpha\rightarrow\triple'\leftarrow\alpha_3'$ for either choice of an embedding $\alpha$ for $w$.
    \end{enumerate}
\end{dfn}

\begin{lem}\label{lem:consist_bad}
    No inconsistent obstruction admits a simultaneous list-embedding.
\end{lem}
\begin{pf}
    Suppose that $\tseq{G}$ is a temporal triple admitting a simultaneous embedding. For each minor $(w_1\minorm\triple\minorp w_3)$ of $\tseq{G}$, with $w_1$ a $3$-block of $G_1$ and $w_3$ a $3$-block of $G_3$, the blocks $w_1$ and $w_3$ are linked and thus coupled by \autoref{not-disagree}. For every simultaneous embedding of $\tseq{G}$, we have that $\beta_1\rightarrow\triple\leftarrow\beta_3$ where $\beta_1$ is the embedding induced on $w_1$ and $\beta_3$ is the embeddings induced on $w_3$ by the simultaneous embedding. We observe that these are never present in a simultaneous list-embedding whose list satisfies condition (i) of the definition of an inconsistent obstruction. The proof for condition (ii) is analogous.
\end{pf}

\subsection{Uncombinable obstructions}\label{subsec:uncombinable}
One may assign multiple lists to a temporal sequence. In this case, a simultaneous list-embedding is required to satisfy all of them at once. Given a pair of lists, we would like to know from this data alone, whether it is possible to satisfy both at once. Towards this, we make the following definitions.

\begin{dfn}
    Let $(G_1\minorm\ldots\minorp G_n)$ be a temporal sequence along with a pair of lists $L$ and $M$. We say that $L$ and $M$ are \emph{combinable} if there exist lists $L'$ and $M'$, equivalent to $L$ and $M$, respectively, so that $L'$ and $M'$ specify identical embeddings for each $3$-block of $(G_1\minorm\ldots\minorp G_n)$.
\end{dfn}

\begin{lem}\label{lem:sim_gives_combinable}
    Let $(G_1\minorm\ldots\minorp G_n)$ be a temporal sequence along with a pair of lists $L$ and $M$. Let $\Phi$ be a simultaneous embedding of $(G_1\minorm\ldots\minorp G_n)$ and suppose that it satisfies both $L$ and $M$. Then $L$ and $M$ are combinable.
\end{lem}
\begin{pf}
    Since $\Phi$ satisfies $L$, each part of $L$ contains either the embeddings induced by $\Phi$ on its $3$-blocks or all of their reorientations. Hence there is a list $L'$, equivalent to $L$, so that the embeddings in $L'$ are exactly the embeddings induced on the the $3$-blocks by $\Phi$. Let $M'$ be obtained from $M$ analogously. Now observe that the embeddings in $L'$ are identical to the embeddings in $M'$. This completes our proof.
\end{pf}

\begin{dfn}[Uncombinable obstruction]\label{obstruction:uncombinable}
    An \emph{uncombinable obstruction} is a temporal sequence $(G_1\minorm\ldots\minorp G_n)$ along with a pair of lists $L$ and $M$ so that $L$ and $M$ are not combinable.
\end{dfn}

\begin{cor}\label{lem:uncombine_bad}
    No uncombinable obstruction admits a simultaneous embedding satisfying both of its lists. \qed
\end{cor}

\begin{dfn}[Combination]
    Let $(G_1\minorm\ldots\minorp G_n)$ be a temporal sequence along with a pair of lists $L$ and $M$. Suppose $L$ and $M$ are combinable. A \emph{combination} of $L$ and $M$ is a list $N$ defined as follows: since $L$ and $M$ are combinable there exists lists $L'$ and $M'$, equivalent to $L$ and $M$, respectively, so that $L'$ and $M'$ specify identical embeddings. We let the embeddings specified by our new list $N$ be those of $L'$ (or equivalently $M'$). Finally, let the partition of $N$ be the finest partition that is a coarsening of the partition of $L$ and the partition $M$. That is, its parts are the union of the parts $l$ of $L$ and $m$ of $M$ so that $l$ and $m$ have nonempty intersection.
\end{dfn}

\begin{lem}
    Let $(G_1\minorm\ldots\minorp G_n)$ be a temporal sequence along with a pair of lists $L$ and $M$. Suppose $L$ and $M$ are combinable. Then any two combinations of $L$ and $M$ are equivalent.
\end{lem}
\begin{pf}
    Suppose that $L'$ and $M'$ are lists, equivalent to $L$ and $M$, respectively, so that $L'$ and $M'$ specify identical embeddings. Suppose that $L''$ and $M''$ are also lists with this property. We are required to show that defining the combination of $L$ and $M$ via $L'$ is equivalent to defining it via $L''$. We know $L'$ and $L''$ are equivalent lists. Hence $L''$ can be obtained from $L'$ be reorienting all of the embeddings in some of the parts of $L'$. The list $M''$ is obtained from $M'$ by reorienting these exact embeddings. Hence the union of the parts of $L'$ that get reoriented are a union of parts of $M'$. This implies that they are a union of parts of $N'$, the combination of $L$ and $M$ using $L'$. That means we can obtain $M''$, the combination of $L$ and $M$ using $L''$, by reorienting all of the embeddings in these parts of $M'$. Hence $N'$ and $N''$ are equivalent.
\end{pf}

By the equivalence of combinations given by the above lemma, we often refer to a combination of a pair of lists $L$ and $M$ as \emph{the combination of $L$ and $M$}.

\begin{pf}[Proof of \autoref{lem:list-combining_machine}]
    Let $L$ and $M$ be a pair of lists for a temporal sequence $(G_1\minorm\ldots\minorp G_n)$ and suppose that together they do not form an uncombinable obstruction. Then we are required to find a list $N$ for $(G_1\minorm\ldots\minorp G_n)$ so that every simultaneous embedding of $(G_1\minorm\ldots\minorp G_n)$ satisfies $N$ if and only if it satisfies both $L$ and $M$.
    
    Let $N$ be the combination of $L$ and $M$. Let $L'$ and $M'$ be the equivalent lists to $L$ and $M$, respectively, whose embeddings are identical, that are used to define the combination $N$.
    Suppose first that $\Phi$ satisfies $N$. So, for each part $p$ of $N$, either all of the embeddings in $p$ are the ones induced by $\Phi$ or all of the embeddings in $p$ are the reorientation of those induced by $\Phi$. Let $N'$ be the equivalent list to $N$ obtained from $N$ by reorienting all of the embeddings in the parts that are of the latter type. Let $L''$ and $M''$ be obtained from $L'$ and $M'$, respectively, by reorienting all of the embeddings in all of the parts which are subsets of the parts of $N$ which had their embeddings reoriented. In this way, the embeddings specified by $N'$ are identical to the embeddings specified by $L''$ and identical to the embeddings specified by $M''$. But observe that these embeddings are exactly the ones induced by $\Phi$. Hence $\Phi$ satisfies both $L''$ and $M''$, and hence satisfies both $L$ and~$M$.
    
    Now suppose that $\Phi$ satisfies $L$ and $M$, and hence $L'$ and $N'$. Let $a$ and $b$ denote $3$-blocks in the same part of $N$ and $\alpha$ and $\beta$ their embeddings given by $N$, respectively. Since $a$ and $b$ are in the same part of $N$, $a$ is in a part $l$ of $L'$ which overlaps a part $m$ of $M'$ containing $b$. So let $c$ be a $3$-block in both $l$ and $m$ and $\gamma$ its embedding in $N$. Now since $\Phi$ satisfies $L'$, either $\Phi$ induces both $\alpha$ and $\gamma$, or it induces their reorientations. Similarly, since $\Phi$ satisfies $M'$, either it induces both $\gamma$ and $\delta$, or it induces their reorientations. Hence we get that $\Phi$ induces both $\alpha$ and $\beta$, or induces their reorientations. This completes our proof.
\end{pf}

\subsection{Obstructed}
We introduce the following definition to gather all five of our obstructions.
\begin{dfn}[Obstructed]
    Let $\Gcal=(G_1\minorm\ldots\minorp G_n)$ be a $2$-connected temporal sequence along with a pair of lists $L$ and $M$. We say that $\Gcal$ is \emph{obstructed} if any of the following five conditions hold:
    \begin{enumerate}
        \item one of the temporal triples in $\Gcal$ is a disagreeable obstruction (\ref{obs:disagreeable}); or \vspace{-.25cm}
        \item one of the temporal triples in $\Gcal$ is an amazingly-bad obstruction (\ref{obstruction:amazing}); or \vspace{-.25cm}
        \item one of the temporal triples in $\Gcal$ is a gaudy obstruction (\ref{obstruction:gaudy}); or \vspace{-.25cm}
        \item one of the temporal triples in $\Gcal$, along with the restriction of one of the lists $L$ or $M$, is an inconsistent obstruction (\ref{obstruction:consistent}); or \vspace{-.25cm}
        \item the temporal sequence $\Gcal$ along with $L$ and $M$ is an uncombinable obstruction (\ref{obstruction:uncombinable}).
    \end{enumerate}
\end{dfn}

\begin{lem}\label{bad_is_bad}
    Let $\Gcal=(G_1\minorm\ldots\minorp G_n)$ be a $2$-connected temporal sequence along with a pair of lists $L$ and $M$. If $\Gcal$ along with $L$ and $M$ is obstructed, then it has no simultaneous embedding that respects both $L$ and $M$.
\end{lem}
\begin{pf}
    We check each of the five obstructions in turn. We respectively apply \autoref{cor:disagreeable_bad}, \autoref{lem:amazingly_bad}, \autoref{cor:gaudy_bad}, \autoref{lem:consist_bad} and \autoref{lem:uncombine_bad}.
\end{pf}

We would ideally like to build a converse for the above lemma. Sadly this is not possible using obstructions which are temporal triples. Indeed, there are longer temporal sequences, all of whose temporal triples are simultaneously embeddable, but which are not themselves simultaneously embeddable. However, the main result of this paper, \autoref{thm:main}, shows that if none of these five obstructions are present, then there is a way to reduce the complexity of the temporal sequence, in such a way that simultaneous embeddability is preserved. This directly implies an algorithm for determining simultaneous embeddability.

\section{Making things a little nicer}\label{sec:nice}

\begin{rem}
In this remark, we shall briefly discuss why we introduced the notion of nice. We have already seen that for the algorithmic application \autoref{cor:poly}, we just need to make the sequence nice once and then this property is preserved throughout the process. All definitions for nice temporal sequences made in this paper, can also be generalised so that niceness would not be required, yet these generalisations would be much more technical, and the constructions in the proofs of the corresponding lemmas would blow up accordingly. Thus we took the route to first make the sequence nice, which is indeed straightforward to achieve in our context. 
And this way we can work with simpler definitions and have \lq nicer\rq\ proofs in lots of places, where a lot of technicalities are moved to the part where we show that niceness is preserved along the way. This strategy helps to separate ideas from technical computations, and we invite the reader to skip in a first read all the lemmas that show that niceness is preserved. 
\end{rem}

\subsection{Towards the proof of \autoref{lem:niceness_main}}
In this subsection we move towards the proof of \autoref{lem:niceness_main}.

\begin{dfn}[Accurate]
    Let $G$ and $X$ be a pair of graphs so that $G$ is a minor of $X$. Then for every $3$-block $w$ of $G$, there is a $3$-block $w'$ of $X$ which contains it as a minor. We call $w'$ a \emph{lift} of $w$. When every $3$-block in $G$ admits a unique lift up to Tutte-equivalence, then we say that $X$ is \emph{accurate for $G$}.
\end{dfn}

\begin{eg}\label{eg:3-block_shenanigans}
    Let $G$ and $X$ be planar graphs so that $G$ is a minor of $X$. Let $w$ be a $3$-block of $G$ and $w'$ its lift in $X$. Let $\alpha$ be one of the embeddings of $w$. Then one of the embeddings of $w'$ induces $\alpha$ on $w$ and the other induces its reorientation. In this sense, the embeddings of $w$ and $w'$ are coupled together. This coupling will be exploited by \autoref{lem:list_lifter}.
\end{eg}

\begin{lem}\label{lem:list_lifter}
    Let $\Hcal=(H_1\minorm\ldots\minorp H_n)$ be a temporal sequence of planar graphs and $\Gcal=(G_1\minorm\ldots\minorp G_n)$ be a temporal minor of the same length. Suppose that $L$ is a list for $\Gcal$. If for each $i\in[n]$, we have that $H_i$ is accurate for $G_i$, then there exists a list $L^*$ for $\Hcal$ so that the following holds: Let $\Sigma$ be a simultaneous embedding of $\Hcal$ and $\Sigma'$ be the simultaneous embedding it induces on $\Gcal$. Then $\Sigma'$ satisfies $L$ if and only if $\Sigma$ satisfies $L^*$.
\end{lem}
\begin{pf}
    We first deal with the case where $n=1$. So let $G$ and $H$ by graphs so that $G\minorp H$, let $L$ be a list for $G$ and suppose that $G$ is accurate for $H$.
    Each part $W$ of the list $L$ lifts uniquely to a set $W^*$ of $3$-blocks of $H$. We define a list $L^*$ for $G$ as follows. Take the parts of $L^*$ to be the sets $W^*$ and, for every $3$-block which is not a lift of a $3$-block of $G$, we include its Tutte-equivalence class as its own part. The embeddings that $L$ specifies for the $3$-blocks of $G$ lift to embeddings for the $3$-blocks of $H$ in the sets $W^*$. Let $L^*$ specify these embeddings for the $3$-blocks which are lifts and, for those which are not lifts, we pick an arbitrary embeddings for each Tutte-equivalence class. It is easy to verify that any embedding of $H$ satisfies $L^*$ if and only if the embedding it induces on $G$ satisfies $L$. For temporal sequences of longer length than one, the proof is completely analogous.
\end{pf}

We will see a proof of the following lemma in \autoref{sec:nice_graph}.
\begin{lem}\label{lem:nice_graph_main}
    Let $G$ be $2$-connected graph. Then there exists a $2$-connected graph $H\minorm G$ so that $H$ is nice and every embedding $G$ lifts to an embedding of $H$.
    Furthermore $H$ is accurate for $G$ and $|E(H)|\leq100\cdot(|E(G)|)^2$.
\end{lem}

\begin{pf}[Proof of \autoref{lem:niceness_main} assuming \autoref{lem:nice_graph_main}]
    Let $\Gcal$ be a temporal sequence of $2$-connected graphs and $L$ a list thereof. We are required to show that there exists an analogue $(\Hcal,L')$ of $(\Gcal,L)$ so that $\Hcal$ is nice and furthermore $\text{size}(\Hcal)\leq 100\cdot(\text{size}(\Gcal))^2$.
    We obtain the sequence $\Hcal$ by applying \autoref{lem:nice_graph_main} to each graph in $\Gcal$. Then our result follows by \autoref{lem:list_lifter}.
\end{pf}

\subsection{Stretching a graph}\label{sec:nice_graph}
We say that a graph $G$ is \emph{really simple} if it is both simple (in particular this means that $G$ has no parallel edges) and no $2$-separator of $G$ is spanned by an edge. 
\begin{rem}
In the remainder of this section, we often work with really simple graphs. This reduces the case analysis and makes our definitions cleaner. We note that this is justified, since it is easy to build a really simple graph from one which is not really simple by subdividing some edges.
\end{rem}

A \emph{subgraph separation} of a graph $G$ is a pair of subgraphs $(A,B)$ so that $V(G)=V(A)\cup V(B)$ and the edges of $A$ and $B$ partition the edge set of $G$.
A side $A$ of a subgraph separation $(A,B)$ of a graph $G$ is \emph{componental} if the graph $A\sm B:=(V(A)\sm V(B),E(A)\sm E(B))$ is connected. We say that a subgraph separation $(A,B)$ is \emph{componental} if one of its sides is componental.

Let $G$ be a $2$-connected graph $x$ be a vertex in $G$. Let $\Tcal_x$ be the tree decomposition obtained from the set of componental $2$-separations of $G$ whose separators contain the vertex $x$, and let $T_x$ be its decomposition tree. We call $\Tcal_x$ the \emph{local tree decomposition} of $G$ at $x$ and $T_x$ the \emph{local tree} of $G$ at $x$.
Observe that every edge of $G$ is contained in a unique bag of the local tree decomposition~$\Tcal_x$.

\begin{eg}
If a vertex $x$ is not in any $2$-separator of a graph $G$, then the local tree $T_x$ at $x$ is the trivial tree consisting of one node.
\end{eg}

\begin{dfn}[Stretching]
Let $G$ be a really simple $2$-connected graph. The \emph{stretching} of $G$ is a graph defined as follows.
The vertex set of the stretching is the disjoint union of the vertex sets of the local trees $T_x$ for every vertex $x$ of $G$. 
For a vertex $x$ of $G$ and each bag $b$ of the local tree decomposition $\Tcal_x$ at $x$ of $G$, denote by $x[b]$ the copy of $x$ in $T_x$ in the stretching of $G$.
For each edge $e$ between vertices $u$ and $v$ in $G$, denote the unique bag\footnote{The bag is unique since we use subgraph separations instead of the usual pair of subsets of the vertex set of $G$ to define our tree decomposition.} containing $e$ in the local tree decomposition $\Tcal_x$ at $x$ by $\beta(e)$ and add an edge labeled $e$ in the stretching between
$u[\beta(e)]$ and $v[\beta(e)]$. All remaining edges of the stretching are the edges of the local trees $T_x$.
\end{dfn}

\begin{eg}
    Let $G$ be a really simple $2$-connected graph and $H$ its stretching. Then $G\minorp H$. Indeed, it is clear that if we contract all the trees in $H$ corresponding to local trees $T_x$ then we recover $G$.
\end{eg}


\begin{lem}\label{lem:tutti_nice}
    Let $G$ be a really simple $2$-connected graph. Then the stretching of $G$ is nice.
\end{lem}

In preparation for the above, we prove the following.

\begin{lem}\label{lem:local_leaves}
    Let $G$ be a really simple $2$-connected graph and $H$ its stretching. Let $x$ be a vertex of $G$. Consider the local tree $T_x$ as a subgraph of $H$. Then every leaf of $T_x$ is incident in $H$ with an edge from $E(G)$.
\end{lem}
\begin{pf}
    Let $l$ be a leaf of $T_x$. Denote the unique edge in $T_x$ that is incident with $l$ by $f$. Then $f$ corresponds to a componental $2$-separation $(A,B)$ of $G$, where the componental side $B$ is equal to the leaf bag of $T_x$ corresponding to $l$. That is, $l=x[B]$. Each leaf bag of the local decomposition at $x$ contains at least one edge that is incident with $x$, and we pick such an edge $g$. Note that $g$ is an edge of $G$. By the definition of the stretching, one of the endvertices of $g$ in $H$ is $x[B]=l$.
\end{pf}

\begin{dfn}[Coaddition]
    Let $G$ and $J$ be graphs and $x$ a vertex of $G$. Let $\theta$ be a function from the set $E(x)$ of edges incident with $x$ in $G$ to the vertex set of $J$. Then a graph $H$ is \emph{obtained from $G$ by coadding $J$} with respect to $\theta$ if $J$ is a subgraph of $H$, we have $G=H/J$ and each edge $e$ of $E(x)$ is incident with the vertex $\theta(e)$ and the unique endvertex of $e$ from $G$ other than $x$. This graph $H$ is uniquely determined and we refer to it as \lq the\rq\ graph obtained from $G$ by coadding $J$ with respect to $\theta$.
\end{dfn}

\begin{lem}\label{lem:2-con_lifts}
    Let $G$ be a $2$-connected graph and $H$ a graph obtained from $G$ by coadding a connected graph $Y$ at some vertex $x$ of $G$ with respect to a function $\theta$ from the edges incident with $x$ in $G$ to the vertices of $Y$. If for every vertex $y$ in $Y$, each connected component of $Y-y$ contains an element in the image of $\theta$, then $H$ is $2$-connected.
\end{lem}
\begin{pf}
    Suppose for the contrapositive that $y$ is a cutvertex of $H$. Firstly observe that $y$ must be in $Y$ since otherwise each of the components of $H-y$ would have nonempty shadow in $G$ and $x$ would be a cutvertex thereof. Now let $A_1,\ldots,A_n$ be the $n\geq 2$ connected graphs so that $E(A_1)\cup\ldots\cup E(A_n)=E(G)$, $V(A_i)\cap V(A_j)=\{y\}$ for each $i,j$ and each $A_i$ contains exactly one component of $H-y$. Since $x$ is not a cutvertex of $G$, at most one $A_i$ contains edges from $G$. So there is some $k\in[n]$ so that $A_k$ contains only edges from $Y$. Since otherwise there would be an edge from outside of $Y$ in $A_k$, we get that $A_k$ contains at most one element in the image of $\theta$ (and if it does, this element is $y$). We conclude that $A_i-y$ is a connected component of $Y-y$ which contains no element in the image of $\theta$.
\end{pf}

\begin{cor}\label{cor:tutti_2-con}
    Let $G$ be a really simple $2$-connected graph. Then the stretching of $G$ is $2$-connected.
\end{cor}
\begin{pf}
    We combine \autoref{lem:local_leaves} and \autoref{lem:2-con_lifts}.
\end{pf}

Along with \autoref{cor:tutti_2-con}, the following \autoref{lem:tutti_comps_lift} and \autoref{lem:tutti_comps_proj} are the lemmas required to prove \autoref{lem:tutti_nice}. They characterise the tuplets of the stretching in terms of the tuplets of the original graph.

\begin{lem}\label{lem:tutti_comps_lift}
    Let $G$ be a really simple $2$-connected graph and $H$ its stretching. Suppose that $t=\{x,y\}$ is a $2$-separator in $G$ with $m\ge3$ many $2$-components $C_1,\ldots,C_m$. Then there exists a unique $2$-separator $s$ in $H$ with $m$ many $2$-components $D_1,\ldots,D_m$ so that $C_i\minorp D_i$ for all $i$. Furthermore, $s=\{x[t],y[t]\}$ and for all $i$, and $x[t]$ and $y[t]$ have unique neighbours in $D_i-x$ and $D_i-y$, respectively.
\end{lem}
\begin{pf}
    Let $i\in[m]$. Consider the pair $(C_i,A_i)$ where $A_i$ is obtained from $G$ by deleting the vertex set $C_i\sm \{x,y\}$. By assumption, $(C_i,A_i)$ is a componental $2$-separation of $G$ with separator $t$ and corresponds to some edge in $T_x$; denote it by $f_i^x$. The forest $T_x-f_i^x$ has two components; denote the component that contains $x[t]$ by $A_i^x$ and the other component by $C_i^x$. We define $A_i^y$ and $C_i^y$ like $A_i^x$ and $C_i^x$, respectively, but with `$y$' in place of `$x$'.

    Let $B_i$ be the branch set of $A_i\sm t$ along with $A_i^x$ and $A_i^y$. Let $D_i$ be the branch set of $C_i\sm t$ along with $C_i^x$ and $C_i^y$, the vertices $x[t]$ and $y[t]$, and the edges from each of $C_i^x$ and $C_i^y$ to $x[t]$ and~$y[t]$, respectively.

    \begin{claim}
        The pair $(D_i,B_i)$ is a $2$-separation of $H$ with a separator $\{x[t],y[t]\}$.
    \end{claim}
    \begin{cproof}
        The vertex sets $D_i$ and $B_i$ meet exactly in $\{x[t],y[t]\}$. We claim that no edges of $H$ join vertices in different sides of $(D_i,B_i)$. First note that the edges in the branch sets of $A_i\sm t$ and $C_i\sm t$ do not cross the separation. It remains to check the edges in the trees $T_x$ and $T_y$. Since $f_i^x$ corresponds to the $2$-separation $(C_i,A_i)$, the subtree $C^x_i$ meets only the branch set of $C_i$ and the subtree $A^x_i$ meets only the branch set of $A_i$. An analogue of this statement also holds with `$y$' in place of `$x$'. Finally, the edges $f^x_i$ and $f^y_i$ are incident with $x[t]$ and $y[t]$, respectively, and so also do not cross the separation.
    \end{cproof}

    The above claim completes the proof that there exists a $2$-separator $s$ in $H$ with $2$-components $D_1,\ldots,D_m$ so that $C_i\minorp D_i$ for all $i$. It remains to show that $s$ is the unique $2$-separator with this property. So consider an arbitrary $2$-separator $s'=\{p,q\}$ in $H$.
    If the branch set $T_x$ of $x$ in $H$ is not separated by the set $\{p,q\}$, then one of the $2$-components at $s'$ contains all of $T_x$. If we further have that at most one of $p$ or $q$ is a leaf of $T_x$, then this $2$-component also contains all but one of the minors $C_1,\ldots,C_m$, and so does not satisfy the required property. The analogue of this holds with `$y$' in place of `$x$' and hence $\{p,q\}$ must separate both $T_x$ and $T_y$. In particular, each of $T_x$ and $T_y$ must contain one of $p$ or $q$.
    
    Without loss of generality, assume that $p$ separates $T_x$ and $q$ separates $T_y$. Suppose that $p$ is not equal to $x[t]$. Let $i$ be the unique choice of index so that the edge $f_i^x$ separates $p$ from $x[t]$. Then we know that there is there is a $2$-component at $s'$ that contains all of $A_i^x$, and so contains an edge from each of $C_1,\ldots,C_m$ except $C_i$. So we must have that $p=x[t]$ and, by symmetry, $q=y[t]$ and we get $s'=\{x[t],y[t]\}=s$.
\end{pf}

\begin{lem}\label{lem:tutti_comps_proj}
    Let $G$ be a really simple $2$-connected graph and $H$ be its stretching. Let $s$ be a $2$-separator in $H$ with $m\ge3$ many $2$-components $D_1,\ldots,D_m$. Then there exists a separator $t$ in $G$ with $m$ many $2$-components $C_1,\ldots,C_m$ so that $s=\{x[t],y[t]\}$ and $C_i\minorp D_i$ for all $i$.
\end{lem}
\begin{pf}
    Let $t=\pi s$ be the shadow of $s$ in $G$. Denote by $p$ and $q$ the two vertices of $s$. For all $i$, let $C_i=\pi D_i$ be the shadow of $D_i$.
    \begin{claim}\label{claim:all_but_one}
        For all but at most one $i$, the graph $C_i\sm t$ is nonempty.
    \end{claim}
    \begin{cproof}
        Choose a $p$--$q$ path in $D_i$ that does not use the torso edge of $D_i$ from $p$ to $q$; denote it by $P_i$. If $C_i\sm t$ is empty, then $P_i$ is contained in a single local tree. Thus the endvertices $p$ and $q$ of $P_i$ are both in that local tree. Since no local tree contains a cycle and different local trees are vertex disjoint in $H$, there is at most one $i$ such that $C_i\sm t$ is empty.
    \end{cproof}

    \begin{claim}
        The set $t$ has two elements and, for all $i$, the graph $C_i\sm t$ is nonempty.
    \end{claim}
    \begin{cproof}
        By \autoref{claim:all_but_one} and the fact that $m\ge 3$, we have that $t$ is a separator of $G$ of order at most $2$. Since $G$ is $2$-connected the order is exactly $2$. So no $D_i$ contains a $p$--$q$ path contained in a single local tree. So $C_i\sm t$ is nonempty for every $i$.
    \end{cproof}

    By the last claim, $t$ is a $2$-separator with $2$-components $C_1,\ldots,C_m$ so that $C_i\minorp D_i$ for every $i$. By \autoref{lem:tutti_comps_lift}, the separator $\{x[t],y[t]\}$ is the unique $2$-separator in $H$ with $2$-components each containing exactly one of $C_1,\ldots,C_m$. Hence $s=\{x[t],y[t]\}$.
\end{pf}

\begin{pf}[Proof of \autoref{lem:tutti_nice}]
    Let $G$ be a really simple $2$-connected graph. We are required to show that the stretching $H$ of $G$ is nice.
    By \autoref{lem:tutti_comps_proj} and \autoref{lem:tutti_comps_lift}, each tuplet $s$ in $H$ is of the form $s=\{x[t],y[t]\}$ for a unique choice of tuplet $t=\{x,y\}$ in~$G$ that depends on~$s$. Furthermore, the choices $t$ lead to tuplets $\{x[t],y[t]\}$ which are disjoint. We also have, for a tuplet $t$ in $G$, the vertices $x[t]$ and $y[t]$ are the centres of stars whose leaves $x[b]$ and $y[b]$ arise from bags $b$ that are not tuplets. By another application of \autoref{lem:tutti_comps_proj} and \autoref{lem:tutti_comps_lift}, we conclude that $x[b]$ and $y[b]$ are not part of any tuplets in $H$. Hence all tuplets in $H$ have distance at least two, as desired.
\end{pf}

Now we prove the required topological properties of the stretching.

\begin{lem}\label{lem:tutti_expand}
    Let $G$ be a really simple $2$-connected graph and $H$ be its stretching. Then for every embedding $\phi$ of $G$, there exists an embedding $\psi$ of $H$ which induces it.
\end{lem}

To prove the above, we make the following preparations.

\begin{lem}\label{lem:block_boundaries}
    Let $\phi$ be an embedding of a really simple $2$-connected graph $G$ and let $x$ be a vertex of $G$. Then for each $2$-block $B$ of $G-x$, there is a unique cycle $C$ in $B$, bounding a disk $D$, so that the closed disc $D$ contains exactly the block $B$ and no other vertices or edges of $G$.
\end{lem}
\begin{pf}
    Let $B$ be a $2$-block of $G-x$. Let $\phi_B$ be the embedding induced on $B$ by $\phi$. Since $B$ is $2$-connected, each face in $\phi_B$ is bounded by a cycle (\cite{graphsonsurfaces2001} Prop.\ 2.1.5). Let $C$ be the cycle which bounds the face in $\phi_B$ from which $x$ was removed. We claim that $C$ satisfies our requirements. Since otherwise we would be done, we assume there is a vertex $y$ of $G-x$ that is not contained in $B$. Since $G$ is $2$-connected, there are two $y$--$C$ paths in $G$ that share only the vertex $y$. Since $y$ is not in $B$, we have that $y$ is separated from $B$ by a cutvertex in $G-x$, and so one of these paths must use $x$. Hence, $x$ and $y$ embed into the same side of $C$ in $\phi$. By assumption, one side of $C$ contains all of $B$, and we have shown that the other contains everything else.

    The cycle $C$ is unique with this property because any other cycle of $B$ would have a side in $\phi$ whose interior contains both an edge of $B$ and the vertex $x$.
\end{pf}

Let $G$ be a really simple $2$-connected graph. Let $T_x$ be the local tree at $x$ of $G$. Let $H_x$ be the graph on the vertex set $V(T_x)\cup\{x\}$ obtained from $T_x$ by adding the edges $e$ incident with $x$ in $G$ between the vertex $x$ and the node of $T_x$ which represents the unique bag containing $e$.

\begin{lem}\label{lem:block_stars}
    Let $G$ be a really simple $2$-connected graph and $x$ a vertex in $G$. Given an embedding $\phi$ of $G$, the graph $H_x$ has an embedding $\phi_x$ so that the rotator $\phi_x(x)$ at $x$ in $H_x$ is equal to the rotator $\phi(x)$ at $x$ in $G$.
\end{lem}
\begin{pf}
    For each $2$-block $B$ of $G-x$, let $C_B$ denote the unique cycle in $B$ given by \autoref{lem:block_boundaries}. Denote by $D_B$ the closed disk in $\phi$ bounded by $C_B$ that contains exactly $B$.

    Let $b$ be a bag of the local tree decomposition $\Tcal_x$ of $G$ at $x$. Consider the graph $b-x$. If $b-x$ contains a cutvertex $c$, then the set $\{c,x\}$ is a $2$-separator of $G$ which contains $x$, contradicting the fact that $b$ is a bag of the local tree decomposition. Hence $b-x$ is either $2$-connected, a single edge or single vertex. We get that $b-x$ is a $2$-block, bridge or cutvertex of $G-x$.

    Let $G^*$ be the graph obtained from $G$ by performing the following operations. For every bag $b$ of $\Tcal_x$ so that $b-x$ is a cutvertex $c$ of $G-x$, relabel the cutvertex $c$ to $b$. For every bag $b$ of $\Tcal_x$ so that $b-x$ is a $2$-block $B$ of $G-x$, replace $B$ by a wheel $W_b$ with rim the cycle $C_B$ and central vertex labelled by $b$. For every bag $b$ of $\Tcal_x$ so that $b-x$ is a bridge $B$ of $G-x$, replace $B$ by a path $P_b$ of length two with central vertex labelled by $b$.

    We construct an embedding $\phi^*$ of $G^*$ from $\phi$ in the following way. Each path $P_b$ we embed in the image of the bridge $b-x$. For each wheel $W_b$, we delete everything in the $2$-block $B=b-x$ except the boundary cycle $C_B$. Then $C_B$ bounds a face (the interior of $D_B$). Into this face we embed the central vertex $b$ of the wheel, along with all its incident edges.

    Now observe that for each pair $b,c$ of adjacent bags in $\Tcal_x$, there exists a unique path $P(b,c)$ in $G^*-x$ between the vertices labelled by $b$ and $c$ having length at most two. Also observe that for each edge $e$ incident with $x$ in $G$, contained in the unique bag $b$ of $\Tcal_x$, there is a unique path $P(e)$ in $G^*$ between $x$ and the vertex labelled by $b$ having length exactly two. In particular, the path $P(e)$ contains the edge $e$.

    Importantly, the paths $P(bc)$ and $P(e)$ are all edge-disjoint. Now let $G'$ be the graph obtained from $G^*$ in the following way. Contract each path $P(bc)$ to a single edge and labelling it by $bc$. Then contract the edge in each path $P(e)$ which is not $e$.
    By construction, we get that $G'$ contains a copy of $H_x$. So let $\psi$ be the embedding of $H_x$ that $\phi^*$ induces on this copy and observe that we never modified the rotator at $x$. Hence, the rotator $\psi(x)$ at $x$ in $H_x$ is equal to the rotator $\phi(x)$ at $x$ in $G$. This completes the proof.
\end{pf}

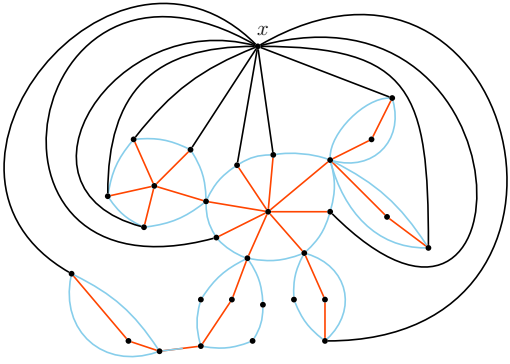
\begin{figure}
    \vspace{-1.2cm}
    \begin{minipage}{0.49\textwidth}
    \Large
        \scalebox{0.55}{\hspace{-2cm}\begin{tikzpicture}
	\begin{pgfonlayer}{nodelayer}
		\node [style=none] (128) at (-0.25, 8.75) {$x$};
		\node [style=black dot] (132) at (0, 0) {};
		\node [style=black dot] (134) at (0.25, 2.75) {};
		\node [style=black dot] (135) at (3, 2.5) {};
		\node [style=black dot] (136) at (3, 0) {};
		\node [style=black dot] (137) at (1.75, -2) {};
		\node [style=black dot] (138) at (-1, -2.25) {};
		\node [style=black dot] (139) at (-2.5, -1.25) {};
		\node [style=black dot] (140) at (-3, 0.5) {};
		\node [style=black dot] (141) at (-1.5, 2.25) {};
		\node [style=black dot] (142) at (6, 5.5) {};
		\node [style=black dot] (143) at (5, 3.5) {};
		\node [style=black dot] (144) at (-0.5, 8) {};
		\node [style=black dot] (145) at (-5.5, 1.25) {};
		\node [style=black dot] (146) at (-3.75, 3) {};
		\node [style=black dot] (147) at (-6.5, 3.5) {};
		\node [style=black dot] (148) at (-7.75, 0.75) {};
		\node [style=black dot] (149) at (-6, -0.75) {};
		\node [style=black dot] (150) at (-3.25, -4.25) {};
		\node [style=black dot] (151) at (-3.25, -6.5) {};
		\node [style=black dot] (152) at (-0.75, -6.25) {};
		\node [style=black dot] (153) at (-0.25, -4.5) {};
		\node [style=black dot] (154) at (-1.75, -4.25) {};
		\node [style=black dot] (155) at (-5.25, -6.75) {};
		\node [style=black dot] (156) at (-6.75, -6.25) {};
		\node [style=black dot] (157) at (-9.5, -3) {};
		\node [style=black dot] (158) at (2.75, -4.25) {};
		\node [style=black dot] (159) at (2.75, -6.25) {};
		\node [style=black dot] (160) at (1.25, -4.25) {};
		\node [style=black dot] (161) at (5.75, -0.25) {};
		\node [style=black dot] (162) at (7.75, -1.75) {};
	\end{pgfonlayer}
	\begin{pgfonlayer}{edgelayer}
		\draw [style=c2] (141) to (132);
		\draw [style=c2] (132) to (134);
		\draw [style=c2] (135) to (132);
		\draw [style=c2] (132) to (136);
		\draw [style=c2] (137) to (132);
		\draw [style=c2] (132) to (138);
		\draw [style=c2] (132) to (140);
		\draw [style=c2] (139) to (132);
		\draw [style=c1, bend right=15, looseness=0.75] (139) to (138);
		\draw [style=c1, bend right=15] (138) to (137);
		\draw [style=c1, bend right=15] (137) to (136);
		\draw [style=c1, bend right=15] (136) to (135);
		\draw [style=c1, bend right=15, looseness=0.75] (135) to (134);
		\draw [style=c1, bend right=15] (134) to (141);
		\draw [style=c1, bend right=15] (141) to (140);
		\draw [style=c1, bend right=15] (140) to (139);
		\draw [style=c2] (135) to (143);
		\draw [style=c2] (143) to (142);
		\draw [style=c1, bend left=60, looseness=1.25] (142) to (135);
		\draw [style=c1, bend left=45, looseness=0.75] (135) to (142);
		\draw [style=c0] (142) to (144);
		\draw [style=c0] (144) to (134);
		\draw [style=c0] (141) to (144);
		\draw [style=c1, bend right=15] (140) to (146);
		\draw [style=c1, bend right=15] (146) to (147);
		\draw [style=c1, bend right=15] (147) to (148);
		\draw [style=c1, bend right=15] (148) to (149);
		\draw [style=c1, bend right=15] (149) to (140);
		\draw [style=c2] (145) to (140);
		\draw [style=c2] (146) to (145);
		\draw [style=c2] (145) to (149);
		\draw [style=c2] (148) to (145);
		\draw [style=c2] (145) to (147);
		\draw [style=c0] (146) to (144);
		\draw [style=c0, bend right=15] (144) to (147);
		\draw [style=c2] (138) to (154);
		\draw [style=c2] (151) to (154);
		\draw [style=c1, bend right=15] (138) to (150);
		\draw [style=c1, bend right=15] (150) to (151);
		\draw [style=c1, bend right=15] (151) to (152);
		\draw [style=c1, bend right=15] (152) to (153);
		\draw [style=c1, bend right=15] (153) to (138);
		\draw [style=c2] (151) to (155);
		\draw [style=c2] (156) to (155);
		\draw [style=c2] (156) to (157);
		\draw [style=c0, in=150, out=135, looseness=1.75] (144) to (157);
		\draw [style=c0, bend left=45, looseness=1.25] (148) to (144);
		\draw [style=c0, in=165, out=165, looseness=1.75] (149) to (144);
		\draw [style=c1, bend left=15] (157) to (155);
		\draw [style=c1, bend left=60, looseness=1.25] (155) to (157);
		\draw [style=c1 dashed] (155) to (151);
		\draw [style=c2] (137) to (158);
		\draw [style=c2] (158) to (159);
		\draw [style=c1, bend right=15] (137) to (160);
		\draw [style=c1, bend right=15] (160) to (159);
		\draw [style=c1, bend right=60, looseness=1.25] (159) to (137);
		\draw [style=c0, in=30, out=0, looseness=2.50] (159) to (144);
		\draw [style=c0, in=15, out=-45, looseness=4.00] (136) to (144);
		\draw [style=c0, in=150, out=-165, looseness=3.50] (139) to (144);
		\draw [style=c2] (135) to (161);
		\draw [style=c2] (161) to (162);
		\draw [style=c1, in=285, out=-180] (162) to (135);
		\draw [style=c1, bend right=15] (162) to (135);
		\draw [style=c0, in=0, out=90, looseness=1.50] (162) to (144);
	\end{pgfonlayer}
\end{tikzpicture}}
    \end{minipage}
    \hfill%
    \begin{minipage}{0.49\textwidth}
        \caption{An illustration of the proof of \autoref{lem:block_stars}. We modify each part of the block-cutvertex decomposition of $G-x$ so that the resulting graph contains a copy of the decomposition tree. We do it in such a way that the node representing a particular bag of the decomposition tree is placed inside its modified bag.}
    \label{fig:block_stars}
    \end{minipage}
    \vspace{-1cm}
\end{figure}

Recall the following definition from \cite{carmesin2017embeddingI}. Let $G_1$ and $G_2$ be graphs sharing one vertex $v$ so that the set of labels $E_v$ of edges incident with $v$ in $G_1$ is identical to the set of labels of edges incident with $v$ in $G_2$. The \emph{vertex sum of $G_1$ and $G_2$ over $v$} is the graph obtained from the disjoint union of $G_1-v$ and $G_2-v$ by adding each edge $e\in E_v$ between the vertex $x$ incident with $e$ in only $G_1$ and the vertex $y$ incident with $e$ in only $G_2$. Observe that for $2$-connected $G_1$ and $G_2$, we have that each $G_i$ is a minor of the vertex sum. Indeed, we may contract all edges in $G_2$ that are not in $G_1$ to obtain $G_1$, and vice versa. We will need the following complementary lemma.

\begin{lem}[Lemma 3.2 in \cite{carmesin2017embeddingI}]\label{lem:vertex_sum}
    Let $G_1$ and $G_2$ be $2$-connected graphs sharing one vertex $v$ so that the set of labels of edges incident with $v$ in $G_1$ is identical to the set of labels of edges incident with $v$ in $G_2$. Let $\iota_1$ and $\iota_2$ be embeddings of $G_1$ and $G_2$, respectively, so that $\iota_1(v)=\iota_2(v)$. Then there exists an embedding of the vertex sum of $G_1$ and $G_2$ over $v$. In particular, this embedding induces $\iota_1$ on $G_1$ and the reorientation of $\iota_2$ on $G_2$. \qed
\end{lem}

\begin{pf}[Proof of \autoref{lem:tutti_expand}] 
    Let $G$ be a really simple $2$-connected graph and $H$ be its stretching. Let $\phi$ be an embedding of $G$. We are required to show that there exists an embedding $\psi$ of $H$ which induces $\phi$.
    We simply observe that the stretching $H$ is exactly the vertex sum of $G$ and the graph $H_x$ for each vertex $x$ of $G$. Then by \autoref{lem:block_stars}, there is an embedding $\phi_x$ of each $H_x$ so that the rotators $\phi(x)$ and $\phi_x(x)$ are identical. Finally, we apply \autoref{lem:vertex_sum} to obtain an embedding $\psi$ of $H$ which induces $\phi$ on $G$.
\end{pf}

\begin{lem}\label{lem:tutti_accurate}
    Let $G$ be a really simple $2$-connected graph and $H$ be its stretching. Then $H$ is accurate for $G$.
\end{lem}
\begin{pf}
    Let $c$ be a $3$-block of $H$ and $b$ and $b'$ two $3$-blocks of $G$ that are minors of $c$. We are required to show that $b$ and $b'$ are torsos of the same bag of the Tutte decomposition of $G$. Let $B$ and $B'$ denote the bags of the Tutte decomposition of $G$ corresponding to $b$ and~$b'$, respectively. Suppose for a contradiction that $B$ and $B'$ can be distinguished by some totally-nested $2$-separator $\{x,y\}$ of $G$.
    Let $(P,Q)$ be a componental $2$-separation whose separator is $\{x,y\}$. Let $e_x$ and $e_y$ denote the edges in the local trees $T_x$ at $x$ and $T_y$ at $y$, respectively, which correspond to the componental $2$-separation $(P,Q)$. By construction, $\{e_x,e_y\}$ is cut of $H$ into two components. One of the components of $H-\{e_x,e_y\}$ contains all but one edge of $b$ and the contains all but one edge of $b'$. This contradicts the fact that $c$ is $3$-connected, proving our result.
\end{pf}

Now we can prove our desired result regarding niceness for graphs.
\begin{pf}[Proof of \autoref{lem:nice_graph_main}]
    Given a $2$-connected graph $G$, we are required to show that there exists a $2$-connected graph $H\minorm G$ so that $H$ is nice and, for every embedding $\iota$ of $G$, there is an embedding of $H$ that induces $\iota$ on $G$ and furthermore $H$ is accurate for $G$ and $\text{size}(H)\leq 100\cdot(\text{size}(G))^2$. We may assume that $G$ is really simple, or else we can easily make it so by subdividing some edges. Now we take $H$ to be the stretching of $G$. By \autoref{lem:tutti_expand}, every embedding of $G$ lifts to an embedding of $H$. By \autoref{lem:tutti_nice}, the graph $H$ is nice and by \autoref{lem:tutti_accurate}, the graph $H$ is accurate for $G$. Furthermore, the number of edges in $H$ is the number of edges in $G$ plus the number of edges in each local tree. It is easy to see that the number of edges in any local tree is at most the number of edges in the graph $G$. Hence $\text{size}(H)\leq 100\cdot(\text{size}(G))^2$ as required.
\end{pf}

\section{Amazing analogues}\label{sec:amazing_analogues}
This section is dedicated to proving \myref{lem:A}.
This lemma was already stated in \autoref{sec:main} (where we showed that this lemma together with other lemmas implies the main result). For convenience of the reader we recall 
\begingroup
\def\thelem{\ref{lem:A}}
\begin{lem}[Amazing Analogue Lemma]
    Let $\Gcal=(G_1\minorm\ldots\minorp G_n)$ be a nice $2$-connected temporal sequence with a list $L$, and suppose that $A$ is a last upleaf of $\Gcal$ rooted at $i$ for some $i\in[n]$. Suppose further that some edge of $A$ appears in $G_{i+1}$. If $\Gcal$ is not a disagreeable obstruction or an amazingly-bad obstruction and the pair $(\Gcal,L)$ is not an inconsistent obstruction, then $(\Gcal,L)$ admits an amazing analogue for the pair $(A,i)$.
\end{lem}
\addtocounter{lem}{-1}
\endgroup

\subsection{Towards a proof of The \nameref{lem:A}}
In this subsection we state two facts and prove that they imply \myref{lem:A}. These facts will then be proved in the subsequent subsections.

\begin{dfn}[$A$-mazing]
    Given a temporal triple $\tseq{G}$ and an edge set $A$, we say that $\tseq{G}$ is \emph{$A$-mazing} (read: `amazing') if $A$ is an upleaf rooted at $i=1$ and for every tuplet $s$ of $G_3$, exactly one $2$-component at $s$ contains edges of $A$ (that is, is $A$-ffected).
\end{dfn}

\begin{thm}\label{thm-list}
Let $\Gcal=(G_1\minorm\ldots\minorp G_n)$ be a temporal sequence with list $L$, and $A$ an upleaf rooted at $i\in[n]$ that corresponds to a $3$-block. Suppose that the temporal triple $(G_i\minorm G_{i+1}\minorp G_{i+2})$ is $A$-mazing. Let $b=G_i[[A:e]]$ be the $3$-connected torso for~$A$. Assume that $\Gcal$ does not contain a disagreeable obstruction and $(G,\Lcal)$ is not an inconsistent obstruction.
Then there is a list $L'$ for the temporal sequence\vspace{-.25cm}
\[
(G_1\minorm\ldots\minorm G_{i-1} \minorp G_i[[A\ct:e]]\minorm G_{i+1}[[\pi A\ct:e]]\minorp G_{i+2}\minorm\ldots\minorp G_n)\vspace{-.25cm}
\]
so that this temporal sequence together with the list~$L'$ admits a simultaneous list-embedding if and only if $(G_1\minorm\ldots\minorp G_n)$ admits a simultaneous list-embedding for $L$. 
\end{thm}
The proof of \autoref{thm-list} will follow in \autoref{sec:amazing_reductions}.

\begin{figure}
    \centering
    {\Large \scalebox{0.7}{\begin{tikzpicture}
	\begin{pgfonlayer}{nodelayer}
		\node [style=none] (341) at (-6, 3.75) {$x$};
		\node [style=none] (342) at (-6, -3.75) {$y$};
		\node [style=black dot] (423) at (-6, 3) {};
		\node [style=black dot] (424) at (-6, -3) {};
		\node [style=black dot] (425) at (-7, 0) {};
		\node [style=black dot] (426) at (-5, 0) {};
		\node [style=black dot] (427) at (-3.5, 0) {};
		\node [style=black dot] (428) at (-8.5, 0) {};
		\node [style=black dot] (429) at (-2, 0) {};
		\node [style=none] (430) at (-1, 0) {$\succeq$};
		\node [style=none] (432) at (2.5, 3.75) {$x$};
		\node [style=none] (433) at (2.5, -3.75) {$y$};
		\node [style=black dot] (440) at (2.5, 3) {};
		\node [style=black dot] (441) at (2.5, -3) {};
		\node [style=black dot] (442) at (1.5, 0) {};
		\node [style=black dot] (443) at (3.5, 0) {};
		\node [style=black dot] (444) at (5, 0) {};
		\node [style=black dot] (445) at (0, 0) {};
		\node [style=black dot] (446) at (6.5, 0) {};
		\node [style=none] (447) at (7.5, 0) {$\preceq$};
		\node [style=black dot] (456) at (11, 3) {};
		\node [style=black dot] (457) at (11, -3) {};
		\node [style=black dot] (458) at (10, 0) {};
		\node [style=black dot] (459) at (15, 0) {};
		\node [style=black dot] (460) at (16.5, 0) {};
		\node [style=black dot] (461) at (8.5, 0) {};
		\node [style=black dot] (462) at (17.5, 0) {};
		\node [style=black dot] (464) at (14, 3) {};
		\node [style=black dot] (465) at (14, -3) {};
		\node [style=black dot] (466) at (12, 3) {};
		\node [style=black dot] (467) at (13, 3) {};
		\node [style=black dot] (468) at (12, -3) {};
		\node [style=black dot] (469) at (13, -3) {};
		\node [style=black dot] (470) at (12.5, 2) {};
		\node [style=black dot] (471) at (12.5, 4) {};
		\node [style=black dot] (472) at (12.5, -2) {};
		\node [style=black dot] (473) at (12.5, -4) {};
		\node [style=none] (476) at (-6, -5) {$G_1$};
		\node [style=none] (477) at (2.5, -5) {$G_2$};
		\node [style=none] (478) at (12.5, -5) {$G_3$};
		\node [style=black dot] (484) at (24.5, 3) {};
		\node [style=black dot] (485) at (24.5, -3) {};
		\node [style=black dot] (486) at (23.5, 0) {};
		\node [style=black dot] (487) at (28.5, 0) {};
		\node [style=black dot] (488) at (30, 0) {};
		\node [style=black dot] (489) at (22, 0) {};
		\node [style=black dot] (492) at (27.5, 3) {};
		\node [style=black dot] (493) at (27.5, -3) {};
		\node [style=black dot] (494) at (25.5, 3) {};
		\node [style=black dot] (495) at (26.5, 3) {};
		\node [style=black dot] (496) at (25.5, -3) {};
		\node [style=black dot] (497) at (26.5, -3) {};
		\node [style=black dot] (498) at (26, 2) {};
		\node [style=black dot] (499) at (26, 4) {};
		\node [style=black dot] (500) at (26, -2) {};
		\node [style=black dot] (501) at (26, -4) {};
		\node [style=none] (504) at (26, -5) {$J_3$};
		\node [style=small back] (505) at (-7.5, 1.25) {};
		\node [style=none] (506) at (-7.5, 1.25) {$1$};
		\node [style=small back] (507) at (-6.5, 1.25) {};
		\node [style=none] (508) at (-6.5, 1.25) {$2$};
		\node [style=small back] (509) at (-5.5, 1.25) {};
		\node [style=none] (510) at (-5.5, 1.25) {$3$};
		\node [style=small back] (511) at (-4.5, 1.25) {};
		\node [style=none] (512) at (-4.5, 1.25) {$4$};
		\node [style=small back] (513) at (-6, 0) {};
		\node [style=none] (514) at (-6, 0) {$p$};
		\node [style=medium back] (515) at (-3.25, 2.25) {};
		\node [style=none] (516) at (-3.25, 2.25) {$d_1$};
		\node [style=medium back] (517) at (-3.25, -2.25) {};
		\node [style=none] (518) at (-3.25, -2.25) {$d_2$};
		\node [style=small back] (519) at (1, 1.25) {};
		\node [style=none] (520) at (1, 1.25) {$1$};
		\node [style=small back] (521) at (2, 1.25) {};
		\node [style=none] (522) at (2, 1.25) {$2$};
		\node [style=small back] (523) at (3, 1.25) {};
		\node [style=none] (524) at (3, 1.25) {$3$};
		\node [style=small back] (525) at (4, 1.25) {};
		\node [style=none] (526) at (4, 1.25) {$4$};
		\node [style=medium back] (528) at (5.25, 2.25) {};
		\node [style=none] (529) at (5.25, 2.25) {$d_1$};
		\node [style=medium back] (530) at (5.25, -2.25) {};
		\node [style=none] (531) at (5.25, -2.25) {$d_2$};
		\node [style=small back] (532) at (9.5, 1.25) {};
		\node [style=none] (533) at (9.5, 1.25) {$1$};
		\node [style=small back] (534) at (10.5, 1.25) {};
		\node [style=none] (535) at (10.5, 1.25) {$2$};
		\node [style=small back] (536) at (14.5, 1.25) {};
		\node [style=none] (537) at (14.5, 1.25) {$3$};
		\node [style=small back] (538) at (15.5, 1.25) {};
		\node [style=none] (539) at (15.5, 1.25) {$4$};
		\node [style=medium back] (540) at (16.5, 2.25) {};
		\node [style=none] (541) at (16.5, 2.25) {$d_1$};
		\node [style=medium back] (542) at (16.5, -2.25) {};
		\node [style=none] (543) at (16.5, -2.25) {$d_2$};
		\node [style=medium back] (544) at (12.5, 0) {};
		\node [style=none] (545) at (12.5, 0) {$a_1$};
		\node [style=medium back] (546) at (19, 0) {};
		\node [style=none] (547) at (19, 0) {$a_2$};
		\node [style=small back] (548) at (23, 1.25) {};
		\node [style=none] (549) at (23, 1.25) {$1$};
		\node [style=small back] (550) at (24, 1.25) {};
		\node [style=none] (551) at (24, 1.25) {$2$};
		\node [style=small back] (552) at (28, 1.25) {};
		\node [style=none] (553) at (28, 1.25) {$3$};
		\node [style=small back] (554) at (29, 1.25) {};
		\node [style=none] (555) at (29, 1.25) {$4$};
		\node [style=medium back] (556) at (26, 0) {};
		\node [style=none] (557) at (26, 0) {$a_1$};
		\node [style=medium back] (558) at (32.5, 0) {};
		\node [style=none] (559) at (32.5, 0) {$a_2$};
	\end{pgfonlayer}
	\begin{pgfonlayer}{edgelayer}
		\draw [style=c0] (428) to (423);
		\draw [style=c0] (423) to (425);
		\draw [style=c0] (425) to (424);
		\draw [style=c0] (424) to (426);
		\draw [style=c0] (426) to (423);
		\draw [style=c0] (423) to (427);
		\draw [style=c0] (427) to (424);
		\draw [style=c0] (424) to (428);
		\draw [style=c0, bend left=45] (423) to (429);
		\draw [style=c0, bend left=45] (429) to (424);
		\draw [style=c0] (428) to (425);
		\draw [style=c0] (425) to (426);
		\draw [style=c0] (426) to (427);
		\draw [style=c0] (445) to (440);
		\draw [style=c0] (440) to (442);
		\draw [style=c0] (442) to (441);
		\draw [style=c0] (441) to (443);
		\draw [style=c0] (443) to (440);
		\draw [style=c0] (440) to (444);
		\draw [style=c0] (444) to (441);
		\draw [style=c0] (441) to (445);
		\draw [style=c0, bend left=45] (440) to (446);
		\draw [style=c0, bend left=45] (446) to (441);
		\draw [style=c0] (445) to (442);
		\draw [style=c0] (443) to (444);
		\draw [style=c0] (461) to (456);
		\draw [style=c0] (456) to (458);
		\draw [style=c0] (458) to (457);
		\draw [style=c0] (457) to (461);
		\draw [style=c0] (461) to (458);
		\draw [style=c0] (459) to (460);
		\draw [style=c0] (464) to (459);
		\draw [style=c0] (460) to (464);
		\draw [style=c0] (459) to (465);
		\draw [style=c0] (465) to (460);
		\draw [style=c1] (456) to (470);
		\draw [style=c1] (470) to (464);
		\draw [style=c1] (464) to (471);
		\draw [style=c1] (471) to (456);
		\draw [style=c1] (456) to (466);
		\draw [style=c1] (466) to (467);
		\draw [style=c1] (467) to (464);
		\draw [style=c1] (471) to (466);
		\draw [style=c1] (466) to (470);
		\draw [style=c1] (470) to (467);
		\draw [style=c1] (467) to (471);
		\draw [style=c1] (457) to (473);
		\draw [style=c1] (473) to (465);
		\draw [style=c1] (465) to (472);
		\draw [style=c1] (472) to (457);
		\draw [style=c1] (457) to (468);
		\draw [style=c1] (468) to (472);
		\draw [style=c1] (472) to (469);
		\draw [style=c1] (469) to (473);
		\draw [style=c1] (473) to (468);
		\draw [style=c1] (468) to (469);
		\draw [style=c1] (469) to (465);
		\draw [style=c0, bend left=45] (464) to (462);
		\draw [style=c0, bend right=45] (465) to (462);
		\draw [style=c0] (470) to (472);
		\draw [style=c0, bend right=105, looseness=2.75] (473) to (471);
		\draw [style=c0] (489) to (484);
		\draw [style=c0] (484) to (486);
		\draw [style=c0] (486) to (485);
		\draw [style=c0] (485) to (489);
		\draw [style=c0] (489) to (486);
		\draw [style=c0] (487) to (488);
		\draw [style=c0] (492) to (487);
		\draw [style=c0] (488) to (492);
		\draw [style=c0] (487) to (493);
		\draw [style=c0] (493) to (488);
		\draw [style=c1] (484) to (494);
		\draw [style=c2] (494) to (495);
		\draw [style=c1] (495) to (492);
		\draw [style=c2] (494) to (498);
		\draw [style=c2] (495) to (499);
		\draw [style=c1] (485) to (496);
		\draw [style=c2] (500) to (497);
		\draw [style=c2] (501) to (496);
		\draw [style=c2] (496) to (497);
		\draw [style=c1] (497) to (493);
		\draw [style=c2] (498) to (500);
		\draw [style=c2, bend right=105, looseness=2.75] (501) to (499);
	\end{pgfonlayer}
\end{tikzpicture}}}
    \caption{A temporal sequence $\tseq{G}$ and a subgraph $J_3$ of $G_3$. The graph $G_2$ is obtained from $G_1$ by deleting the edge $p$ and obtained from $G_3$ by contracting the \bcol\ edges and deleting the edges $a_1$ and $a_2$. The graph $J_3$ is obtained from $G_3$ by deleting some of the \bcol\ edges. The edges of a cycle $o$ in $J_3$ have been coloured \rcol.}
    \label{fig:expansion_impossible}
\end{figure}
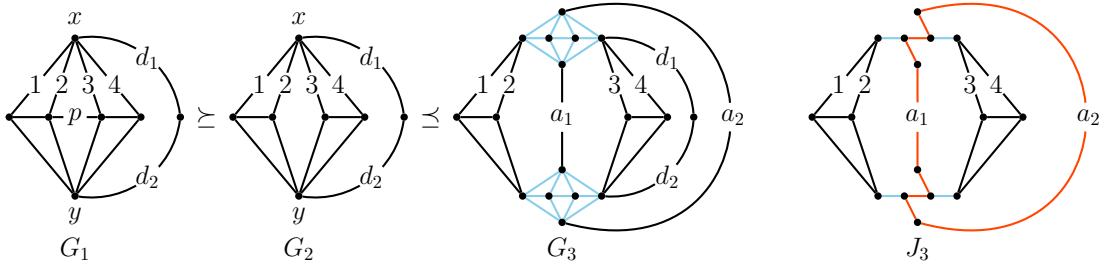

\begin{eg}\label{eg:expansion_impossible}
We shall explain how \autoref{fig:expansion_impossible} demonstrates that lists are necessary for our approach to work. 
Imagine we are given a temporal sequence $(G_1,G_2,G_3,..., G_n)$ starting with the temporal triple $(G_1,G_2,G_3)$ of \autoref{fig:expansion_impossible}. In our approach via \autoref{thm-list}, we would replace the graphs $G_1$ and $G_2$ by smaller graphs, which in this case would be cycles (if we chose $A=E(G_1)\sm\{d_1,d_2\}$), and replace $G_3$ by a larger graph $H_3$ that together with the list $L'$ captures all the embedding information of the temporal triple $(G_1,G_2,G_3)$. 
In fact in this example having nontrivial lists is necessary as the following fact shows.
\begin{prop}\label{eg-prop}
There is no planar graph $H_3\minorm G_3$ so that the embeddings $\phi_3$ that $H_3$ induces on $G_3$ are exactly those so that there are embeddings $\phi_1$ and $\phi_2$ so that $(\phi_1,\phi_2,\phi_3)$ is a simultaneous embedding of $(G_1,G_2,G_3)$.
\end{prop}

\begin{pf}
Suppose for a contradiction that there is such a graph $H_3$. 
Let $J_3$ be the subgraph of $G_3$ also depicted in \autoref{fig:expansion_impossible} and let $o$ be the cycle in $J_3$ containing the edges $a_1$ and $a_2$. Since $J_3$ has maximum degree three, there exists a subdivision $J_3'$ of $J_3$ as a subgraph of $H_3$. Let $o'$ be the subdivision of $o$ in $J_3'$. 
 
Since $o$ is a separating cycle of $G_3$, $o'$ must be a separating cycle of the planar graph $H_3$; and in any embedding of $H_3$ it must separate the edges labelled $1$ and $3$. 
We refer to the component of $H_3\sm o'$ containing $1$ as the \emph{positive side} and to the other component as the \emph{negative side}.
 
There are exactly two simultaneous embeddings of $\tseq{G}$, up to reorientation; one is depicted in \autoref{fig:expansion_impossible} and the other is obtained from the depicted one by reversing the rotators at $x$ and $y$ in $G_1$ and in $G_3$ reversing the rotators at every vertex where a \bcol\ and a \xcol\ edge meet. In particular, they induce either the rotator $1234d_1$ or $4321d_1$ at $x$ in $G_2$. 

Hence, there must exist embeddings $\iota_1$ and $\iota_2$ of $H_3$ so that the following holds: the embedding that $\iota_1$ induces on $J_3$ has the edges $1$ and $a_1$ in a shared face boundary, and $a_1$ and $3$ in another shared face boundary. The embedding that $\iota_2$ induces on $J_3$ has $2$ and $a_1$ in a shared face boundary, and $4$ and $a_1$ in another shared face boundary.
Thus the positive side and the negative side of $H_3$ each do not have unique embeddings inside, and thus each must contain a $2$-separator. However, then we can flip only the one on the positive side. This induces an embedding on $G_3$ which induces the rotator $2134d_1$ at $x$ in $G_2$, contrary to our assumptions. Thus there is no such graph $H_3$.
\end{pf}
\end{eg}

Given the machinery in the statement of \autoref{thm-list}, we would like a way to make a triple in our sequence $A$-mazing without changing its embeddability properties too much, or increasing its flexibility. To that end, we introduce the following definitions.

\begin{dfn}[Expansion]
    Given a temporal triple $\tseq{G}$, an \emph{expansion} $X_3\minorm G_3$ is a graph such that, for every simultaneous embedding $\tseq{\phi}$ of $\tseq{G}$, there exists an embedding $\psi_3$ of $H_3$ that induces $\phi_3$ on $G_3$.
    We shall refer to the triple $(G_1\minorm G_2\minorp H_3)$ as the \emph{expansion-triple}.
\end{dfn}

We say that an expansion $H_3\minorm G_3$ of a temporal triple $\tseq{G}$ is an \emph{accurate expansion} if $H_3$ is accurate for $G_3$.

\begin{eg}
    Let $\tseq{G}$ be a temporal triple and $H_3$ an expansion thereof. Clearly every simultaneous embedding of $(G_1\minorm G_2\minorp H_3)$ induces a simultaneous embedding of\\ $(G_1\minorm G_2\minorp G_3)$. However, we also have that, for every simultaneous embedding $\tseq{\phi}$ of $\tseq{G}$, there is an embedding $\psi_3$ of $H_3\minorm G_3$ so that $(\phi_1\minorm  \phi_2\minorp \psi_3)$ is a simultaneous embedding of $(G_1\minorm  G_2\minorp H_3)$. In this sense, $\tseq{G}$ and $(G_1\minorm G_2\minorp H_3)$ have the same embeddability properties.
\end{eg}

\begin{lem}\label{lem:amazing_expansion_main}
    Let $\tseq{G}$ be a temporal triple of planar graphs and $A$ a Tutte leaf of~$G_1$. Suppose that $G_3$ contains no upleaves, that $G_3$ is nice and that $\tseq{G}$ is not an amazingly-bad obstruction. Then there exists an nice accurate planar expansion $H_3$ of $\tseq{G}$ so that the expansion triple $(G_1\minorm G_2\minorp H_3)$ is $A$-mazing. Furthermore we have that\vspace{-.25cm}
    \[
    \Cbb(X_{i+2})\leq \Cbb(G_{i+2})\quad\text{and}\quad |E(X_{i+2})|\leq |E(G_{i+2})| + 4\;\Cbb(G_{i+2}).\vspace{-.25cm}
    \]
\end{lem}
The proof of \autoref{lem:amazing_expansion_main} will follow in \autoref{sec:amazing_expansions}.

\begin{pf}[Proof of \myref{lem:A} assuming \autoref{thm-list} and \autoref{lem:amazing_expansion_main}]
Let $\Gcal=(G_1\minorm\ldots\minorp G_n)$ be a nice temporal sequence of $2$-connected graphs along with a list $L$ and suppose that $A$ is a last upleaf rooted at some $i\in[n]$. Suppose further that some edge of $A$ appears in $G_{i+1}$. Then we are required to exhibit graph $H_{i+2}\minorm G_{i+2}$ and a list $L'$ for the temporal sequence\vspace{-.25cm}
    \[
    \Gcal':=(G_1\minorm\ldots\minorm G_{i-1}\minorp G_i[[A\ct:e]]\minorm G_{i+1}[[\pi A\ct:e]]\minorp H_{i+2}\minorm G_{i+3}\minorp\ldots\minorp G_n)\vspace{-.25cm}
    \]
    so that $(\Gcal',L')$ is an analogue of $(\Gcal,L)$ and furthermore, we have that\vspace{-.25cm}
    \[
    \Cbb(H_{i+2})\leq \Cbb(G_{i+2})\quad\text{and}\quad |E(H_{i+2})|\leq |E(G_{i+2})| + 4\;\Cbb(G_{i+2}).\vspace{-.25cm}
    \]
    If $A$ corresponds to a path in $G_i$, then our result follows immediately after setting $H_{i+2}=G_{i+2}$. So suppose that $A$ corresponds to a $3$-block.
    Now we apply \autoref{lem:amazing_expansion_main} to the temporal triple $(G_i\minorm G_{i+1}\minorp G_{i+2})$ to obtain an $A$-mazing expansion triple $(G_i\minorm G_{i+1}\minorp H_{i+2})$ so that $H_{i+2}$ is nice and accurate for $G_{i+2}$. Denote by $\Gcal^*$ the temporal sequence\vspace{-.25cm}
    \[
    (G_1\minorm\ldots\minorm G_{i+1}\minorp H_{i+2}\minorm G_{i+3}\minorp G_n).\vspace{-.25cm}
    \]
    Since both $H_{i+2}$ is accurate for $G_{i+2}$ and $H_{i+2}$ is a planar expansion, the list $L$ lifts via \autoref{lem:list_lifter} to a list $L^*$ for $\Gcal^*$ so that $(\Gcal,L)$ admits a simultaneous list-embedding if and only if $(\Gcal^*,L^*)$ admits a simultaneous list-embedding. Furthermore,\vspace{-.25cm}
    \[
    \Cbb(H_{i+2})\leq \Cbb(G_{i+2})\quad\text{and}\quad |E(H_{i+2})|\leq |E(G_{i+2})| + 6\;\Cbb(G_{i+2}).\vspace{-.25cm}
    \]
    Finally, applying \autoref{thm-list} to the pair $(\Gcal^*,L^*)$ yields our required result.
\end{pf}

\subsection{Amazing reductions}\label{sec:amazing_reductions}
In this subsection, we prove \autoref{thm-list}. In order to do so, we make the following preparations.

\begin{obs}\label{obs:greeable}
    Let $\tseq{G}$ be a temporal triple, along with an upleaf $A$ rooted at~$1$. Suppose that $A$ has a $3$-connected torso $G[[A:e]]$. Let $\alpha$ be an embedding of $G[[A:e]]$ and $\psi_3$ an embedding of $G_3$. If $\alpha \rightarrow G_2[[\pi A:e]]\leftarrow \psi_3$ holds, then the embedding $\psi_3$ is  $\alpha$-greeable.
\end{obs}
\begin{pf}[Proof:]
    follows immediately from the definition of $\alpha$-greeable.
\end{pf}

The converse of \autoref{obs:greeable} need not be true, as shown in the following example.

\begin{figure}
    \centering
    {\Large\scalebox{0.7}{\begin{tikzpicture}
	\begin{pgfonlayer}{nodelayer}
		\node [style=none] (341) at (-6, 3.75) {$x$};
		\node [style=none] (342) at (-6, -3.75) {$y$};
		\node [style=none] (357) at (-9, 2.5) {$A$};
		\node [style=black dot] (423) at (-6, 3) {};
		\node [style=black dot] (424) at (-6, -3) {};
		\node [style=black dot] (425) at (-6.5, 0) {};
		\node [style=black dot] (426) at (-4.75, 0) {};
		\node [style=black dot] (427) at (-3.25, 0) {};
		\node [style=black dot] (428) at (-8.75, 0) {};
		\node [style=black dot] (429) at (-1.5, 0) {};
		\node [style=none] (430) at (-0.5, 0) {$\succeq$};
		\node [style=none] (447) at (10.5, 0) {$\preceq$};
		\node [style=none] (476) at (-6, -4.75) {$G_1$};
		\node [style=black dot] (478) at (-10.5, 0) {};
		\node [style=none] (479) at (5, 3.75) {$x$};
		\node [style=none] (480) at (5, -3.75) {$y$};
		\node [style=black dot] (486) at (5, 3) {};
		\node [style=black dot] (487) at (5, -3) {};
		\node [style=black dot] (488) at (4.5, 0) {};
		\node [style=black dot] (489) at (6.25, 0) {};
		\node [style=black dot] (490) at (7.75, 0) {};
		\node [style=black dot] (491) at (2.5, 0) {};
		\node [style=black dot] (492) at (9.5, 0) {};
		\node [style=none] (494) at (5, -4.75) {$G_1$};
		\node [style=black dot] (495) at (0.5, 0) {};
		\node [style=none] (496) at (16, 3.75) {$x$};
		\node [style=none] (497) at (16, -3.75) {$y$};
		\node [style=black dot] (502) at (16, 3) {};
		\node [style=black dot] (503) at (16, -3) {};
		\node [style=black dot] (504) at (15.5, 0) {};
		\node [style=black dot] (505) at (17.25, 0) {};
		\node [style=black dot] (506) at (18.75, 0) {};
		\node [style=black dot] (507) at (13.5, 0) {};
		\node [style=black dot] (508) at (20.5, 0) {};
		\node [style=none] (509) at (16, -4.75) {$G_1$};
		\node [style=black dot] (510) at (11.5, 0) {};
		\node [style=small back] (511) at (-8.75, 1.25) {};
		\node [style=none] (512) at (-8.75, 1.25) {$1$};
		\node [style=small back] (513) at (-7.5, 1.25) {};
		\node [style=none] (514) at (-7.5, 1.25) {$2$};
		\node [style=small back] (515) at (-6.25, 1.25) {};
		\node [style=none] (516) at (-6.25, 1.25) {$3$};
		\node [style=small back] (517) at (-5.25, 1.25) {};
		\node [style=none] (518) at (-5.25, 1.25) {$4$};
		\node [style=small back] (519) at (-7.5, 0) {};
		\node [style=none] (520) at (-7.5, 0) {$e$};
		\node [style=small back] (529) at (13.25, 1.25) {};
		\node [style=none] (530) at (13.25, 1.25) {$3$};
		\node [style=small back] (531) at (14.5, 1.25) {};
		\node [style=none] (532) at (14.5, 1.25) {$4$};
		\node [style=small back] (533) at (15.75, 1.25) {};
		\node [style=none] (534) at (15.75, 1.25) {$1$};
		\node [style=small back] (535) at (16.75, 1.25) {};
		\node [style=none] (536) at (16.75, 1.25) {$2$};
		\node [style=small back] (537) at (2.25, 1.25) {};
		\node [style=none] (538) at (2.25, 1.25) {$1$};
		\node [style=small back] (539) at (3.5, 1.25) {};
		\node [style=none] (540) at (3.5, 1.25) {$2$};
		\node [style=small back] (541) at (4.75, 1.25) {};
		\node [style=none] (542) at (4.75, 1.25) {$3$};
		\node [style=small back] (543) at (5.75, 1.25) {};
		\node [style=none] (544) at (5.75, 1.25) {$4$};
	\end{pgfonlayer}
	\begin{pgfonlayer}{edgelayer}
		\draw [style=c1] (428) to (423);
		\draw [style=c1] (423) to (425);
		\draw [style=c1] (425) to (424);
		\draw [style=c1] (424) to (426);
		\draw [style=c1] (426) to (423);
		\draw [style=c0] (423) to (427);
		\draw [style=c0] (427) to (424);
		\draw [style=c1] (424) to (428);
		\draw [style=c1] (425) to (426);
		\draw [style=c0] (429) to (423);
		\draw [style=c0] (429) to (427);
		\draw [style=c0] (424) to (429);
		\draw [style=c1] (424) to (478);
		\draw [style=c1] (478) to (423);
		\draw [style=c1] (428) to (478);
		\draw [style=c1] (428) to (425);
		\draw [style=c1] (491) to (486);
		\draw [style=c1] (486) to (488);
		\draw [style=c1] (488) to (487);
		\draw [style=c1] (487) to (489);
		\draw [style=c1] (489) to (486);
		\draw [style=c0] (486) to (490);
		\draw [style=c0] (490) to (487);
		\draw [style=c1] (487) to (491);
		\draw [style=c1] (488) to (489);
		\draw [style=c0] (492) to (486);
		\draw [style=c0] (492) to (490);
		\draw [style=c0] (487) to (492);
		\draw [style=c1] (487) to (495);
		\draw [style=c1] (495) to (486);
		\draw [style=c1] (491) to (495);
		\draw [style=c1] (507) to (502);
		\draw [style=c1] (502) to (504);
		\draw [style=c1] (504) to (503);
		\draw [style=c1] (503) to (505);
		\draw [style=c1] (505) to (502);
		\draw [style=c0] (502) to (506);
		\draw [style=c0] (506) to (503);
		\draw [style=c1] (503) to (507);
		\draw [style=c1] (504) to (505);
		\draw [style=c0] (508) to (502);
		\draw [style=c0] (508) to (506);
		\draw [style=c0] (503) to (508);
		\draw [style=c1] (503) to (510);
		\draw [style=c1] (510) to (502);
		\draw [style=c1] (507) to (510);
	\end{pgfonlayer}
\end{tikzpicture}}}
    \caption{A temporal triple $\tseq{G}$ along with a sequence of embeddings which is not a simultaneous embedding. The non-torso edges $A$ of a $3$-block of $G_1$ have been coloured \bcol. The graph $G_2$ is identical to the graph $G_3$ (although its illustrated embedding differs).}
    \label{fig:agreeable_not_amazing}
\end{figure}
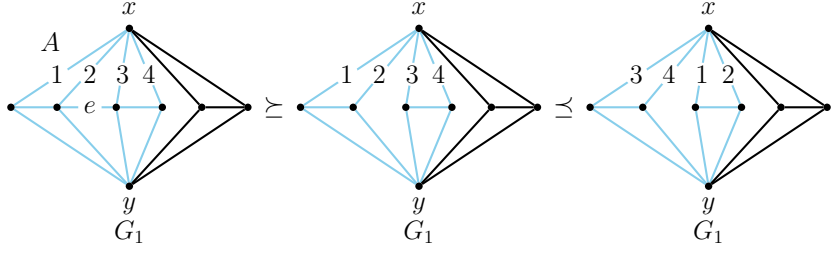

\begin{eg}\label{eg:agreeable_but}
    Consider the temporal triple $\tseq{G}$ depicted in \autoref{fig:agreeable_not_amazing}.  Let $\alpha$ denote the depicted embedding for the $3$-block corresponding to the \bcol\ edges of $G_1$. Then the depicted embedding of $G_3$ is $\alpha$-greeable. Indeed, we observe that the orientation of the $3$-blocks of $G_1$ and $G_4$ `agree' on  $G_2$. However, the depicted embeddings do not form a simultaneous embedding, even when restricted to the temporal triple $(G_1[[A:e]]\minorm G_2[[\pi A:e]]\minorp G_3)$. This is because the cyclic order of the edges $1,2,3,4$ in the rotator at $x$ in $G_1$ is not compatible with their cyclic order in the rotator at $x$ in $G_3$.
\end{eg}

Despite no converse to \autoref{obs:greeable} in general, \autoref{amazing_lock} provides a converse in the case where $\tseq{G}$ is $A$-mazing. It is for this reason that we defined $A$-mazing and used it in the statement of \autoref{thm-list}.

\begin{lem}\label{amazing_lock}
Let $\tseq{G}$ be a temporal triple, along with an upleaf $A$ rooted at $1$. Suppose that $A$ has a $3$-connected torso $G[[A:e]]$. Let $\alpha$ be an embedding of $G_1[[A:e]]$, and let $\psi_3$ be an embedding of $G_3$ that is $\alpha$-greeable. If $\tseq{G}$ is $A$-mazing, then $\alpha \rightarrow G_2[[\pi A:e]]\leftarrow \psi_3$. 
\end{lem}

In order to prove \autoref{amazing_lock}, we prepare as follows.

\begin{lem}\label{amazing-prop}
    Let $G$ be a $2$-connected graph and $x,y,z\in E(G)$ distinct.
    Then either there is a totally nested $2$-separator $s$ of $G$ so that all of $x$, $y$ and $z$ are in different $2$-components at $s$ or else there is a cyclic or $3$-connected torso\footnote{Recall that a torso consists of edges from its bag in the Tutte decomposition as well as `torso edges' which are from elsewhere. In this lemma, the edges $x$, $y$ and $z$ may be either of these types.} of the Tutte decomposition of $G$ containing $x$, $y$ and $z$.
\end{lem}
\begin{pf}
Consider the decomposition tree $T$ of the Tutte decomposition of $G$. Denote by $u_x$, $u_y$ and $u_z$ nodes of $T$ corresponding to the bags containing the edges $x$, $y$ and $z$, respectively. We emphasise that $u_x$, $u_y$ and $u_z$ may not all be distinct. 
Let $r$ be the endvertex of the path in $T$ from $u_z$ to the path\footnote{Given a tree $T$ with notes $u$ and $r$, by $uTr$ we denote the unique path in $T$ from $u$ to $r$.} $u_xTu_y$. 
This choice of $r$ ensures that the paths $u_xTr$, $u_yTr$ and $u_zTr$ are disjoint except at $r$. By our choice of $r$, the edges $x,y,z$ that are not contained in the bag corresponding to $r$ are in distinct branches. Let $W$ be the torso corresponding to $r$, where we label a torso edge by $x$, $y$ or $z$ if it corresponds to the branch containing $x$, $y$ or $z$, respectively. It follows from $2$-connectivity that such a torso with the given labelling exists. If $W$ is a parallel graph, then $W$ corresponds to a $2$-separator $s$ of $G$ so that all of $x$, $y$ and $z$ are in different 2-components at $s$. Otherwise, the torso is cyclic or $3$-connected, as desired.
\end{pf}

\begin{lem}[Amazing Lifting Lemma]\label{amazing_lift}
Let $\tseq{G}$ be a temporal triple, along with an upleaf $A$ rooted at $1$. Assume that $\tseq{G}$ is $A$-mazing and let $s$ be a totally nested $2$-separator of $G_2[[\pi A:e]]$ and $x$, $y$ and $z$
be edges in different $2$-components at $s$. There is a torso $b$ of a $3$-block of $G_3$ containing all of $x$, $y$ and $z$.
\end{lem}

In the proof of \autoref{amazing_lift} we shall use the following lemma.

\begin{lem}\label{overlap-cross}
Let $G$ be a $2$-connected graph with three subgraphs $A_1$, $A_2$ and $A_3$ so that $A_1\cup A_2\cup A_3=G$.
Assume there are distinct vertices $u_1$, $u_2$ and $u_3$ so that $A_i\cap A_j=\{u_k\}$ for $\{i,j,k\}=\{1,2,3\}$.
Then no bond contains edges of all three subgraphs $A_1$, $A_2$ and $A_3$.
\end{lem}

\begin{pf}
Let an arbitrary bond of $G$ be given; denote it by $b$.
Two vertices of $u_1$, $u_2$ and $u_3$ are on the same side of the bond $b$, say these are $u_1$ and $u_2$.
Denote the side of $b$ that does not contain $u_1$ by $X$. 
Consider the (possibly improper) separation $(V(A_3), V(A_1\cup A_2))$. Its separator is $\{u_1,u_2\}$, and hence it is a 2-separation. 

Suppose for a contradiction that each subgraph $A_i$ for $i\in \{1,2,3\}$ contains an edge of $b$. 
Since every edge in $b$ has an endvertex in $X$, every subgraph $A_i$ for $i\in \{1,2,3\}$ contains a vertex of $X$.
Note that since $G$ is $2$-connected, sides of bonds are connected; in particular $X$ is connected. 
So there is a path contained in the connected set $X$ from a vertex of $A_3$ to a vertex of $A_1$; pick such a path and denote it by $P$.
The separation $(V(A_3), V(A_1\cup A_2))$ contains an endvertex of $P$ in each of its two sides. Thus its separator $\{u_1,u_2\}$ contains a vertex of $P$. Thus $u_1$ or $u_2$ is in $X$, a contradiction to the choice of $X$. Hence 
no bond $b$ contains edges of all three sets $A_1$, $A_2$ and $A_3$.
\end{pf}

\begin{pf}[Proof of \autoref{amazing_lift}]
Since $x$, $y$ and $z$ are edges of $G_2$, they are also edges of $G_3$. 
By \autoref{amazing-prop} there is a cyclic or $3$-connected torso of the Tutte decomposition of $G_3$ containing all three of $x,y,z$ or else there is a totally nested $2$-separator $s'$ of $G_3$ so that $x,y,z$ are in different 2-components at $s'$. The latter is not possible as at least two of $x,y,z$ are in $A$ and $\tseq{G}$ is $A$-mazing. So there is a cyclic or $3$-connected torso $W$ of the Tutte decomposition of $G_3$ containing all three of $x,y,z$.
\begin{claim}\label{nice-prop}
The torso $W$ is not a cycle.
\end{claim}

\begin{cproof}
Suppose for a contradiction that the torso $W$ is a cycle.
Pick vertices $u_1$, $u_2$ and $u_3$ on the cycle $W$ so that each of the segments of $W$ between two vertices $u_i$ contains exactly one of the three edges $x$, $y$ and $z$. By definition of torsos, the vertices of $W$ are vertices of $G_3$; and we shall sometimes interpret the vertices $u_i$ as vertices of $G_3$. 
Pick $i,j,k$ so that $\{i,j,k\}=\{1,2,3\}$. 
Since $W$ is a cyclic torso of the Tutte decomposition, each pair $\{u_j,u_k\}$ is a 2-separator of $G_3$. 
For each $i\in [3]$ we set $A_i$ to be the set  $\{u_j,u_k\}$ together with all components of $G_3- \{u_j,u_k\}$ that do not contain the vertex $u_i$. 

Observe that each $A_i$ contains exactly one of the edges $x,y,z$.
We claim that $A_1\cup A_2\cup A_3=G_3$. For this we show that every vertex and edge of $G_3$ is in one of the sets $A_i$ for $i\in [3]$.
Clearly the vertices $u_1$, $u_2$ and $u_3$ and the edges between two of these vertices are in some set $A_i$.
Each component of $G_3-u_1-u_2-u_3$ has only two neighbours in $\{u_1,u_2,u_3\}$, and thus occurs in some set $A_i$.
This implies that every edge that has an endvertex not of the form $u_i$ is in a set $A_j$. 
Thus we have verified the assumptions of \autoref{overlap-cross}. Hence by that lemma no bond contains edges of all three sets $A_1$, $A_2$ and $A_3$. Since each of the edges $x$, $y$ and $z$ in the torso $W$ is in a different segment between vertices $u_i$, in the graph $G_3$ each of the edges $x$, $y$ and $z$ is in a different set $A_i$. Hence no bond of $G_3$ contains all three edges $x$, $y$ and $z$. 

By assumption, the graph $G_2$ has a theta-graph as a minor whose three edges have the labels $x$, $y$ and $z$.
Hence also $G_3$ has such a minor. Since the theta-graph has a bond containing $x$, $y$ and $z$, also $G_3$ has a bond containing the edges $x$, $y$ and $z$. This is in contradiction to the conclusion of the previous paragraph. Hence $W$ cannot be a cycle. 
\end{cproof}

Since the torso $W$  is not a cycle by \autoref{nice-prop}, it is $3$-connected by the paragraph just above that claim.
Hence $W$ is a $3$-block that contains the edges $x$,  $y$ and $z$; this completes our proof.
\end{pf}

\begin{pf}[Proof of \autoref{amazing_lock}]
Let $\tseq{G}$ be a temporal triple, along with an upleaf $A$ rooted at $1$. Suppose that $A$ has a $3$-connected torso $G[[A:e]]$. Let $\alpha$ be an embedding of $G_1[[A:e]]$, and let $\psi_3$ be an embedding of $G_3$ that is $\alpha$-greeable. Assume that $\tseq{G}$ is $A$-mazing. We are required to show that $\alpha \rightarrow G_2[[\pi A:e]]\leftarrow \psi_3$ holds.

Let $\iota$ be the embedding of $G_2[[\pi A:e]]$ induced by the embedding $\alpha$. 
Similarly, let $\dot \iota$ be the embedding of $G_2[[\pi A:e]]$ induced by the embedding $\psi_3$. Our goal is to show that $\iota=\dot\iota$. Suppose for a contradiction that this doesn't hold. 
Then by \autoref{determine_embed}, either there is $3$-block $b_2$ of $G_2[[\pi A:e]]$
so that the embeddings of $b_2$ induced by $\iota$ and $\dot\iota$ do not coincide, or there is a totally nested $2$-separator $s$
of $G_2[[\pi A:e]]$, and three $2$-components $c_1$, $c_2$ and $c_3$ at $s$ receiving distinct orders in the cyclic orders induced on the $2$-components at $s$ by $\iota$ and $\dot\iota$.

Assume we are in the first case. Since $G_2\minorp G_3$, there is a $3$-block $b$ of $G_3$ so that $b_2$ is a minor of $b$. Since $b_2$ is $3$-connected, it has a $K_4$-minor by a theorem of Tutte, so it has a theta graph $\triple$ as a minor.

Each embedding of $b$ induces an embedding of $\triple$ and this defines a bijection between the two embeddings of $b$ and the two embeddings of $\triple$. Indeed, at least one embedding of $\triple$ can be obtained since $b$ is planar and hence both can, since one is the reorientation of the other.
From this we get that the embeddings $\iota$ and $\dot\iota$ induce different embeddings on $\triple$, which is a contradiction to the assumption that $\psi_3$ is $\alpha$-greeable. In particular, the required property is not satisfied on the temporal triple $(G_1[[A:e]]\minorm \triple \minorp b)$.

Having derived a contradiction in the first case, we now consider the second case. 
Take a theta-minor of $G_2[[\pi A:e]]$ that contains an edge from each of the components  $c_1$, $c_2$ and $c_3$; denote this theta graph by $\triple$.
By \myref{amazing_lift}, there is a $3$-block $b$ of $G_3$ that minors down to $\triple$. 
The embeddings $\iota$ and $\dot\iota$ induce different embeddings on $\triple$. This is a contradiction to the assumption that $\psi_3$ is $\alpha$-greeable. In particular we have that $\alpha\rightarrow\triple\leftarrow \psi_3\restric b$ does not hold.
\end{pf}

\begin{lem}\label{H-red}
Let $\tseq{G}$ be an $A$-mazing temporal triple for a Tutte leaf $A$ of $G_1$. Suppose that $A$ has a $3$-connected torso $G_1[[A:e]]$ and let $\alpha$ be an embedding of $G_1[[A:e]]$. Let $(\psi_1,\psi_2,\psi_3)$ be a simultaneous embedding of $(G_1[[A\ct:e]]\minorm G_2[[\pi A\ct:e]] \minorp G_3)$. If $\psi_3$ is $\alpha$-greeable, then there is a simultaneous embedding $(\phi_1,\phi_2,\psi_3)$ of $(G_1\minorm G_2 \minorp G_3)$ so that $\phi_i$ induces $\psi_i$ for $i=1,2$ and $\phi_1$ induces $\alpha$ on~$G_1[[A:e]]$.
\end{lem}

\begin{pf}
    The graph $G_1$ is a $2$-sum of the graphs $G_1[[A:e]]$ and $G_1[[A\ct:e]]$ along the gluing edge $e$. Similarly the graph $G_2$ is a $2$-sum of the graphs $G_2[[\pi A:e]]$ and $G_2[[\pi A\ct:e]]$ along the gluing edge $e$. Denote by $\phi_1$ the embedding obtained by applying \autoref{lem:2-sum_embeddings} to $\alpha$ and $\psi_1$. Similarly let $\phi_2$ denote the embedding obtained by applying \autoref{lem:2-sum_embeddings} to $\alpha\restric G_2[[\pi A:e]]$ and $\psi_2$. By construction, $\phi_1$ induces $\phi_2$ on $G_2$.
    
    By \autoref{amazing_lock}, we have that $\alpha \rightarrow G_2[[\pi A:e]]\leftarrow \psi_3$. Hence $\phi_2$ induces the same embedding as $\psi_3$ on both $G_2[[\pi A:e]]$ and $G_2[[\pi A\ct:e]]$. Hence we have that $\psi_3$ induces $\phi_2$ on $G_2$ and this completes the proof that $(\phi_1,\phi_2,\psi_3)$ is a simultaneous embedding.
\end{pf}

In the following, we consider two $3$-blocks to be the same if they correspond to the same bag of the Tutte decomposition. Recall that in a temporal triple $\tseq{G}$, we say that a $3$-block $b_1$ of $G_1$ is linked to a $3$-block $b_3$ of $b_3$ if $\tseq{G}$ contains a minor of the form $(b_1\minorm\triple\minorp b_3)$.
\begin{lem}\label{extend-list-3}
    Let $\tseq{G}$ be a temporal triple and $A$ a Tutte leaf of $G_1$. Let $b:=G_i[[A:e]]$ denote its $3$-connected torso and let $\alpha$ be an embedding of $b$. Let $L$ be a list for $\tseq{G}$ that does not form an inconsistent obstruction. If $\tseq{G}$ is not a disagreeable obstruction, then there exists a list $L^+$ of $\tseq{G}$ so that every simultaneous embedding of $\tseq{G}$ that respects $L$ also respects $L^+$ and, for every embedding $\psi$ of $G_3$ which respects $L^+$, we have that either $\psi$ or its reorientation is $\alpha$-greeable. Furthermore, the part of $L^+$ containing $b$ contains all of the $3$-blocks of $G_3$ which are linked to $b$.
\end{lem}
\begin{pf}
    By reorienting a part of $L$ if necessary, we assume that the embedding $\alpha$ is in $L$. Since $\tseq{G}$ is not a disagreeable obstruction, for each $3$-block $w$ in $G_3$ which is linked to $b$, there is an embedding $\alpha_w$ of $w$ so that $\alpha\rightarrow\triple\leftarrow\alpha_w$ holds for each minor of $\tseq{G}$ of the form $(b\minorm\triple\minorp w)$.
    \begin{claim}\label{claim:no-inc-1}
        For each part $l$ of the list $L$ which does not contain $b$, either:
        \begin{enumerate}\vspace{-.15cm}
            \item for every $3$-block $w$ in $l$ from $G_3$ and linked to $b$, the embedding $\alpha_w$ appears in $L$; or\vspace{-.15cm}
            \item for every $3$-block $w$ in $l$ from $G_3$ and linked to $b$, the embedding $(\alpha_w)^-$ appears in~$L$.\vspace{-.15cm}
        \end{enumerate}
    \end{claim}
    \begin{cproof}
        Suppose for a contradiction that there are $3$-blocks $u$ and $w$ of $G_3$ which are linked to $b$, in the same part of $L$, so that $L$ assigns to $u$ the embedding $\alpha_u$ and to $w$ the embedding $(\alpha_w)^-$. Then $L$ forms an inconsistent obstruction.
    \end{cproof}

    \begin{claim}\label{claim:no-inc-2}
        Let $l$ be the part of the list $L$ which contains $b$. Then for every $3$-block $w$ of $G_3$ in $l$ which is linked to $b$, we have that $L$ specifies embedding $\alpha_w$ for $w$.
    \end{claim}
    \begin{cproof}
        Suppose for a contradiction that this is not the case. Then again we have an inconsistent obstruction.
    \end{cproof}
    
    Given \autoref{claim:no-inc-1}, let $L'$ be the equivalent list obtained from $L$ by reorienting all of the embeddings in the parts of $L$ which assign to some $3$-block $w$ of $G_3$ which is linked to $b$ the embedding $(\alpha_w)^-$. Along with \autoref{claim:no-inc-2}, we observe that for every $3$-block $w$ of $G_3$ which is linked to $b$ is assigned $\alpha_w$ by $L'$.

    Let $M$ be the list on $\tseq{G}$ which has one non-trivial part $m$ consisting of the linked $3$-blocks of $G_3$ along with the $3$-block $b$, which assigns to each $w\in m$ the embedding $\alpha_w$. By the above analysis, we can choose $M$ so that it specifies identical embeddings to $L'$. Hence $M$ and $L'$ are combinable and we let $L^+$ be their combination.

    \begin{claim}
        Every embedding $\psi$ of $G_3$ that satisfies the list $M$ is $\alpha$-greeable.
    \end{claim}
    \begin{cproof}
        This follows from the construction of the list $M$ and the definition of our embeddings $\alpha_w$.
    \end{cproof}
    
    By the above claim, we have that every simultaneous embedding of $\tseq{G}$ satisfies $M$. By \autoref{lem:list-combining_machine}, we get that a simultaneous embedding of $\tseq{G}$ satisfies $L$ if and only if it satisfies $L^+$.
\end{pf}

\begin{cor}\label{extend-list}
    The result \autoref{extend-list-3} extends to temporal sequences of arbitrary length with the graphs $G_1,G_2,G_3$ replaced by a consecutive triple $G_i,G_{i+1},G_{i+2}$.
\end{cor}
\begin{pf}
    The proof of \autoref{extend-list-3} easily extends to temporal sequences of arbitrary length. Indeed, we work in the triple $(G_i\minorm G_{i+1}\minorp G_{i+2})$ except when reorienting parts of our lists, in which case we perform this operation globally.
\end{pf}

We summarise this section as follows.
\begin{pf}[Proof of \autoref{thm-list}]
    Let $\Gcal=(G_1\minorm\ldots\minorp G_n)$ be a temporal sequence with list $L$, and $A$ an upleaf rooted at some $i\in[n]$ that corresponds to a $3$-block. Suppose that the temporal triple $(G_i\minorm G_{i+1}\minorp G_{i+2})$ is $A$-mazing. Let $b:=G_i[[A:e]]$ be the $3$-connected torso for $A$. Assume that $\Gcal$ does not contain a disagreeable obstruction, and along with $L$ is not an inconsistent obstruction.
    For brevity, we write\vspace{-.25cm}
    \begin{equation*}
        \Hcal=(H_1\minorm\ldots\minorp H_n)\quad:=\quad(G_1\minorm\ldots\minorp G_i[[A\ct:e]]\minorm G_{i+1}[[\pi A\ct:e]]\minorp G_{i+2}\minorp\ldots\minorp G_n).\vspace{-.25cm}
    \end{equation*}
    We are required to show that there is a list $L'$ of the temporal sequence $\Hcal$ so that $(\Hcal,L')$ admits a simultaneous list-embedding if and only if $(\Gcal,L)$ admits a simultaneous list-embedding. 
    
Let $\alpha$ be one of the embeddings of the $3$-block $b$. By applying \autoref{extend-list} if necessary, we may assume that any embedding $\psi$ of $G_{i+2}$ satisfying $L$ is such that $\psi$ or its reorientation $(\psi)^-$ is $\alpha$-greeable and furthermore the part of $L$ containing $b$ also contains all of the $3$-blocks of $G_{i+2}$ which are linked to $b$. By reorienting the embeddings in the part of $L$ containing $b$ if necessary, we assume that $L$ specifies the embedding $\alpha$ for $b$.

Observe that the $3$-blocks of $\Hcal$ are exactly the $3$-blocks of $\Gcal$ without $b$. Let $L'$ be the restriction of $L$ to the $3$-blocks of $\Hcal$. It is immediate that any simultaneous list-embedding of $(\Gcal,L)$ induces a simultaneous list-embedding of $(\Hcal,L')$. So let a simultaneous list-embedding $(\psi_1,\ldots,\psi_n)$ of $(\Hcal,L')$ be given. There are two cases.

\textbf{Case 1.} The graph $G_{i+2}$ contains no $3$-block which is linked to $b$. Then immediately from the definition of $\alpha$-greeable, we get that both of the embeddings $\psi_{i+2}$ and its reorientation $(\psi_{i+2})^-$ are $\alpha$-greeable. By reorienting it if necessary we may assume that the simultaneous embedding $(\psi_1,\ldots,\psi_n)$ is such that the following holds: for each $3$-block in the same part of $L$ as $b$, the embedding induced on $w$ is precisely the embedding specified by $L$, and $\psi_{i+2}$ is $\alpha$-greeable. Now by \autoref{H-red}, there is a simultaneous embedding $(\phi_i,\phi_{i+1},\psi_{i+2})$ of $(G_i \minorm G_{i+1} \minorp G_{i+2})$
that induces $(\psi_i,\psi_{i+1},\psi_{i+2})$ and is such that $\phi_i$ induces $\alpha$ on $b$.
Consider $(\psi_1,\ldots,\psi_{i-1},\phi_i,\phi_{i+1},\psi_{i+2}, \ldots,\psi_n)$. Since $A$ is an upleaf, this is a simultaneous embedding of $\Gcal$. Furthermore, since every $3$-block $w$ in the same part as $b$ has the embedding specified by $L$ and so does $b$, we have that $(\psi_1,\ldots,\psi_n)$ satisfies $L$. In this case, we are done.

\textbf{Case 2.} Not case 1. Thus, there exists a minor of $(G_i\minorm G_{i+1}\minorp G_{i+2})$ of the form $(b\minorm\triple\minorp w)$ where $w$ is a $3$-block of $G_{i+2}$ and so that $\alpha\rightarrow\triple\leftarrow\psi_{i+2}\restric w$. By a previous assumption, we have that $\psi_{i+2}$ or its reorientation is $\alpha$-greeable. Suppose for a contradiction that the reorientation $(\psi_{i+2})^-$ is $\alpha$-greeable. This implies $\alpha\rightarrow\triple\leftarrow (\psi_{i+2})^-\restric w=(\psi_{i+2}\restric w)^-$, contradicting our assumption.
So we get that $\psi_{i+2}$ is $\alpha$-greeable.

Now by \autoref{H-red}, there is a simultaneous embedding $(\phi_i,\phi_{i+1},\psi_{i+2})$ of $(G_i \minorm G_{i+1} \minorp G_{i+2})$
that induces $(\psi_i,\psi_{i+1},\psi_{i+2})$ and is such that $\phi_i$ induces $\alpha$ on $b$.
Consider $(\psi_1,\ldots,\psi_{i-1},\phi_i,\phi_{i+1},\allowbreak\psi_{i+2},\ldots,\psi_n)$. Since $A$ is an upleaf, this is a simultaneous embedding of $\Gcal$. All that remains is to show that this simultaneous embedding satisfies the list $L$. Clearly, since $(\psi_1,\ldots,\psi_{i-1},\phi_i,\phi_{i+1},\allowbreak\psi_{i+2}, \ldots,\psi_n)$ induces $(\psi_1,\ldots,\psi_n)$, we have that all of the parts except perhaps the one containing $b$, say $l$, are satisfied. We have to check that, for every $3$-block $u$ in $l$, the embedding given by $L$ for $u$ is satisfied. Let $\alpha_w$ and $\alpha_u$ be the embeddings of $w$ and $u$ specified by $L$, respectively. We have that $\alpha_w$ is induced by $\psi_{i+2}$ since it is $\alpha$-greeable. But then since our simultaneous embedding induces $(\psi_1,\ldots,\psi_n)$, which satisfies $L'$, we have that our simultaneous embedding also induces $\alpha_u$. This completes this case, and our proof.
\end{pf}

\subsection{Becoming amazing}\label{sec:amazing_expansions}

In \autoref{sec:amazing_reductions} we proved \autoref{H-red} and \autoref{thm-list}, which make use of the assumption that $\tseq{G}$ is $A$-mazing. In this subsection, we prove \autoref{lem:amazing_expansion_main}, which lets us build an $A$-mazing graph from one which may not be $A$-mazing.

Recall that an upleaf of a graph is a last upleaf if no graph with higher index in the temporal sequence has an upleaf rooted at it. For a tuplet $t$ of a $2$-connected graph $G$, let $\Ccal(t)$ denote the set of $2$-components at $t$.

The following extends \autoref{lem:affected_G_2} and \autoref{lem:affected_G_2_cycle} to full temporal triples (not just pairs of graphs related via the minor relation).

\begin{lem}\label{lem:A-ffected}
    Let $\tseq{G}$ be a $2$-connected temporal triple that is not an amazingly-bad obstruction and suppose $A$ is a last upleaf rooted at $1$. Let $t$ be a $2$-separator of $G_3$. Then in every simultaneous embedding of $\tseq{G}$, the $A$-ffected $2$-components at $t$ in $G_3$ form a segment in the cyclic structure induced on $\Ccal(t)$ by this simultaneous embedding. Furthermore, the (cyclic or linear) structure of this segment is independent of the choice of simultaneous embedding.
\end{lem}
\begin{pf}
    Since otherwise there would be nothing to prove, we assume that $A$ has nonempty shadow $\pi A$ in $G_2$ and that $t$ is a tuplet.
    \begin{claim}\label{claim:shadow_2-sep}
        The shadow $t'=\pi t$ of $t$ in $G_2$ is also a tuplet and, for each $2$-component $c'$ at $t'$, there is a unique $2$-component $c$ at $t$ which contains $c'$ as a minor. Furthermore, this assignment is surjective.
    \end{claim}
    \begin{cproof}
        Since $G_3$ has no upleaves, for each $2$-component $c$ at $t$, it's shadow $\pi c$ in $G_2$ contains at least one edge. Hence $t'$ is a separator of $G_2$ and since $G_2$ is $2$-connected, we conclude that $t'$ must be a tuplet in $G_2$. Furthermore, for each component $c'$ at $t'$, there is a unique $2$-component $c$ at $t$ which contains $c'$ as a minor. This assignment is surjective since the shadow of each $2$-component $t$ contains at least one edge.
    \end{cproof}

    By \autoref{claim:shadow_2-sep}, we get a surjective function $f:\Ccal(t')\rightarrow\Ccal(t)$ assigning each $2$-component $c'\in\Ccal(t')$ to its corresponding $2$-component $c\in\Ccal(t)$. Suppose that $(\phi_1,\phi_2,\phi_3)$ is a simultaneous embedding of $\tseq{G}$ and observe that $f$ is order preserving between the cyclic structure induced on $\Ccal(t')$ by $\phi_2$ and the cyclic structure induced on $\Ccal(t)$ by $\phi_3$. Further observe that $f$ preserves the property of being $A$-ffected.
    
    Now if $\pi A^C$ contains at most one edge, then we are be done, since $G_2$ is a minor of $G_1[[A:e]]$ for an appropriate choice of $e$, which is a minor of $G_1$ and has a unique embedding up to reorientation.
    So we assume that $\pi A^C$ contains at least two edges. Since $G_2$ is $2$-connected, the pair $(\pi A,\pi A^C)$ induces a nontrivial $2$-separation of $G_2$.
    Hence we may apply \autoref{lem:affected_G_2} and \autoref{lem:affected_G_2_cycle} to obtain that the $A$-ffected $2$-components form a segment in the cyclic structure induced on $\Ccal(t')$ by $\phi_2$, and the segments' (cyclic or linear) structure is independent of our choice of simultaneous embedding. That is, the so-called (cyclic or linear) structure suggested by $A$.
    
    Since $\tseq{G}$ is not an amazingly-bad obstruction, each $2$-component $c$ at $t$ in $G_3$ minors to a collection of $2$-components $x(c)$ at $t'$ in $G_2$ that is affectively connected. Restricting to the $A$-ffected $2$-components at $t'$ only, each $x(c)$ forms a segment in the (cyclic or linear) structure suggested by $A$. Hence the (cyclic or linear) structure suggested by $A$ lifts to a (cyclic or linear) structure on the $A$-ffected $2$-components at $t$ in $G_3$. Furthermore this (cyclic or linear) structure forms a segment in the cyclic structure induced on the $2$-components at $t$ in $G_3$ by any simultaneous embedding.
\end{pf}

In the context of \autoref{lem:A-ffected}, we call the unique (cyclic or linear) structure induced on the $A$-ffected $2$-components at $t$, the \emph{structure suggested by $A$}.

\begin{lem}\label{lem:tuplets_of_tree_coadd}
    Let $H$ be a $2$-connected graph and $G$ a graph obtained from $H$ by contracting a tree $T$ in $H$ to a single vertex $x$. Let $t$ be a tuplet of $H$. If $G$ is $2$-connected, then either $t$ is also a tuplet in $G$ or is of the form $\{j,y\}$ for some vertex $j$ in $T$ and some vertex $y$ disjoint from $T$. In the latter case, $\{x,y\}$ forms a tuplet of $G$.
\end{lem}
\begin{pf}
    Let $t$ be a tuplet of $H$. Since otherwise we would be done, assume that one of the vertices of $t$ is in $T$. First suppose that both of the vertices of $t$ lie in $T$. Then there is a $2$-component $c$ at $t$ containing a path in $T$ between the vertices of $t$. Hence every other $2$-component at $t$ contains at least one edge from outside of $T$, or else $T$ would contain a cycle. So let $c'$ and $c''$ be $2$-components at $t$ with nonempty shadow in $G=H/T$. But then the shadow of $c'$ is separated from the shadow of $c''$ by the vertex $x$, a contradiction to the fact that $G$ is $2$-connected. So exactly one vertex $j$ of $t$ lies in $T$ and the other, say $y$, is disjoint from $T$. Now each $2$-component at $t$ contains some edge of $G$ and hence has nonempty shadow in $G$. We conclude that $\{x,y\}$ is a tuplet of $G$.
\end{pf}

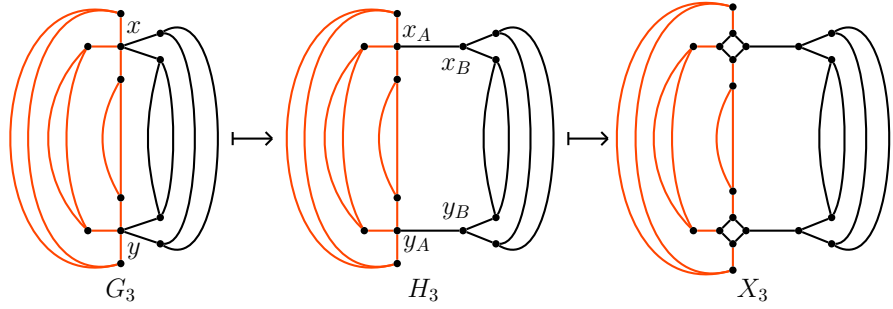
\begin{figure}
    \centering
    {\Large\scalebox{0.7}{\begin{tikzpicture}
	\begin{pgfonlayer}{nodelayer}
		\node [style=none] (341) at (6.5, 4.25) {$x$};
		\node [style=none] (342) at (6.5, -4.25) {$y$};
		\node [style=none] (476) at (6, -5.75) {$G_3$};
		\node [style=black dot] (479) at (6, 3.5) {};
		\node [style=black dot] (480) at (6, -3.5) {};
		\node [style=black dot] (481) at (6, -2.25) {};
		\node [style=black dot] (482) at (4.75, -3.5) {};
		\node [style=black dot] (483) at (6, -4.75) {};
		\node [style=black dot] (484) at (6, 4.75) {};
		\node [style=black dot] (485) at (4.75, 3.5) {};
		\node [style=black dot] (486) at (6, 2.25) {};
		\node [style=black dot] (487) at (7.5, 3) {};
		\node [style=black dot] (488) at (7.5, 4) {};
		\node [style=black dot] (489) at (7.5, -4) {};
		\node [style=black dot] (490) at (7.5, -3) {};
		\node [style=none] (491) at (10.25, 0.25) {};
		\node [style=none] (492) at (10.25, -0.25) {};
		\node [style=none] (493) at (10.25, 0) {};
		\node [style=none] (494) at (11.75, 0) {};
		\node [style=none] (495) at (11.5, 0.25) {};
		\node [style=none] (496) at (11.5, -0.25) {};
		\node [style=none] (497) at (17.25, 4) {$x_A$};
		\node [style=none] (498) at (17.25, -4) {$y_A$};
		\node [style=none] (499) at (17.5, -5.75) {$H_3$};
		\node [style=black dot] (500) at (16.5, 3.5) {};
		\node [style=black dot] (501) at (16.5, -3.5) {};
		\node [style=black dot] (502) at (16.5, -2.25) {};
		\node [style=black dot] (503) at (15.25, -3.5) {};
		\node [style=black dot] (504) at (16.5, -4.75) {};
		\node [style=black dot] (505) at (16.5, 4.75) {};
		\node [style=black dot] (506) at (15.25, 3.5) {};
		\node [style=black dot] (507) at (16.5, 2.25) {};
		\node [style=black dot] (508) at (20.25, 3) {};
		\node [style=black dot] (509) at (20.25, 4) {};
		\node [style=black dot] (510) at (20.25, -4) {};
		\node [style=black dot] (511) at (20.25, -3) {};
		\node [style=none] (512) at (23, 0.25) {};
		\node [style=none] (513) at (23, -0.25) {};
		\node [style=none] (514) at (23, 0) {};
		\node [style=none] (515) at (24.5, 0) {};
		\node [style=none] (516) at (24.25, 0.25) {};
		\node [style=none] (517) at (24.25, -0.25) {};
		\node [style=black dot] (518) at (19, 3.5) {};
		\node [style=black dot] (519) at (19, -3.5) {};
		\node [style=none] (520) at (18.75, 2.75) {$x_B$};
		\node [style=none] (521) at (18.75, -2.75) {$y_B$};
		\node [style=none] (524) at (30, -5.75) {$X_3$};
		\node [style=black dot] (527) at (29.25, -2) {};
		\node [style=black dot] (528) at (27.75, -3.5) {};
		\node [style=black dot] (529) at (29.25, -5) {};
		\node [style=black dot] (530) at (29.25, 5) {};
		\node [style=black dot] (531) at (27.75, 3.5) {};
		\node [style=black dot] (532) at (29.25, 2) {};
		\node [style=black dot] (533) at (33, 3) {};
		\node [style=black dot] (534) at (33, 4) {};
		\node [style=black dot] (535) at (33, -4) {};
		\node [style=black dot] (536) at (33, -3) {};
		\node [style=black dot] (537) at (31.75, 3.5) {};
		\node [style=black dot] (538) at (31.75, -3.5) {};
		\node [style=black dot] (539) at (29.75, 3.5) {};
		\node [style=black dot] (540) at (28.75, 3.5) {};
		\node [style=black dot] (541) at (29.25, 4) {};
		\node [style=black dot] (542) at (29.25, 3) {};
		\node [style=black dot] (543) at (28.75, -3.5) {};
		\node [style=black dot] (544) at (29.25, -3) {};
		\node [style=black dot] (545) at (29.75, -3.5) {};
		\node [style=black dot] (546) at (29.25, -4) {};
	\end{pgfonlayer}
	\begin{pgfonlayer}{edgelayer}
		\draw [style=c2] (484) to (479);
		\draw [style=c2] (479) to (485);
		\draw [style=c2] (479) to (486);
		\draw [style=c0] (479) to (487);
		\draw [style=c0] (488) to (479);
		\draw [style=c2] (481) to (480);
		\draw [style=c2] (480) to (482);
		\draw [style=c2] (480) to (483);
		\draw [style=c0] (480) to (489);
		\draw [style=c0] (480) to (490);
		\draw [style=c2] (486) to (481);
		\draw [style=c2, bend left=45] (482) to (485);
		\draw [style=c2, bend right=105, looseness=1.25] (484) to (483);
		\draw [style=c0, bend left, looseness=0.50] (487) to (490);
		\draw [style=c0, bend right=120] (489) to (488);
		\draw [style=c0, bend left=15] (490) to (487);
		\draw [style=c0, bend left=90, looseness=0.50] (488) to (489);
		\draw [style=c2, bend right] (486) to (481);
		\draw [style=c2, bend left, looseness=0.75] (482) to (485);
		\draw [style=c2, bend right=105, looseness=1.50] (484) to (483);
		\draw [style=c0] (491.center) to (492.center);
		\draw [style=c0] (493.center) to (494.center);
		\draw [style=c0] (494.center) to (495.center);
		\draw [style=c0] (494.center) to (496.center);
		\draw [style=c2] (505) to (500);
		\draw [style=c2] (500) to (506);
		\draw [style=c2] (500) to (507);
		\draw [style=c2] (502) to (501);
		\draw [style=c2] (501) to (503);
		\draw [style=c2] (501) to (504);
		\draw [style=c2] (507) to (502);
		\draw [style=c2, bend left=45] (503) to (506);
		\draw [style=c2, bend right=105, looseness=1.25] (505) to (504);
		\draw [style=c0, bend left, looseness=0.50] (508) to (511);
		\draw [style=c0, bend right=120] (510) to (509);
		\draw [style=c0, bend left=15] (511) to (508);
		\draw [style=c0, bend left=90, looseness=0.50] (509) to (510);
		\draw [style=c2, bend right] (507) to (502);
		\draw [style=c2, bend left, looseness=0.75] (503) to (506);
		\draw [style=c2, bend right=105, looseness=1.50] (505) to (504);
		\draw [style=c0] (512.center) to (513.center);
		\draw [style=c0] (514.center) to (515.center);
		\draw [style=c0] (515.center) to (516.center);
		\draw [style=c0] (515.center) to (517.center);
		\draw [style=c0] (500) to (518);
		\draw [style=c0] (518) to (509);
		\draw [style=c0] (508) to (518);
		\draw [style=c0] (501) to (519);
		\draw [style=c0] (519) to (511);
		\draw [style=c0] (510) to (519);
		\draw [style=c2] (532) to (527);
		\draw [style=c2, bend left=45] (528) to (531);
		\draw [style=c2, bend right=105, looseness=1.25] (530) to (529);
		\draw [style=c0, bend left, looseness=0.50] (533) to (536);
		\draw [style=c0, bend right=120] (535) to (534);
		\draw [style=c0, bend left=15] (536) to (533);
		\draw [style=c0, bend left=90, looseness=0.50] (534) to (535);
		\draw [style=c2, bend right] (532) to (527);
		\draw [style=c2, bend left, looseness=0.75] (528) to (531);
		\draw [style=c2, bend right=105, looseness=1.50] (530) to (529);
		\draw [style=c0] (537) to (534);
		\draw [style=c0] (533) to (537);
		\draw [style=c0] (538) to (536);
		\draw [style=c0] (535) to (538);
		\draw [style=c2] (531) to (540);
		\draw [style=c0] (540) to (542);
		\draw [style=c2] (542) to (532);
		\draw [style=c0] (542) to (539);
		\draw [style=c0] (539) to (537);
		\draw [style=c0] (539) to (541);
		\draw [style=c2] (541) to (530);
		\draw [style=c0] (541) to (540);
		\draw [style=c2] (528) to (543);
		\draw [style=c0] (543) to (544);
		\draw [style=c2] (544) to (527);
		\draw [style=c0] (544) to (545);
		\draw [style=c0] (545) to (538);
		\draw [style=c0] (545) to (546);
		\draw [style=c2] (546) to (529);
		\draw [style=c0] (546) to (543);
	\end{pgfonlayer}
\end{tikzpicture}}}
    \caption{The graph $G_3$ is a nice graph with a tuplet $\{x,y\}$. Some of its $2$-components have been coloured \rcol. The graph $H_3$ is obtained from $G_3$ by coadding an edge at $x$ with endvertices $x_A$ and $x_B$ and at $y$ with endvertices $y_A$ and $y_B$. The graph $X_3$ is obtained from $H_3$ by coadding a $4$-cycle at the vertices $x_A$ and $y_A$.}
    \label{fig:full_enrichment}
\end{figure}

Instead of proving \autoref{lem:amazing_expansion_main} directly, we will split its proof into two steps. For that, we shall make use of the following definition.
\begin{dfn}[almost $A$-mazing]
    Given a temporal triple $\tseq{G}$ and an edge set $A$, we say that $\tseq{G}$ is \emph{almost $A$-mazing} if $A$ is an upleaf rooted at $1$ and for every tuplet $s$ of $G_3$, either there is at most one $A$-ffected $2$-component at $s$ or at most one $2$-component at $s$ which is not $A$-ffected.
\end{dfn}

In the following lemma, we coadd an edge at some vertices in tuplets. For an illustration of this, see \autoref{fig:full_enrichment}.
\begin{lem}\label{lem:almost_amazing}
    Let $\tseq{G}$ be a temporal triple of planar graphs and $A$ a Tutte leaf of $G_1$. Suppose that $G_3$ contains no upleaves, that $G_3$ is nice and that $\tseq{G}$ is not an amazingly-bad obstruction. Then there exists a nice accurate planar expansion $H_3$ of\\ $\tseq{G}$ so that the expansion triple $(G_1\minorm G_2\minorp H_3)$ is almost $A$-mazing and not an amazingly-bad obstruction. Furthermore we have that\vspace{-.25cm}
    \[
    \Cbb(H_3)\leq\Cbb(G_3)\quad\text{and}\quad |E(H_3)|\leq |E(G_3)|+2\,|\Cbb(G_3)|.\vspace{-.25cm}
    \]
\end{lem}
\begin{pf}
    For each vertex $x$ in a tuplet $t$ of $G_3$ with more than one $A$-ffected $2$-component and more than one $2$-component which is not $A$-ffected, let $J_x$ be the graph with one edge, whose vertices are labeled $x_A$ and $x_B$. Let $\theta_x$ be the function from the edges incident with $x$ in $G_3$ to $V(J_x)$ so that the edges from $A$-ffected $2$-components are sent to $x_A$ and the edges from $2$-components which are not $A$-ffected are sent to $x_B$. Let $H_3$ be the graph obtained from $G_3$ coadding $J_x$ at $x$ in $G$ for each $x$. By \autoref{lem:2-con_lifts}, the graph $H_3$ is $2$-connected.
    \begin{claim}\label{claim:almost_amazing_tuplets}
        All tuplets of $H_3$ are either tuplets of $G_3$ or of the form $\{x_A,y_A\}$ or $\{x_B,y_B\}$ for some tuplet $\{x,y\}$ of $G_3$.
    \end{claim}
    \begin{cproof}
        By repeated application of \autoref{lem:tuplets_of_tree_coadd} we conclude that all tuplets of $H_3$ are either tuplets of $G_3$ or of the form $\{x_P,y_Q\}$ for some tuplet $\{x,y\}$ of $G_3$ and $P,Q\in\{A,B\}$. In the latter case, we must have $P=Q$ since otherwise the $2$-separator $\{x_P,y_Q\}$ has only two $2$-components.
    \end{cproof}

    The fact that the graph $H_3$ is nice follows from \autoref{claim:almost_amazing_tuplets} and the fact that $H_3$ is $2$-connected. Observe that by construction all but one of the $2$-components at a tuplet $\{x_A,y_A\}$ is $A$-ffected and all but one of the $2$-components at a tuplet $\{x_B,y_B\}$ are not $A$-ffected. Furthermore, if $\{x,y\}$ is a tuplet of both $G_3$ and $H_3$, then it has by assumption either at most one $A$-ffected $2$-component or at most one $2$-component which is not $A$-ffected. Hence by \autoref{claim:almost_amazing_tuplets}, we get that $H_3$ is almost $A$-mazing.

    \begin{claim}\label{claim:almost_lift}
         Let $\phi_3$ be an embedding of $G_3$ so that for every tuplet $t$ of $G_3$, the $A$-ffected $2$-components at $t$ form a segment in the cyclic structure induced on them by $\phi_3$. Then $\phi_3$ lifts to an embedding of $H_3$.
    \end{claim}
    \begin{cproof}
        We get that at each vertex $x$ of $t$, the edges from $A$-ffected $2$-components form a segment in the rotator $\phi_3(x)$. Let $\phi_3'$ be the embedding obtained from $\phi_3$ by performing the following operation at each vertex $x$ of a tuplet $t$ of $G_3$ that has at least two $A$-ffected $2$-components and at least two $2$-components which are not $A$-ffected: First subdivide each edge incident with $x$ in $G_3$ so that the subdivision vertices lie on the boundary of a small disk around $x$ and then delete $x$ and its incident edges. Now by our previous observation, the subdivision vertices from edges in $A$-ffected components form a segment in their cyclic order in the boundary of the small disk. We identify these vertices into a new vertex $x_A$ by a contraction of the plane inside the boundary of the small disk. We identify all other vertices similarly into a new vertex $x_B$. Then we add an edge between $x_A$ and $x_B$.
        Observe that the embedding $\phi_3$ is an embedding of $H_3$ which induces $\phi_3$ on $G_3$.
    \end{cproof}
    
    \begin{claim}
        The graph $H_3$ is planar.
    \end{claim}
    \begin{cproof}
        Since $G_3$ is planar, there exists an embedding $\phi_3$ of $G_3$ so that the $A$-ffected $2$-components form a segment in the cyclic structure induced by $\phi_3$ on the $2$-components at $t$. By \autoref{claim:almost_lift}, the embedding $\phi_3$ lifts to an embedding of $H_3$ and we conclude that $H_3$ is planar.
    \end{cproof}

    \begin{claim}
        The graph $H_3$ is an expansion of $\tseq{G}$.
    \end{claim}
    \begin{cproof}
        Let $(\phi_1,\phi_2,\phi_3)$ be a simultaneous embedding of $\tseq{G}$. For each tuplet $t$ of $G_3$, we have by \autoref{lem:A-ffected} that $A$-ffected $2$-components form a segment in the cyclic structure induced on them by $\phi_3$. By \autoref{claim:almost_lift}, $\phi_3$ lifts to an embedding of $H_3$. This completes the proof of our claim.
    \end{cproof}
    
    \begin{claim}
        The graph $H_3$ is accurate for $G_3$.
    \end{claim}
    \begin{cproof}
        Let $c$ be a $3$-block of $H_3$ and $b$ and $b'$ two $3$-blocks of $G_3$ that are minors of $c$. We are required to show that $b$ and $b'$ are torsos of the same bag of the Tutte decomposition of $G_3$. Let $B$ and $B'$ denote the bags of the Tutte decomposition of $G_3$ corresponding to $b$ and~$b'$, respectively. By \autoref{claim:almost_amazing_tuplets}, if $B$ and $B'$ can be distinguished by some $2$-separator of $G_3$, they can be separated by a pair of vertices in $H_3$. But this cannot happen since both $B$ and $B'$ (after possibly relabelling some vertices $x$ to $x_A$ or $x_B$ in each $B,B'$) are contained in the bag of the Tutte decomposition of $H$ corresponding to~$c$. Hence $B=B'$ and this completes our proof.
    \end{cproof}

    \begin{claim}
        The triple $(G_1\minorm G_2\minorp H_3)$ is not an amazingly-bad obstruction.
    \end{claim}
    \begin{cproof}
        Let $c$ be a $2$-component at some tuplet $s$ of $H_3$. We are required to show that the set $x(c)$ (of $2$-components in the shadow of $c$) in $G_2$ is affectively-connected. There are two cases.
        
        \textbf{Case 1.} The edges in $c$ meeting $s$ are of the form $x_Ax_B$ and $y_Ay_B$ for some tuplet $t=\{x,y\}$ of $G_3$. Then $x(c)$ either consists of all of the $2$-components at the shadow of $s$ in $G_2$ which are $A$-ffected or consists of those which are not $A$-ffected.
        In either case, we get that $x(c)$ is affectively-connected.

        \textbf{Case 2.} The edges in $c$ meeting $s$ correspond to the the edges of a $2$-component $c'$ at the shadow of $s$ in $G_3$ which meet the shadow of $s$. Now we simply apply the fact that $\tseq{G}$ is not an amazingly-bad obstruction to get that $x(c)$ is affectively-connected.
    \end{cproof}

    \begin{claim}
        We have that $\Cbb(H_3)\leq\Cbb(G_3)$.
    \end{claim}
    \begin{cproof}
        Let $t$ be a tuplet in $G_3$ and suppose that $t$ has $p$ many $A$-ffected $2$-components and $q$ many $2$-components which are not $A$-ffected. Then if both $p,q\geq 2$ we have that $t$ is replaced in $H_3$ by a pair of tuplets, with $p+1$ and $q+1$ many $2$-components, respectively. Since $(p+1)-2+(q+1)-2=(p+q)-2$, we get $\beta(H_3)=\beta(G_3)$ and this completes the proof of our claim.
    \end{cproof}
    
    Since to obtain $H_3$ we coadd at most two edges for each tuplet in $G_3$ and there are at most $\Cbb(G_3)$ tuplets in $G_3$, we obtain $|E(H_3)|\leq|E(G_3)|+2\,\Cbb(G_3)$. Combining this and the above claims complete our proof.
\end{pf}

In the following two lemmas, we coadd a cycle at a vertex in a tuplet. For an illustration of our construction, see \autoref{fig:full_enrichment}.
Instead of rewriting a long hypothesis, we introduce the following setting.
\begin{setting}\label{set:coadd_cycle}
    Let $G$ be a nice graph and $\{x,y\}$ a tuplet of $G$ with $2$-components $c_1,\ldots,c_k$. Let $C_k$ be a $k$-cycle, denote its vertices by $x_1,\ldots,x_k$ (in its cyclic ordering) and let $\theta$ be a function mapping the unique edge incident with $x$ in $c_i$ to $x_i$. Let $H$ be obtained from $G$ by coadding $C_k$ at $x$ with respect to $\theta$.
\end{setting}
\begin{lem}\label{lem:coadd_cycle}
    Assume \autoref{set:coadd_cycle}. Then every embedding of $G$ where the cyclic structure of the $2$-components at $t$ is $(c_1,c_2,\ldots,c_k)$ is induced by an embedding of~$H$.
\end{lem}
\begin{pf}
    Let $\iota$ be an embedding of $G$ where the cyclic structure of the $2$-components at $t$ is $(c_1,\ldots,c_k)$. Now let $\iota'$ be the embedding obtained from $\iota$ by subdividing the edges incident with $x$ and then adding the edges of the cycle $C_k$ between the subdivision vertices on the boundary of a small disc around $x$. We label the subdivision vertices $x_1,\ldots,x_n$ according to the $2$-components and place each edge of $C_k$ so that we end up with a copy of $C_k$ in our embedded graph. Now let $\iota''$ be the embedding obtained from $\iota'$ by deleting $x$. Observe that $\iota''$ is an embedding of $H$ and furthermore $\iota''$ induces $\iota$ on~$G$.
\end{pf}

\begin{lem}\label{lem:coadd_accurate}
    Assume \autoref{set:coadd_cycle}. Then $H$ is accurate for $G$.
\end{lem}
\begin{pf}
    Let $c$ be a $3$-block of $H$ and $b$ and $b'$ two $3$-blocks of $G$ that are minors of $c$. We are required to show that $b$ and $b'$ are torsos of the same bag of the Tutte decomposition of $G$. Let $B$ and $B'$ denote the bags of the Tutte decomposition of $G$ corresponding to $b$ and~$b'$, respectively. Let $C'$ denote the bag of the Tutte decomposition of $H$ corresponding to $c$. Suppose for a contradiction that $B$ and $B'$ can be distinguished by some $2$-separator of $G$. First suppose that $B$ and $B'$ are distinguished by $\{x,y\}$. Then there exists $2$-components $c_i$ and $c_j$ at $\{x,y\}$, containing $B$ and $B'$, respectively. We observe that $\{x_i,y\}$ is a $2$-separator of $H$ which separates the bag $C$, contradicting the fact that it corresponds to a $3$-block. So assume instead that $B$ and $B'$ are distinguished by some $2$-separator $t$ of $G$ but are contained in the same $2$-component $c_i$ at $\{x,y\}$ (and so $t\neq \{x,y\}$). Suppose first that $x\in t$. Then by replacing $x$ by $x_i$ in $t$, we obtain a $2$-separator of $H$ which separates $C$, a contradiction. If instead $t$ does not contain $x$, then $t$ is a $2$-separator of $H$, separating $C$, which is similarly a contradiction.
    Hence $B=B'$ and this completes our proof.
\end{pf}

\begin{lem}\label{lem:coadd_tuplet_correspondence}
    Let $G$ be a nice graph and $\{x,y\}$ a tuplet of $G$ with $2$-components $c_1,\ldots,c_k$. Let $C_x$ be a $k$-cycle, denote its vertices by $x_1,\ldots,x_k$ (in its cyclic ordering) and $\theta_x$ a function mapping the unique edge incident with $x$ in $c_i$ to $x_i$. Let $C_y$ and $\theta_y$ be defined similarly except with `$y$' in place of `$x$'. If $H$ is obtained from $G$ by coadding $C_x$ at $x$ with respect to $\theta_x$ and $C_y$ at $y$ with respect to $\theta_y$, then every tuplet of $H$ is a tuplet of $G$.
\end{lem}
\begin{pf}
    Observe that no pair of vertices $\{x_i,x_j\}$ separates the graph $H$. Hence, every $2$-separator of $H$ is a separator of $G$ or of the form $\{x_i,y\}$, $\{x,y_i\}$ or $\{x_i,y_i\}$ for some $2$-separator $\{x,y\}$ of $G$ and choice of $i$. Now observe that each of these type of $2$-separator have only two $2$-components in $H$. This completes the proof of our lemma.
\end{pf}

\begin{pf}[Proof of \autoref{lem:amazing_expansion_main}]
    Let $\tseq{G}$ be a temporal triple of planar graphs and $A$ a Tutte leaf of $G_1$. Suppose that $G_3$ contains no upleaves, that $G_3$ is nice and that $\tseq{G}$ is not an amazingly-bad obstruction. Then we are required to show that there exists an nice accurate planar expansion $H_3$ of $\tseq{G}$ so that the expansion triple $(G_1\minorm G_2\minorp H_3)$ is $A$-mazing and furthermore
    \[
    \Cbb(X_{i+2})\leq \Cbb(G_{i+2})\quad\text{and}\quad |E(X_{i+2})|\leq |E(G_{i+2})| + 4\;\Cbb(G_{i+2}).
    \]

    Assume first that $\tseq{G}$ is almost $A$-mazing. Below, we deduce this special case from the general case.

    So let $t$ be a tuplet of $G_3$ with more than one $A$-ffected $2$-component. Since $\tseq{G}$ is almost $A$-mazing, $t$ has at most one $2$-component that is not $A$-mazing. Hence the (cyclic or linear) structure suggested by $A$ on the $A$-ffected $2$-components at $t$ extends uniquely to a cyclic structure on the set of all $2$-components at $t$. For each such tuplet $t$ of $G_3$, we let $c_1,\ldots,c_k$ denote its $2$-components in an order which induces this unique cyclic structure.
    
    Now for each vertex $x$ in a tuplet $t$ of $G_3$ with more than one $A$-ffected $2$-component, let $C_x$ be a $k$-cycle with $k$ the number of $2$-components at $t$, denote its vertices by $x_1,\ldots,x_k$. Let $\theta_x$ be the function assigning each the unique edge in each $2$-component $c_i$ meeting $x$ to the vertex~$x_i$.

    Let $H_3$ be the graph obtained from $G_3$ by coadding each $C_x$ with respect to the functions~$\theta_x$. By \autoref{lem:coadd_accurate}, we have that $H_3$ is accurate for $G_3$. By \autoref{lem:coadd_cycle} and \autoref{lem:A-ffected}, the graph $H_3$ is an expansion for $\tseq{G}$. By \autoref{lem:2-con_lifts}, $H_3$ is $2$-connected. We verify the following claim.

    \begin{claim}
        All tuplets of $H_3$ are tuplets of $G_3$.
    \end{claim}
    \begin{cproof}
        This claim follows via repeated applications of \autoref{lem:coadd_tuplet_correspondence}.
    \end{cproof}

    By the above claim, the only tuplets in $H_3$ are those of $G_3$ which have at most one $A$-ffected $2$-component or with at most one $2$-component which is not $A$-ffected. Furthermore, a $2$-component at a tuplet of $H_3$ is $A$-ffected if and only if its corresponding minor at a tuplet in $G_3$ is $A$-ffected. We conclude that $H_3$ is $A$-mazing.
    \begin{claim}\label{claim:coadd_cycle_flex}
        We have that $\Cbb(H_3)\leq\Cbb(G_3)$.
    \end{claim}
    \begin{cproof}
        Let $t$ be a tuplet of $G_3$, suppose $t$ has $p\geq 3$ many $2$-components, all but one of which are $A$-ffected. Then $t$ is replaced in $H_3$ by $p$ many totally-nested $2$-separators having exactly two $2$-components. We observe that $\alpha(H_3)+3\beta(H_3)\leq \alpha(G_3)+3\beta(G_3)$, as required.
    \end{cproof}
    
    Since $H_3$ is obtained from $G_3$ by coadding at most as many edges at each vertex in a tupet of $G_3$ as there are $2$-components at the tuplet, we obtain $|E(H_3)|\leq|E(G_3)|\leq 2\;\Cbb(G_3)$. This completes our proof in the case where $\tseq{G}$ is almost $A$-mazing.

    Now suppose that $\tseq{G}$ is not almost $A$-mazing. Then we may apply \autoref{lem:almost_amazing} and then proceed analogously to above. In this case we end up with\vspace{-.25cm}
    \[
    |E(H_3)|\leq|E(G_3)|+2\,\Cbb(G_3)+2\,\Cbb(G_3)=|E(G_3)|+4\,\Cbb(G_3)\vspace{-.25cm}
    \]
    as required.
\end{pf}

\section{Delightful analogues}\label{sec:delightful_analogues}

This section is dedicated to proving the following, which we recall from \autoref{sec:main}.
\begingroup
\def\thelem{\ref{lem:B_2}}
\begin{lem}[Delightful Analogue Lemma]
    Let $\Gcal=(G_1\minorm\ldots\minorp G_n)$ be a nice $2$-connected temporal sequence  with a list $L$,
and $A$ a last upleaf of $\Gcal$ consisting of a single edge and rooted at $i\in [n]$. If $\Gcal$ contains no gaudy obstruction, then $(\Gcal,L)$ admits a delightful analogue for $(d,i)$.
\end{lem}
\addtocounter{lem}{-1}
\endgroup

\subsection{Towards the proof of  The \nameref{lem:B_2}}

We will first derive \myref{lem:B_2} from the following weaker version, which the rest of \autoref{sec:delightful_analogues} is devoted to proving.
\begin{lem}[Almost Delightful Lemma]\label{lem:B_1}
    Let $\Gcal=(G_1\minorm\ldots\minorp G_n)$ be a nice temporal sequence of $2$-connected graphs along with a list $L$ and suppose that $A$ is a last upleaf rooted at some $i\in[n]$. Suppose that $A$ consists of a single edge $d$, spanning a $2$-separator $\{x,y\}$ in $G_i$. Suppose that at least three $2$-components at $\{x,y\}$ in $G_i$ have nontrivial shadow in $G_{i+1}$. If $\Gcal$ contains no gaudy obstruction, then there exists a nice graph $H_{i+2}\minorm G_{i+2}$ and a pair of lists $\{M,N\}$ for the temporal sequence\vspace{-.25cm}
    \[
    \Gcal':=(G_1\minorm\ldots\minorm G_{i-1}\minorp G_i-d\minorm G_{i+1}-d\minorp H_{i+2}\minorm G_{i+1}\minorp\ldots\minorp G_n)\vspace{-.25cm}
    \]
    so that $(\Gcal',\{M,N\})$ is an analogue of $(\Gcal,L)$ and furthermore, we have that\vspace{-.25cm}
    \[
    \Cbb(H_{i+2})\leq \Cbb(G_{i+2})\quad\text{and}\quad |E(H_{i+2})|\leq |E(G_{i+2})| + 12\,(\Cbb(G_{i+2}))^2.\vspace{-.25cm}
    \]
\end{lem}

\begin{pf}[Proof of the \nameref{lem:B_2} assuming the \nameref{lem:B_1}]
    Let $\Gcal=$\\ $(G_1\minorm\ldots\minorp G_n)$ be a nice temporal sequence of $2$-connected graphs along with a list $L$. Suppose that a last upleaf $A$ consists of a single edge $d$ and is rooted in $i\in [n]$. Assume that $\Gcal$ contains no gaudy obstruction. We are required to show that there exists a temporal sequence\vspace{-.25cm}
    \[
    \Gcal':=(G_1\minorm\ldots\minorm G_{i-1}\minorp G_i-d\minorm H_{i+1}\minorp H_{i+2}\minorm G_{i+1}\minorp\ldots\minorp G_n)\vspace{-.25cm}
    \]
    along with a pair of lists $M$ and $N$ so that $H_{i+1}\minorm G_{i+1}-d$ and $H_{i+2}\minorm G_{i+2}$ are nice, we have $(\Gcal',\{M,N\})$ is an analogue of $(\Gcal,L)$ and furthermore\vspace{-.25cm}
    \[
        |E(H_{i+1})|\leq |E(G_{i+1})|+1,\vspace{-.25cm}
    \]
    \[
    \Cbb(H_{i+2})\leq \Cbb(G_{i+2})+5\quad\text{and}\quad |E(H_{i+2})|\leq |E(G_{i+2})| + 12\,(\Cbb(G_{i+2}))^2+1.\vspace{-.25cm}
    \]
    
    Denote by $t$ the $2$-separator in $G_i$ spanned by $d$. Since one of the $2$-components at $t$ is the single edge $d$, the $2$-separator $t$ has at least three $2$-components (is a tuplet). If at least three $2$-components at $t$ have nonempty shadow in $G_{i+1}$, then we are done by \myref{lem:B_1}, where we set $H_{i+1}=G_{i+1}$. Furthermore, if $d$ does not appear in $G_{i+1}$ then we are done immediately by setting $H_{i+1}=G_{i+1}$ and $H_{i+2}=G_{i+2}$. So assume that $d$ appears in $G_{i+1}$ and hence also in $G_{i+2}$ and that there is a $2$-component at $t$ apart from $d$ that has empty shadow in $G_{i+1}$; pick one such $2$-component and denote it by $c$. Let $e$ be an edge from $c$. Now let $H_{i+1}$ be the graph obtained from $G_{i+1}$ by adding the edge $e$ in parallel to $d$. Similarly, let $H_{i+2}$ be the graph obtained from $G_{i+2}$ by adding the edge $e$ in parallel to $d$. Observe that\vspace{-.25cm}
    \[
        \Hcal:=(G_1\minorm\ldots\minorm G_{i-1}\minorp G_i\minorm H_{i+1}\minorp H_{i+2}\minorm G_{i+3}\minorp\ldots\minorp G_n)\vspace{-.25cm}
    \]
    is a nice $2$-connected temporal sequence. Furthermore, since $\Hcal$ has the same $3$-blocks as $\Gcal$, the list $L$ is also a list for $\Hcal$.

    \begin{claim}
        We have that $\Cbb(H_{i+2})\leq \Cbb(G_{i+2})+5$.
    \end{claim}
    \begin{cproof}
        Since we may have created a new totally-nested $2$-separator or added a $2$-component to a tuplet, we get $\alpha(H_{i+2})+3\beta(H_{i+2})\leq \alpha(G_{i+2})+3\beta(G_{i+2})+3$. However we also get in this case that $\gamma(H_{i+2})=\gamma(G_{i+2})+1$.
    \end{cproof}
    
    Furthermore for both vertices $x\in t$, since the edge $e$ receives its own unique colour in the $(t,x)$-centred colouring of $G_{i+2}$ and $\Gcal$ contains no gaudy obstruction, we also get that $\Hcal$ contains no gaudy obstruction. Now we have verified that $(\Hcal,L)$ satisfies the hypotheses of \myref{lem:B_1}, it remains to show that $(\Hcal,L)$ is an analogue of $(\Gcal,L)$.

    Clearly a simultaneous list-embedding of $(\Hcal,L)$ induces a simultaneous embedding of $(\Gcal,L)$. So let a simultaneous list-embedding $(\phi_1,\ldots,\phi_n)$ of $(\Gcal,L)$ be given. Let $\phi_{i+1}'$ be the embedding induced on $H_{i+1}$ by $\phi_i$.
    Since at most one other $2$-component other than $d$ has nonempty shadow in $G_{i+1}$, we get by \autoref{lem:comp_segs} that $e$ appears on one of the two sides of $d$ in the rotator at each endvertex of $e$. Let $x_{i+1}$ be the endvertex of $d$ in $H_{i+1}$ in the rotator of which $e$ appears after $d$. Since $(\phi_1,\ldots,\phi_n)$ is a simultaneous embedding, $\phi_i$ and $\phi_{i+2}$ induce the same embedding $\phi_{i+1}$ on $G_{i+1}$ and since $\phi_{i+1}'$ is induced by $\phi_{i+1}$, it also induces $\phi_{i+1}$ on $G_{i+1}$. Let $x_{i+2}$ be the endvertex of $d$ in $H_{i+2}$ which lies in the branch set of $x_{i+1}$. By adding the edge $e$ to the embedding $\phi_{i+2}$ so that it is in parallel to the edge $d$, and after $d$ in the rotator at $x_{i+2}$, we obtain an embedding $\phi_{i+2}'$ of $H_{i+2}$ that induces $\phi_{i+1}'$ on $H_{i+1}$ and induces $\phi_{i+3}$ on $G_{i+3}$. Hence $\Phi:=(\phi_1,\ldots,\phi_i,\phi_{i+1}',\phi_{i+2}',\phi_{i+3},\ldots,\phi_n)$ is a simultaneous embedding of $\Hcal$. Furthermore, the embeddings $\phi_{i+1}'$ and $\phi_{i+1}$ induce the same embeddings on the $3$-blocks of $H_{i+1}$, which are the same as the $3$-blocks as $G_{i+1}$. Similarly, the embeddings $\phi_{i+2}'$ and $\phi_{i+2}$ induce the same embeddings on the $3$-blocks of $H_{i+2}$, which are the same as the $3$-blocks as $G_{i+2}$. Hence $\Phi$ satisfies the list $L$ and this completes the proof.
\end{pf}

To prove \myref{lem:B_1} we first make the following preparations.

We remind the reader that in the circumstances when we are given a graph $G_2$ and an edge $d$, which is not in $G_2$, we follow the convention that $G_2-d$ denotes a graph that is isomorphic to $G_2$; this way $G_2-d$ is always defined, independent of whether $d\in E(G_2)$.

\begin{lem}\label{lem:put_in_edge}
    Let $\tseq{G}$ be a temporal triple, and let $t=\{x,y\}$ be a tuplet of $G_1$.
    Let $b$ be the $(t,x)$-centred bond of $G_3$ along with the $(t,x)$-centred partial colouring, and assume it has at least three colours.
For every simultaneous embedding $(\psi_i,\psi_{i+1},\psi_{i+2})$ of $(G_1-d\minorp G_{2}-d \minorm G_{3})$ so that $\psi_{i+2}$ is $b$-respectful,
there is a simultaneous embedding $(\phi_i,\phi_{i+1},\phi_{i+2})$ of $\tseq{G}$ that induces $(\psi_i,\psi_{i+1},\psi_{i+2})$ on $(G_1-d\minorp G_{2}-d \minorm G_{3})$. 
\end{lem}

\begin{pf}
    Let a simultaneous embedding $(\phi_1,\phi_2,\phi_3)$ for $(G_1-d\minorm G_2-d\minorp G_3)$ be given. If the edge $d$ is not in $G_2$, then we can easily extend $(\phi_1,\phi_2,\phi_3)$ to a simultaneous embedding of $\tseq{G}$: indeed, just embed $d$ into an arbitrary face of $\phi_1$ whose boundary contains both $x$ and $y$. So from now on, assume that $d$ is an edge of $G_2$. 

    Let $\psi_3=\phi_3$, and let $\psi_2$ be the embedding induced on $G_2$ by $\psi_3$. Note that the embedding $\psi_2$ of $G_2$ induces $\phi_2$ on $G_2-d$ by definition. That is, the embedding $\psi_2$ is obtained from $\phi_2$ by embedding the edge $d$ into one of its faces; denote that face by $f$. 
        Thus it remains to \lq lift\rq\ the embedding $\phi_1$ from $G_1-d$ to an embedding $\psi_1$ of $G_1$ that induces $\psi_2$ on $G_2$. 

Since the embedding $\phi_1$ induces the embedding $\phi_2$, the face $f$ of $\phi_2$ is a union of faces of $\phi_1$. Denote the set of these faces of $\phi_1$ by $S(f)$.
\begin{claim}\label{face-exists}
There is a face $g\in S(f)$ that has both vertices $x$ and $y$ in its boundary.
\end{claim}
\begin{cproof}
    Let $\pi t=\{\pi x,\pi y\}$ denote the shadow of $t$ in $G_2$. Since $G_3$ is coloured using at least three colours from the $(t,x)$-centred colouring, we know that $t$ has at least three $2$-components with nontrivial shadow in $G_2$. Hence $\pi t$ is a tuplet.
    Let $p$ and $q$ denote the $2$-components at the shadow $\pi t=\{\pi x,\pi y\}$ of $t$ in $G_2$ that bound the face $f$ in $\phi_2$. Observe that since $\pi t$ is a tuplet, the $2$-components $p$ and $q$ are distinct. Suppose for a contradiction that $p$ and $q$ are contained as minors in the same $2$-component at $t$ in $G_1$. Hence they inherit the same colour, say \rcol, in the $(t,x)$-centred colouring. Now consider the bond $b$ of $G_3$ along with the cyclic order induced on it by $\phi_3$. We get that the nearest coloured edges on either side of $d$ in this cyclic order have colour \rcol. Furthermore, by assumption there is some other edge with colour distinct from \rcol\ that isn't $d$. From this we get that $\phi_3$ is not a bond-respectful embedding for the bond graph $(G_3,b)$. This contradicts the fact that the $(t,x)$-centred colouring is well-respected.

    So $p$ and $q$ are contained as minors in distinct $2$-components at $t$, say $p'$ and $q'$, respectively. Since $p'$ and $q'$ are distinct $2$-components at $t$, there are no paths in $G_1$ that join $p'$ and $q'$ without using $x$ or $y$. Consider the union of the boundaries of the faces in $S(f)$. This subgraph of $G_1-d$ contains no path from $p'$ to $q'$ that does not use $x$ or $y$. There are hence two faces of $\phi_1$ whose boundaries  both contain $x$ and $y$ and whose interiors along with $x$ and $y$ topologically separate $p'$ from $q'$. One of these such faces is contained in $S(f)$, as desired.
\end{cproof}

Let $g$ be a face as in \autoref{face-exists}. We obtain the embedding $\psi_1$ from the embedding $\phi_1$ by embedding the edge $d$ into the face $g$. By construction, $\psi_1$ induces $\phi_1$, and thus $\phi_2$. So by the choice of $g$, the embedding $\psi_1$ induces $\psi_2$. Hence  $(\psi_1,\psi_2,\psi_3)$ is a simultaneous embedding  for $(G_1\minorm G_2\minorp G_3)$.
\end{pf}

Let $(G,b)$ and $(H,b)$ be a partially-coloured bond graph with the same bond $b$. We say that $(H,b)$ is a \emph{model} of the pair $(G,b)$ if $H\minorm G$ and every $b$-respectful embedding of $G$ is induced by an embedding of $H$. We say that a model $(H,b)$ of a pair $(G,b)$ is an \emph{accurate model} if $H$ is accurate for $G$.

Given a temporal sequence $\tseq{G}$, and an edge $d$ so that $\{d\}$ is a last upleaf of $G_1$, we would ideally like to build a $b$-respectful model of $(G_3,b)$, where $b$ is the $(t,x)$-centred bond of $G_3$ along with the $(t,x)$-centred partial colouring. However, this is not always possible, as shown in \autoref{eg:model_impossible} below.
Instead in \autoref{sec:bond_fitting} we will build a model which is almost bond-respectful. Then in \autoref{sec:del_red} we will bridge the gap by using lists.

Let $(G,b)$ be a bond graph along with a partial \rcol/\bcol-colouring. For each $2$-component $c$ at some $2$-separator in $G$, we say that $c$ \emph{contains colour f} if one of the edges of $c$ is coloured by $f$ in $b$. If $c$ contains exactly one colour, then we say that it is \emph{monochromatic} in this colour. If $c$ contains both \rcol\ and \bcol\ edges, then we say that it is \emph{polychromatic}.

\begin{dfn}[bond-fitting]
    Let $(G,b)$ be a bond graph along with a partial \rcol/\bcol-colouring.
    We say that a tuplet $t$ is \emph{fit} if there exists a choice of colour $f$ so that all but at most one $2$-component at $t$ is monochromatic in~$f$.
    We say that $(G,b)$ is \emph{bond-fitting} if every tuplet of $G$ that crosses~$b$ is fit. 
\end{dfn}

The aim of this \autoref{sec:bond_fitting} will be to prove the following lemma. We say that a graph $G$ is \emph{badly-coloured}, if it has a $2$-separator crossing the bond at which either more than two of its $2$-components are polychromatic or one of its $2$-components has no colour.

\begin{lem}\label{lem:b-model}
    Let $(G,b)$ be a nice bond graph along with a partial \rcol/\bcol-colouring. Suppose that $G$ is not badly-coloured. Then there exists a nice bond-fitting accurate model $(H,b)$ of $(G,b)$ which is not badly-coloured. Furthermore, $\Cbb(H)\leq\Cbb(G)$ and $|E(H)|\leq|E(G)|+12\,\Cbb(G)$. 
\end{lem}

\begin{figure}
    \centering
    {\Large\scalebox{0.7}{\begin{tikzpicture}
	\begin{pgfonlayer}{nodelayer}
		\node [style=none] (341) at (-6, 4.25) {$x$};
		\node [style=none] (342) at (-6, -4.25) {$y$};
		\node [style=black dot] (423) at (-6, 3.5) {};
		\node [style=black dot] (424) at (-6, -3.5) {};
		\node [style=black dot] (425) at (-6.75, 0) {};
		\node [style=black dot] (428) at (-8.75, 0) {};
		\node [style=none] (430) at (-2.25, 0) {$\succeq$};
		\node [style=none] (447) at (5.25, 0) {$\preceq$};
		\node [style=none] (476) at (-6, -6) {$G_1$};
		\node [style=none] (477) at (1.5, -6) {$G_2$};
		\node [style=none] (478) at (10.5, -6) {$G_3$};
		\node [style=none] (505) at (1.5, 4.25) {$x$};
		\node [style=none] (506) at (1.5, -4.25) {$y$};
		\node [style=black dot] (513) at (1.5, 3.5) {};
		\node [style=black dot] (514) at (1.5, -3.5) {};
		\node [style=black dot] (518) at (10.5, -2) {};
		\node [style=black dot] (519) at (10.5, -5) {};
		\node [style=black dot] (520) at (10.5, -3) {};
		\node [style=black dot] (521) at (10.5, -4) {};
		\node [style=black dot] (522) at (9.5, -3.5) {};
		\node [style=black dot] (523) at (11.5, -3.5) {};
		\node [style=black dot] (524) at (10.5, 5) {};
		\node [style=black dot] (525) at (10.5, 2) {};
		\node [style=black dot] (526) at (10.5, 4) {};
		\node [style=black dot] (527) at (10.5, 3) {};
		\node [style=black dot] (528) at (9.5, 3.5) {};
		\node [style=black dot] (529) at (11.5, 3.5) {};
		\node [style=none] (536) at (21.25, -6) {$J_3$};
		\node [style=black dot] (537) at (21.25, -2) {};
		\node [style=black dot] (538) at (21.25, -5) {};
		\node [style=black dot] (539) at (21.25, -3) {};
		\node [style=black dot] (540) at (21.25, -4) {};
		\node [style=black dot] (541) at (20.25, -3.5) {};
		\node [style=black dot] (542) at (22.25, -3.5) {};
		\node [style=black dot] (543) at (21.25, 5) {};
		\node [style=black dot] (544) at (21.25, 2) {};
		\node [style=black dot] (545) at (21.25, 4) {};
		\node [style=black dot] (546) at (21.25, 3) {};
		\node [style=black dot] (547) at (20.25, 3.5) {};
		\node [style=black dot] (548) at (22.25, 3.5) {};
		\node [style=small back] (555) at (-7.75, 0) {};
		\node [style=none] (556) at (-7.75, 0) {$p$};
		\node [style=medium back] (557) at (-7.75, 1.25) {};
		\node [style=none] (558) at (-7.75, 1.25) {$b_1$};
		\node [style=medium back] (559) at (-6.5, 1.25) {};
		\node [style=none] (560) at (-6.5, 1.25) {$b_2$};
		\node [style=medium back] (561) at (-0.75, 0) {};
		\node [style=none] (562) at (-0.75, 0) {$b_1$};
		\node [style=medium back] (563) at (0.75, 0) {};
		\node [style=none] (564) at (0.75, 0) {$b_2$};
		\node [style=small back] (565) at (-5.25, 0) {};
		\node [style=none] (566) at (-5.25, 0) {$d$};
		\node [style=small back] (567) at (-6.5, -1.25) {};
		\node [style=none] (568) at (-6.5, -1.25) {$f$};
		\node [style=small back] (569) at (-7.75, -1.25) {};
		\node [style=none] (570) at (-7.75, -1.25) {$e$};
		\node [style=small back] (571) at (-3.75, 0) {};
		\node [style=none] (572) at (-3.75, 0) {$r$};
		\node [style=small back] (573) at (2.25, 0) {};
		\node [style=none] (574) at (2.25, 0) {$d$};
		\node [style=small back] (575) at (3.75, 0) {};
		\node [style=none] (576) at (3.75, 0) {$r$};
		\node [style=medium back] (577) at (6.75, 0) {};
		\node [style=none] (578) at (6.75, 0) {$b_1$};
		\node [style=medium back] (579) at (9.75, 0) {};
		\node [style=none] (580) at (9.75, 0) {$b_2$};
		\node [style=small back] (581) at (11.25, 0) {};
		\node [style=none] (582) at (11.25, 0) {$d$};
		\node [style=small back] (583) at (14.25, 0) {};
		\node [style=none] (584) at (14.25, 0) {$r$};
		\node [style=medium back] (585) at (8.5, 0) {};
		\node [style=none] (586) at (8.5, 0) {$a_1$};
		\node [style=medium back] (587) at (12.5, 0) {};
		\node [style=none] (588) at (12.5, 0) {$a_2$};
		\node [style=medium back] (589) at (17.5, 0) {};
		\node [style=none] (590) at (17.5, 0) {$b_1$};
		\node [style=medium back] (591) at (20.5, 0) {};
		\node [style=none] (592) at (20.5, 0) {$b_2$};
		\node [style=small back] (593) at (22, 0) {};
		\node [style=none] (594) at (22, 0) {$d$};
		\node [style=small back] (595) at (25, 0) {};
		\node [style=none] (596) at (25, 0) {$r$};
		\node [style=medium back] (597) at (19.25, 0) {};
		\node [style=none] (598) at (19.25, 0) {$a_1$};
		\node [style=medium back] (599) at (23.25, 0) {};
		\node [style=none] (600) at (23.25, 0) {$a_2$};
	\end{pgfonlayer}
	\begin{pgfonlayer}{edgelayer}
		\draw [style=c0] (428) to (423);
		\draw [style=c0] (423) to (425);
		\draw [style=c0] (425) to (424);
		\draw [style=c0] (424) to (428);
		\draw [style=c0] (428) to (425);
		\draw [style=c0, bend left, looseness=0.75] (423) to (424);
		\draw [style=c2, bend left=60, looseness=1.25] (423) to (424);
		\draw [style=c0, bend left, looseness=0.75] (513) to (514);
		\draw [style=c2, bend left=60, looseness=1.25] (513) to (514);
		\draw [style=c0, bend right, looseness=0.75] (513) to (514);
		\draw [style=c0, bend left=60, looseness=1.25] (514) to (513);
		\draw [style=c1] (518) to (522);
		\draw [style=c1] (522) to (519);
		\draw [style=c1] (519) to (523);
		\draw [style=c1] (523) to (518);
		\draw [style=c1] (518) to (520);
		\draw [style=c1] (520) to (521);
		\draw [style=c1] (521) to (519);
		\draw [style=c1] (523) to (520);
		\draw [style=c1] (520) to (522);
		\draw [style=c1] (522) to (521);
		\draw [style=c1] (521) to (523);
		\draw [style=c1] (524) to (528);
		\draw [style=c1] (528) to (525);
		\draw [style=c1] (525) to (529);
		\draw [style=c1] (529) to (524);
		\draw [style=c1] (524) to (526);
		\draw [style=c1] (526) to (527);
		\draw [style=c1] (527) to (525);
		\draw [style=c1] (529) to (526);
		\draw [style=c1] (526) to (528);
		\draw [style=c1] (528) to (527);
		\draw [style=c1] (527) to (529);
		\draw [style=c0, bend left=60, looseness=0.75] (525) to (518);
		\draw [style=c0, bend left=60, looseness=0.75] (518) to (525);
		\draw [style=c0, bend right=60, looseness=0.50] (528) to (522);
		\draw [style=c0, bend left=60, looseness=0.50] (529) to (523);
		\draw [style=c2, bend left=105, looseness=1.25] (524) to (519);
		\draw [style=c0, bend left=255, looseness=1.25] (524) to (519);
		\draw [style=c1] (537) to (539);
		\draw [style=c3] (539) to (540);
		\draw [style=c1] (540) to (538);
		\draw [style=c3] (539) to (541);
		\draw [style=c3] (540) to (542);
		\draw [style=c1] (543) to (545);
		\draw [style=c3] (545) to (546);
		\draw [style=c1] (546) to (544);
		\draw [style=c3] (545) to (547);
		\draw [style=c3] (546) to (548);
		\draw [style=c0, bend left=60, looseness=0.75] (544) to (537);
		\draw [style=c0, bend left=60, looseness=0.75] (537) to (544);
		\draw [style=c3, bend right=60, looseness=0.50] (547) to (541);
		\draw [style=c3, bend left=60, looseness=0.50] (548) to (542);
		\draw [style=c2, bend left=105, looseness=1.25] (543) to (538);
		\draw [style=c0, bend left=255, looseness=1.25] (543) to (538);
	\end{pgfonlayer}
\end{tikzpicture}}}
    \caption{A temporal sequence $\tseq{G}$ and a subgraph $J_3$ of $G_3$. The graph $G_2$ is obtained from $G_1$ by deleting the edge $p$ and contracting the edges $e$ and $f$ and is obtained from $G_3$ by contracting the \bcol\ edges and deleting the edges $a_1$ and $a_2$. The graph $J$ is obtained from $G_3$ by deleting some of the \bcol\ edges. The edges in a cycle $o$ in $J_3$ have been coloured \gcol.}
    \label{fig:model_impossible}
\end{figure}
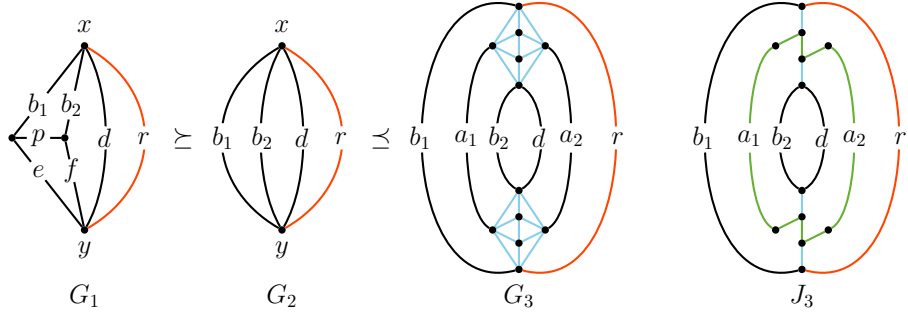

\begin{eg}\label{eg:model_impossible}
This example is similar to \autoref{eg:expansion_impossible} above, and demonstrates the necessity of lists of our approach to work in the context of delightful analogues. 
Imagine we are given a temporal sequence $(G_1,G_2,G_3,..., G_n)$ starting with the temporal triple $(G_1,G_2,G_3)$ of \autoref{fig:model_impossible}. In our approach via \autoref{lem:B_1}, we would replace the graphs $G_1$ and $G_2$ by smaller graphs, which in this case would be $K_4$ and a graph consisting of three parallel edges, respectively (if we chose $d$ as specified in the figure), and replace $G_3$ by a larger graph $H_3$ that together with the pair  of lists $\{M,N\}$ captures all the embedding information of the temporal triple $(G_1,G_2,G_3)$. 
In fact in this example having nontrivial lists is necessary as the following fact shows.
\begin{prop}
There is no planar graph $H_3\minorm G_3$ so that the 
sequences $(G_1,G_2,G_3)$ and $(G_1-d,G_2-d, H_3)$ induce the same simultaneous embeddings on their common minor
$(G_1-d,G_2-d,G_3)$.
\end{prop}

\begin{pf}
    Suppose for a contradiction that there is such a graph $H_3$. Let $J_3$ be the subgraph of $G_3$ also depicted in \autoref{fig:model_impossible} and let $o$ be the cycle in $J_3$ containing the edges $a_1$ and $a_2$. Since $J_3$ has maximum degree three, there exists a subdivision $J_3'$ of $J_3$ as a subgraph of $H_3$. Let $o'$ be the subdivision of $o$ in $J_3'$.

    Since $o$ is a separating cycle of $G_3$, $o'$ must be a separating cyclic of the planar graph $H_3$; and in any embedding of $H_3$ it must separate the edges labelled $b_1$ and $b_2$. We refer to the component of $H_3\sm o'$ containing $b_1$ as the \emph{positive side} and to the other component as the \emph{negative side}.

    There are exactly two simultaneous embeddings of $(G_1,G_2,G_3)$ up to reorientation; one is depicted in \autoref{fig:model_impossible} and the other is obtained from the depicted one by reorienting the embeddings of $G_1$ and $G_2$ and in $G_3$ reversing rotators at the endvertices of $b_1$ and $b_2$. In particular, the simultaneous embeddings induce either the cyclic structure $b_1a_1b_2da_2r$ or $ra_1db_2a_2b_1$ on the bond $\{b_1,a_1,b_2,d,a_2,r\}$ of $G_3$.

    Hence, there must exist embeddings $\iota_1$ and $\iota_2$ of $H_3$ so that the following holds: the embedding $\iota_1$ induces on the bond $\{b_1,a_1,b_2,d,a_2,r\}$ of $J_3'$ the cyclic structure $b_1a_1b_2da_2r$ and the embedding $\iota_2$ induces on the bond $\{b_1,a_1,b_2,d,a_2,r\}$ of $J_3'$ the cyclic structure $ra_1db_2a_2b_1$. Thus the positive side and the negative side of $H_3$ each do not have unique embeddings inside, and thus each must contain a $2$-separator. However, then we can flip only the one on the positive side. This induces an embedding $\kappa$ on $G_3$ in which the bond $\{b_1,a_1,b_2,d,a_2,r\}$ attains the cyclic structure $b_1db_2a_1a_2r$ (that is, obtained from $b_1a_1b_2da_2r$ with the edges $b_2$ and $d$ swapped). By our previous reasoning, $\kappa$ is part of no simultaneous embedding of $(G_1,G_2,G_3)$. However, there exists an embedding of $G_1-d$ so that it induces the same embedding as $\kappa$ on $G_2-d$. Hence we conclude that there is not such graph $H_3$. 
\end{pf}
\end{eg}

\begin{lem}\label{list_for_bonds}
    Let $(G,b)$ be a bond graph along with a partial \rcol/\bcol-colouring. Suppose that this partial colouring is not gaudy and that $G$ is bond-fitting. Then there exists a list $M$ on the $3$-blocks of $G$ so that an embedding of $G$ is bond-respectful if and only if it satisfies $M$.
\end{lem}
The proof of \autoref{list_for_bonds} will follow in \autoref{sec:del_red}. Up to the technical fact that it has a few more lists, and worse constants, the following lemma can be seen as a common generalisation of \autoref{lem:b-model} and \autoref{list_for_bonds}, which at the same time allows for a number of colours different to two. 

\begin{lem}\label{lem:model_list_workhorse}
    Let $(G,b)$ be a bond graph along with a partial colouring $\chi$ using $k\geq 1$ colours. Suppose that $\chi$ is not gaudy, and for every $2$-separator $t$ of $G$ that crosses $b$, there is a coloured edge in every $2$-components at $t$. Then there exists a nice accurate model $(H,b)$ of $(G,b)$ and a collection of lists $\{M_1,\ldots,M_l\}$ for $H$ with $l\leq k$ so that every embedding of $(H,b)$ is bond-respectful if and only if it satisfies all the lists $\{M_1,\ldots,M_l\}$. 
    
    Furthermore, both $\Cbb(H)\leq\Cbb(G)$ and $|E(H)|\leq |E(G)|+12\cdot k\cdot\Cbb(G)$.
\end{lem}
\begin{pf}
    We say that $j\in [k]$ is \emph{good} (with respect to $(G,b)$) if there is a list $M_j$ so that the following holds: consider the partial \rcol/\bcol-colouring $\chi_j$ for $(G,b)$ obtained by recolouring each $j$-coloured edge to \rcol\ and every other coloured edge to \bcol. Then every embedding of $(G,b)$ is bond-respectful with respect to $\chi_j$ if and only if it satisfies the list $M_j$.

    First suppose that all $j\in[k]$ are good. Then every embedding of $(G,b)$ is bond-respectful with respect to its original partial $k$-colouring if and only if it satisfies all the lists $M_1,\ldots,M_k$. In this case, we are done by taking $H=G$.
    Now suppose that not all elements of $[k]$ are good and let $j\in[k]$ be some index which is not good.
    \begin{claim}\label{not-bad3}
        The bond graph $(G,b)$ is not badly-coloured with respect to $\chi_j$.
    \end{claim}
    \begin{cproof}
        Since $\chi$ is not gaudy, the \rcol/\bcol-colouring $\chi_j$ is not gaudy. Then the proof of our claim follows from the assumption that every $2$-component at a $2$-separator of $G$ that crosses $b$ contains a coloured edge.
    \end{cproof}
    
    By \autoref{lem:b-model} and \autoref{not-bad3}, we obtain a nice accurate model $(H,b)$ of $(G,b)$ that is bond-fitting for $\chi_j$ and also not badly-coloured for $\chi_j$.
    
    \begin{claim}
        The bond graph $(H,b)$ is not gaudy for $\chi$.
    \end{claim}
    \begin{cproof}
Since $(G,b)$ is not gaudy for $\chi$, by \autoref{lem:respect_characterisation} there is a bond-respectful embedding $\iota$ of $(G,b)$ for $\chi$. In particular, $\iota$ is a bond-respectful embedding of $(G,b)$ for $\chi_j$.   
        Since $(H,b)$ is a model of $(G,b)$ for $\chi_j$, there is a bond-respectful embedding $\kappa$ of $(H,b)$ for $\chi_j$ so that $\kappa$ induces the embedding $\iota$ on the minor $G$ of $H$. 
        Since the embedding $\kappa$ induces $\iota$ on $G$, and 
    since $(G,b)$ and $(H,b)$ have the same bond $b$, the embedding $\kappa$ is a bond-respectful embedding of $(H,b)$ for not only $\chi_j$ but also for $\chi$. 
        By applying the converse implication from \autoref{lem:respect_characterisation}, we get that $(H,b)$ is not gaudy for $\chi$.
    \end{cproof}

    \begin{claim}
        Every $2$-component at a $2$-separator in $H$ which crosses $b$ has some coloured edge in the colouring $\chi$.
    \end{claim}
    \begin{cproof}
        This follows from the fact that $(H,b)$ is not badly-coloured for $\chi_j$.
    \end{cproof}

    By the above two claims, $(H,b)$ is not badly-coloured for $\chi$.

    \begin{claim}
        If $r\in [k]$ is good with respect to $(G,b)$, then $r$ is good with respect to $(H,b)$. 
    \end{claim}
    \begin{cproof}
        Let $M_r$ be the list witnessing that $r$ is good with respect to $(G,b)$. Since $H$ is accurate for $G$, the list $M_r$ for $G$ lifts via \autoref{lem:list_lifter} to a list $M_r'$ of $H$ so that for every embedding $\iota$ of $G$ that lifts to an embedding $\kappa$ of $H$, we have that $\iota$ satisfies $M_r$ if and only if $\kappa$ satisfies~$M_r'$.
        
        Let $\kappa$ be an embedding of $(H,b)$ that satisfies $M_r'$. Then the embedding $\iota$ that $\kappa$ induces on $(G,b)$ satisfies $M_r$. Since $r$ is good, $\kappa$ is bond-respectful for $\chi_r$. Since the bonds in both bond graphs are identical, $\iota$ is also bond-respectful for $\chi_r$. Hence $M_r'$ witnesses that the index $r$ is good with respect to the new bond graph $(H,b)$.
    \end{cproof}

    \begin{claim}
        The index $j$ is good with respect to $(H,b)$.
    \end{claim}
    \begin{cproof}
        Applying \autoref{list_for_bonds} yields a list $M_j$ for $H$ which witnesses that $j$ is good.
    \end{cproof}

    From the above two claims, we get that more indices are good for the new bond graph $(H,b)$ than were good for $(G,b)$. By induction, we are done. Indeed, at each step we get by \autoref{lem:b-model} that $\Cbb(H)\leq\Cbb(G)$ and $|E(H)|\leq|E(G)|+12\,\Cbb(G)$ and furthermore there are at most $k$ steps.
\end{pf}

\begin{pf}[Proof of \autoref{lem:B_1}]
    Let $\Gcal=(G_1\minorm\ldots\minorp G_n)$ be a nice temporal sequence of $2$-connected graphs along with a list $L$ and suppose that $A$ is a last upleaf rooted at some $i\in[n]$. Suppose that $A$ consists of a single edge $d$, spanning a $2$-separator $\{x,y\}$ in $G_i$. Suppose that at least three $2$-components at $\{x,y\}$ in $G_i$ have nontrivial shadow in $G_{i+1}$. Suppose that $\Gcal$ contains no gaudy obstruction. We are required to show that there exists a nice graph $H_{i+2}\minorm G_{i+2}$ and a pair of lists $\{M,N\}$ for the temporal sequence\vspace{-.25cm}
    \[
    \Gcal':=(G_1\minorm\ldots\minorm G_{i-1}\minorp G_i-d\minorm G_{i+1}-d\minorp H_{i+2}\minorm G_{i+3}\minorp\ldots\minorp G_n)\vspace{-.25cm}
    \]
    so that $(\Gcal',\{M,N\})$ is an analogue of $(\Gcal,L)$ and furthermore, we have that\vspace{-.25cm}
    \[
    \Cbb(H_{i+2})\leq \Cbb(G_{i+2})\quad\text{and}\quad |E(H_{i+2})|\leq |E(G_{i+2})| + 12\,(\Cbb(G_{i+2}))^2.\vspace{-.25cm}
    \]
\begin{claim}\label{boringstuff1}
There exists a nice accurate model $(H_{i+2},b)$ of $(G_{i+2},b)$ and collection of lists $\{M_1,\ldots,M_l\}$ for $H_{i+2}$ with $l\leq k$ so that every embedding of $(H_{i+2},b)$ is bond-respectful if and only if it satisfies all the lists $\{M_1,\ldots,M_l\}$. 

Furthermore, both $\Cbb(H_{i+2})\leq\Cbb(G_{i+2})$ and $|E(H_{i+2})|\leq |E(G_{i+2})|+12\cdot k\cdot\Cbb(G_{i+2})$.
\end{claim}

\begin{cproof}
It suffices to show that the hypotheses of \autoref{lem:model_list_workhorse} are met. For this, consider the temporal triple $(G_i\minorm G_{i+1}\minorp G_{i+2})$. Let $b$ denote the bond of $G_{i+2}$ which is partially coloured by the $(t,x)$-centred partial colouring using $k\geq 3$ colours. Since $\Gcal$ contains no gaudy obstruction, this partial colouring is not gaudy.
    Since $A$ is a last upleaf, $G_{i+2}$ has no upleaves  rooted in it and so all of the $2$-components at $2$-separators of $G_{i+2}$ that cross $b$ contain some coloured edge.
\end{cproof}

\begin{claim}\label{mnotboring}
There exists a list $M_0$ of $H_{i+2}$ so that any given embedding of $H_{i+2}$ satisfies $M_0$ if and only if it satisfies all the lists $M_1,\ldots,M_l$.
\end{claim}
\begin{cproof}
    Since $\Gcal$ contains no gaudy obstruction, $(G_i+2,b)$ is not gaudy. So by \autoref{lem:respect_characterisation} there is a bond-respectful embedding of $(G_i+2,b)$; call it $\iota$. Since $(H_{i+2},b)$ is a model of $(G_i+2,b)$, there is a bond-respectful embedding $\kappa$ of $(H_i+2,b)$ that induces $\iota$ on $G_{i+2}$. Hence $\kappa$ satisfies the lists $M_1,\ldots,M_l$.
    Hence, by \autoref{lem:sim_gives_combinable}, the lists $M_1,\ldots,M_l$ are combinable. From \autoref{lem:list-combining_machine} we obtain a new list $M_0$ so that any given embedding of $H_{i+2}$ satisfies $M_0$ if and only if it satisfies all the lists $M_1,\ldots,M_l$.
\end{cproof}
    
    Let\vspace{-.25cm}
    \[
        \Hcal:=(G_1\minorm\ldots\minorm G_{i-1}\minorp G_i\minorm G_{i+1}\minorp H_{i+2}\minorm G_{i+3}\minorp\ldots\minorp G_n).\vspace{-.05cm}
    \]
\begin{claim}\label{analogue-claim}
There is a list $N$ for $\Hcal$ so that $(\Hcal,N)$ is an analogue for $(\Gcal,L)$.
\end{claim}
\begin{cproof}    
    Since $H_{i+2}$ is accurate for $G_{i+2}$, the list $L$ of $\Gcal$ lifts via \autoref{lem:list_lifter} to a list $N$ for $\Hcal$ so that each simultaneous embedding $\Sigma$ of $\Hcal$ satisfies $N$ if and only if the simultaneous embedding $\Sigma$  induces on $\Gcal$ satisfies $L$.

    By \autoref{lem:b-respect_from_simultaneous}, since $(H_{i+2},b)$ is a model for $(G_{i+2},b)$, we get that every simultaneous embedding of $\Gcal$ lifts to a simultaneous embedding of $\Hcal$. In this way $(\Hcal,N)$ is an analogue of $(\Gcal,L)$.
    \end{cproof}
    Let\vspace{-.25cm}
    \[
        \Gcal':=(G_1\minorm\ldots\minorm G_{i-1}\minorp G_i-d\minorm G_{i+1}-d\minorp H_{i+2}\minorm G_{i+3}\minorp\ldots\minorp G_n).\vspace{-.05cm}
    \]
\begin{claim}\label{crazy}
There is a list $M$ of $\Gcal'$ so that $(\Gcal',M)$ is an analogue of~$(\Hcal,\varnothing)$.
\end{claim}
    \begin{cproof}
    We denote by $M$ the extension of the list $M_0$ for $H_{i+2}$ to $\Gcal'$ in the trivial way; that is, we include all Tutte-equivalence classes not from $H_{i+2}$ as single classes. 
    
    Let a simultaneous list-embedding $(\phi_1,\ldots,\phi_n)$ for $(\Hcal,\varnothing)$ be given.  
Since $(\phi_1,\ldots,\phi_n)$ is a simultaneous embedding, by \autoref{lem:b-respect_from_simultaneous} $\phi_{i+2}$ is $b$-respectful. By \autoref{boringstuff1}, $\phi_{i+2}$ satisfies the lists $M_1,\ldots,M_l$.
Hence by \autoref{mnotboring}, $\psi_{i+2}$ satisfies the list $M_0$. 
    Note that $\Hcal$ and $\Gcal'$ have the same $3$-blocks, and thus $M$ is a list for $\Hcal$. 
    Hence $(\phi_1,\ldots,\phi_n)$ satisfies $M$. 
So the simultaneous embedding that $(\phi_1,\ldots,\phi_n)$ induces on $\Gcal'$ is a list-embedding for $(\Gcal',M)$.

Conversely, let a simultaneous list-embedding $(\psi_1,\ldots,\psi_n)$ for $(\Gcal',M)$ be given. 
Thus $\psi_{i+2}$ satisfies $M_0$. By \autoref{mnotboring} $\psi_{i+2}$ satisfies the lists $M_1,\ldots,M_l$.
Hence by \autoref{boringstuff1}, $\psi_{i+2}$ is $b$-respectful.
Since $k\geq 3$, we can apply \autoref{lem:put_in_edge} to deduce that there is a simultaneous embedding $(\phi_i,\phi_{i+1},\phi_{i+2})$
of $(G_i\minorp G_{i+1} \minorm H_{i+2})$ that induces $(\psi_i,\psi_{i+1},\psi_{i+2})$ on $(G_i-d\minorp G_{i+1}-d \minorm H_{i+2})$. Then $(\psi_1,\ldots, \psi_{i-1},\phi_i,\phi_{i+1},\phi_{i+2}, \psi_{i+3},\ldots,\psi_n)$ is 
a simultaneous embedding of $\Hcal$; and thus a simultaneous list-embedding of~$(\Hcal,\varnothing)$. 

\end{cproof}

Since both $\Hcal$ and $\Gcal'$ have the same $3$-blocks, the list $N$ for $\Hcal$ is also a list for $\Gcal'$.
By \autoref{analogue-claim} and \autoref{crazy}, $(\Gcal',\{M,N\})$ is an analogue for $(\Gcal,N)$. 
    Furthermore, since the number of colours in the $(t,x)$-centred colouring of $G_{i+2}$ is bounded by the number of $2$-components at $t$ in $G_i$, we have that $k\leq \Cbb(G_{i+2})$. Hence
$
        |E(H_{i+2})| \leq |E(G_{i+2})|+12\,(\Cbb(G_{i+2}))^2
$.
\end{pf}

\subsection{Befitting models}
\label{sec:bond_fitting}
In this subsection and the following, we restrict our attention to planar graphs, without the context of temporal sequences. 
The goal of this subsection is to prove \autoref{lem:b-model} stated above. For this, we prepare as follows.

\begin{figure}
    \centering
        {\Large\scalebox{0.7}{\begin{tikzpicture}
	\begin{pgfonlayer}{nodelayer}
		\node [style=none] (341) at (6, 3.75) {$x$};
		\node [style=none] (342) at (6, -3.75) {$y$};
		\node [style=none] (476) at (6, -4.75) {$G$};
		\node [style=black dot] (479) at (6, 2.75) {};
		\node [style=black dot] (480) at (6, -2.75) {};
		\node [style=black dot] (481) at (6, -2) {};
		\node [style=black dot] (482) at (5.25, -3.25) {};
		\node [style=black dot] (485) at (5.25, 3.25) {};
		\node [style=black dot] (486) at (6, 2) {};
		\node [style=black dot] (488) at (6.75, 3.25) {};
		\node [style=black dot] (489) at (6.75, -3.25) {};
		\node [style=none] (491) at (9, 0.25) {};
		\node [style=none] (492) at (9, -0.25) {};
		\node [style=none] (493) at (9, 0) {};
		\node [style=none] (494) at (10.5, 0) {};
		\node [style=none] (495) at (10.25, 0.25) {};
		\node [style=none] (496) at (10.25, -0.25) {};
		\node [style=none] (547) at (13.5, 3.75) {$x_1$};
		\node [style=none] (548) at (14, -3.75) {$y$};
		\node [style=black dot] (551) at (14, -2) {};
		\node [style=black dot] (552) at (13.25, -3.25) {};
		\node [style=black dot] (553) at (12.5, 3.5) {};
		\node [style=black dot] (554) at (14, 1.25) {};
		\node [style=black dot] (555) at (15.5, 3.5) {};
		\node [style=black dot] (556) at (14.75, -3.25) {};
		\node [style=black dot] (558) at (14, -2.75) {};
		\node [style=black dot] (560) at (14, 2) {};
		\node [style=black dot] (561) at (13.25, 3) {};
		\node [style=black dot] (562) at (14.75, 3) {};
		\node [style=none] (563) at (14, -4.75) {$H$};
		\node [style=none] (564) at (24.5, 3.75) {$x$};
		\node [style=none] (565) at (24.5, -3.75) {$y$};
		\node [style=none] (566) at (24.5, -4.75) {$G$};
		\node [style=black dot] (567) at (24.5, 2.75) {};
		\node [style=black dot] (568) at (24.5, -2.75) {};
		\node [style=black dot] (569) at (23.75, -2) {};
		\node [style=black dot] (570) at (23.75, -3.5) {};
		\node [style=black dot] (571) at (23.75, 3.5) {};
		\node [style=black dot] (572) at (23.75, 2) {};
		\node [style=black dot] (573) at (25.25, 3.5) {};
		\node [style=black dot] (574) at (25.25, -3.5) {};
		\node [style=black dot] (575) at (25.25, 2) {};
		\node [style=black dot] (576) at (25.25, -2) {};
		\node [style=none] (577) at (14.5, 3.75) {$x_2$};
		\node [style=none] (578) at (14.75, 1.75) {$x_3$};
		\node [style=none] (579) at (27.75, 0.25) {};
		\node [style=none] (580) at (27.75, -0.25) {};
		\node [style=none] (581) at (27.75, 0) {};
		\node [style=none] (582) at (29.25, 0) {};
		\node [style=none] (583) at (29, 0.25) {};
		\node [style=none] (584) at (29, -0.25) {};
		\node [style=none] (586) at (33.25, -3.75) {$y$};
		\node [style=none] (587) at (33.25, -4.75) {$H$};
		\node [style=black dot] (589) at (33.25, -2.75) {};
		\node [style=black dot] (590) at (32.5, -2) {};
		\node [style=black dot] (591) at (32.5, -3.5) {};
		\node [style=black dot] (592) at (32, 4) {};
		\node [style=black dot] (593) at (32, 1.5) {};
		\node [style=black dot] (594) at (34.5, 4) {};
		\node [style=black dot] (595) at (34, -3.5) {};
		\node [style=black dot] (596) at (34.5, 1.5) {};
		\node [style=black dot] (597) at (34, -2) {};
		\node [style=black dot] (598) at (32.5, 3.5) {};
		\node [style=black dot] (599) at (32.5, 2) {};
		\node [style=black dot] (600) at (34, 2) {};
		\node [style=black dot] (601) at (34, 3.5) {};
		\node [style=none] (602) at (32.75, 4.25) {$x_1$};
		\node [style=none] (603) at (33.75, 4.25) {$x_2$};
		\node [style=none] (604) at (32.75, 1.25) {$x_3$};
		\node [style=none] (605) at (33.75, 1.25) {$x_4$};
	\end{pgfonlayer}
	\begin{pgfonlayer}{edgelayer}
		\draw [style=c0] (479) to (485);
		\draw [style=c0] (479) to (486);
		\draw [style=c0] (488) to (479);
		\draw [style=c0] (481) to (480);
		\draw [style=c0] (480) to (482);
		\draw [style=c0] (480) to (489);
		\draw [style=c0, bend left=15] (486) to (481);
		\draw [style=c0, bend left=45] (482) to (485);
		\draw [style=c0, bend right=60, looseness=0.75] (489) to (488);
		\draw [style=c0, bend left=45, looseness=0.50] (488) to (489);
		\draw [style=c0, bend right=15] (486) to (481);
		\draw [style=c0, bend left, looseness=0.75] (482) to (485);
		\draw [style=c0] (491.center) to (492.center);
		\draw [style=c0] (493.center) to (494.center);
		\draw [style=c0] (494.center) to (495.center);
		\draw [style=c0] (494.center) to (496.center);
		\draw [style=c0, bend left=15] (554) to (551);
		\draw [style=c0, bend left=45] (552) to (553);
		\draw [style=c0, bend right=60, looseness=0.75] (556) to (555);
		\draw [style=c0, bend left=45, looseness=0.50] (555) to (556);
		\draw [style=c0, bend right=15] (554) to (551);
		\draw [style=c0, bend left, looseness=0.75] (552) to (553);
		\draw [style=c0] (553) to (561);
		\draw [style=c0] (561) to (560);
		\draw [style=c0] (560) to (554);
		\draw [style=c0] (560) to (562);
		\draw [style=c0] (562) to (555);
		\draw [style=c0] (561) to (562);
		\draw [style=c0] (551) to (558);
		\draw [style=c0] (558) to (552);
		\draw [style=c0] (558) to (556);
		\draw [style=c0] (567) to (571);
		\draw [style=c0] (567) to (572);
		\draw [style=c0] (573) to (567);
		\draw [style=c0] (569) to (568);
		\draw [style=c0] (568) to (570);
		\draw [style=c0] (568) to (574);
		\draw [style=c2, bend left=15] (572) to (569);
		\draw [style=c1, bend left=45] (570) to (571);
		\draw [style=c1, bend right=45] (574) to (573);
		\draw [style=c1, bend left=45, looseness=0.50] (573) to (574);
		\draw [style=c2, bend right=15] (572) to (569);
		\draw [style=c2, bend left, looseness=0.75] (570) to (571);
		\draw [style=c0] (568) to (576);
		\draw [style=c1, bend right=15, looseness=0.75] (576) to (575);
		\draw [style=c2, bend right=15] (575) to (576);
		\draw [style=c0] (575) to (567);
		\draw [style=c0] (579.center) to (580.center);
		\draw [style=c0] (581.center) to (582.center);
		\draw [style=c0] (582.center) to (583.center);
		\draw [style=c0] (582.center) to (584.center);
		\draw [style=c0] (590) to (589);
		\draw [style=c0] (589) to (591);
		\draw [style=c0] (589) to (595);
		\draw [style=c2, bend right=15, looseness=0.75] (593) to (590);
		\draw [style=c1, bend left=75] (591) to (592);
		\draw [style=c1, bend right=75] (595) to (594);
		\draw [style=c1, bend left=75, looseness=0.75] (594) to (595);
		\draw [style=c2, bend right=45] (593) to (590);
		\draw [style=c2, bend left=60, looseness=0.75] (591) to (592);
		\draw [style=c0] (589) to (597);
		\draw [style=c1, bend right=45] (597) to (596);
		\draw [style=c2, bend left=15, looseness=0.75] (596) to (597);
		\draw [style=c0] (598) to (599);
		\draw [style=c0] (599) to (600);
		\draw [style=c0] (600) to (601);
		\draw [style=c0] (601) to (598);
		\draw [style=c0] (592) to (598);
		\draw [style=c0] (599) to (593);
		\draw [style=c0] (600) to (596);
		\draw [style=c0] (601) to (594);
	\end{pgfonlayer}
\end{tikzpicture}}}
    
        \caption{Coaddition of cycles at a vertex $x$ in a tuplet.}
    \label{fig:coadd_cycles}
\end{figure}
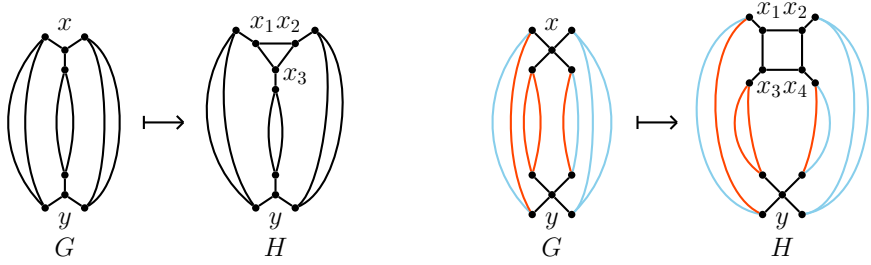

Let $(G,b)$ be a bond graph along with a partial \rcol/\bcol-colouring. We say that a tuplet $t$ of $G$ is \emph{almost fit} if it is fit or there is at most one monochromatic component at $t$ in each colour. 
We say that $(G,b)$ is \emph{almost bond-fitting} if every tuplet of $G$ that crosses~$b$ is almost fit.

Roughly speaking, the next lemma tells us how to improve bond graphs from almost bond-fitting to bond-fitting.
\begin{lem}\label{lem:almost_to_fit}
    Let $(G,b)$ be a bond graph along with a partial \rcol/\bcol-colouring. Suppose that $(G,b)$ is nice and not badly-coloured and almost bond-fitting. Then there exists a nice accurate model $(H,b)$ of $(G,b)$ which is bond-fitting and not badly-coloured.
\end{lem}
\begin{pf}
If every tuplet of $G$ that crosses~$b$ is fit, then we are done. So assume there is such a tuplet that is not fit; call it $t$. By assumption, it is almost fit.
So there is at most one monochromatic component at $t$ in each colour. 
Since $(G,b)$ is not badly coloured, each component at $t$ has at least one colour and at most two of them are polychromatic. Since $t$ is a tuplet, there are at least three components at $t$. 
Thus there are either exactly three components at $t$ or else there are exactly four components at $t$, two of them polychromatic, one monochromatic in \rcol\ and one monochromatic in \bcol.

Next we construct a graph $H$ that has $G$ as a minor. We shall then show that $(H,b)$ is a nice accurate model which is bond-fitting. For each unfit tuplet $t$ of $G$ denote by $c_1,\ldots,c_n$ its $2$-components. When $t$ has exactly four $2$-components, let $c_1$ and $c_3$ be those which are monochromatic. Let $k$ denote the number of $2$-components at $t$. For each vertex $x$ in $t$, let Let $C_k$ be a $k$-cycle, denote its vertices by $x_1,\ldots,x_k$ (in its cyclic ordering) and $\theta_x$ be the function mapping the unique edge incident with $x$ in $c_i$ to the vertex $x_i$. Now let $H$ be the graph obtained from $G$ by coadding $C_k$ with respect to $\theta_x$ for each $x$ in an unfit tuplet $t$ of $G$. The two cases, $k=3$ and $k=4$, are illustrated on the left and right of \autoref{fig:coadd_cycles}, respectively.

Since coadding edges preserves bonds\footnote{This is easy to see. For example notice that bonds are cociruits in the cycle matroid of a graph. So by matroid duality this statement is equivalent to the statement that adding edges preserves circuits.}, $b$ is a bond of $H$.

\begin{claim}\label{cl0}
    Every tuplet of $H$ that crosses $b$ is a fit tuplet of $G$ which crosses $b$.
\end{claim}
\begin{cproof}
    This claim follows via repeated applications of \autoref{lem:coadd_tuplet_correspondence}.
\end{cproof}

\begin{claim}\label{cl1}
    The bond graph $(H,b)$ is bond-fitting.
\end{claim}
\begin{cproof}
    This claim is a direct corollary to \autoref{cl0}. Indeed, every tuplet $t$ of $H$ which cross $b$ is also a tuplet of $G$. By construction, the components at $t$ in $H$ are in bijection with the components at $t$ in $G$ so that their colours are preserved.
\end{cproof}

\begin{claim}\label{cl2}
The graph $H$ is nice.
\end{claim}
\begin{cproof}
    By \autoref{lem:2-con_lifts}, $H$ is $2$-connected. We apply \autoref{cl0} to get that each $2$-component at a tuplet attaches to each vertex of the tuplet by a unique edge. 
\end{cproof}

\begin{claim}\label{cl4}
The bond graph $(H,b)$ is a model of $(G,b)$.
\end{claim}
\begin{cproof}
    We are required to show that every bond-respectful embedding of $(G,b)$ lifts to a bond-respectful embedding of $(H,b)$. By induction on the number of coadditions performed, we assume that $H$ is obtained from $G$ by the single coaddition of a $k$-cycle $C_k$ at an unfit tuplet $t=\{x,y\}$ of $G$.

    We claim that in every bond-respectful embedding of $G$ the cyclic structure of the $2$-components $c_1,\ldots,c_k$ at $t$ is $(c_1,\ldots,c_k)$. For $k=3$ this is immediate since there is a unique cyclic structure on three elements. For $k=4$ observe that by \autoref{lem:bond_segment}, the edges from $b$ in each of the four $2$-components form segments in the cyclic order which $\iota$ induces on the $2$-components at $t$. Hence they must alternate between polychromatic and monochromatic in this cyclic order.

    Now by \autoref{lem:coadd_cycle}, every bond-respectful embedding of $(G,b)$ lifts to an embedding of $(H,b)$.
\end{cproof}

By construction and by \autoref{cl0}, the $2$-components at a tuplet of $H$ that cross $b$ are coloured the same as their counterparts in $G$. Hence $H$ is not badly-coloured. Now combining \autoref{cl1}, \autoref{cl2}, \autoref{cl4} and \autoref{lem:coadd_accurate} completes the proof. 
\end{pf}

\begin{lem}\label{lem:make_almost_fit}
    Let $(G,b)$ be a bond graph along with a partial \rcol/\bcol-colouring. Suppose that $(G,b)$ is nice and not badly-coloured. Then there exists a nice accurate model $(H,b)$ of $(G,b)$ which is almost fit and not badly-coloured.
\end{lem}
\begin{pf}
    Since otherwise we are done, assume that some edge of $b$ receives colour \rcol\ and some edge receives colour \bcol.
    Our strategy will be to coadd a path $J_x$ at each vertex $x$ of a tuplet $t$ of $G$. This path $J_x$ depends on the types of $2$-components at $t$ and is defined according to the following cases.
    
    \textbf{Case 1.} The tuplet $t$ does not cross $b$ or has at most one monochromatic $2$-component in each colour. Then let $J_x$ be the single vertex $x_0$.
    
    \textbf{Case 2.} The tuplet $t$ crosses $b$ and has more than one monochromatic \rcol\ $2$-component and more than one monochromatic \bcol\ $2$-component. Then let $J_x$ be the path $x_R x_r x_0 x_b x_B$.
    
    \textbf{Case 3.} The tuplet $t$ has more than one monochromatic \rcol\ $2$-component but at most one monochromatic \bcol\ $2$-component. Then let $J_x$ be the path $x_R x_r x_0$.
    
    \textbf{Case 4.} The tuplet $t$ has at most one monochromatic \rcol\ $2$-component but has more than one monochromatic \bcol\ $2$-component. Then let $J_x$ be the path $x_0 x_b x_B$.
    
    Now for each vertex $x$ in a tuplet of $G$, let $\theta_x$ be the function from the set of edges of $G$ incident with $x$ to the set of vertices of $J_x$, assigning the unique edge $e$ incident with $x$ in a $2$-component $c$ at $t$ to $x_0,x_R,x_B$ subject to the following rules: If $t$ does not cross $b$ or $c$ is not monochromatic then we assign $e$ to $x_0$. If $c$ is monochromatic but is the only monochromatic component in its colour, then we also assign $e$ to $x_0$. Otherwise $c$ is monochromatic and not the only monochromatic $2$-component at $t$ in its colour. If $c$ is monochromatic in \rcol\ we assign $e$ to $x_R$ and if $c$ is monochromatic in \bcol\ we assign $e$ to $x_B$.

    Let $H$ be the graph obtained from $G$ by coadding $J_x$ with respect to $\theta_x$ for each vertex $x$ in a tuplet of $G$.
    Since coadding edges preserves bonds, $b$ is a bond of $H$. We make the following two structural claims.

    \begin{claim}\label{claim:wobbly_2-con}
        The graph $H$ is $2$-connected.
    \end{claim}
    \begin{cproof}
        Recall that we assumed every tuplet of $G$ has both a \rcol\ $2$-component and a \bcol\ $2$-component. This ensures that, in our construction, the ends of each path $J_x$ meet some edge in $H$ that was incident with $x$ in $G$. Hence by \autoref{lem:2-con_lifts}, the graph $H$ is $2$-connected.
    \end{cproof}
    
    \begin{claim}\label{claim:wobbly_tuplet_struct}
        Every tuplet of $H$ is of the form $\{x_0,y_0\}$, $\{x_R,y_R\}$ or $\{x_B,y_B\}$ for some tuplet $\{x,y\}$ of $G$.
    \end{claim}  
    \begin{cproof}
        Let $t$ be a tuplet of $H$. By repeated application of \autoref{lem:tuplets_of_tree_coadd}, we conclude that $t$ is of the form $\{x_P,y_Q\}$ for $P,Q\in\{0,R,B\}$ and some tuplet $\{x,y\}$ of $G$. Suppose for a contradiction that $P\neq Q$. Both $\{x_P,y_Q\}$ and $\{x_Q,y_P\}$ are $2$-separators of $H$, however they cross and thus by cannot have more than two $2$-components each, as tuplets are always totally nested. This completes the proof of our claim.
    \end{cproof}

    Observe that all of the $2$-components at a tuplet of the form $\{x_R,y_R\}$ are monochromatic in \rcol\ except for one, and hence such tuplets are fit. Similarly, tuplets of the form $\{x_B,y_B\}$ are fit. Furthermore, tuplets of the form $\{x_0,y_0\}$ have at most one monochromatic \rcol\ $2$-component and at most one monochromatic \bcol\ $2$-component, and so such tuplets are almost fit.
    It remains to prove that the bond graph $(H,b)$ is a nice accurate model of $(G,b)$ that is not badly-coloured.

    \begin{claim}\label{claim:wobbly_model}
        The bond graph $(H,b)$ is a model of $(G,b)$.
    \end{claim}
    \begin{cproof}
        Let $\iota$ be an embedding of $G$. For each tuplet $t$ of $G$ and $2$-component $c$ at $t$, we have by \autoref{lem:bond_segment} that the edges of $b$ in $c$ form a segment in the cyclic order induced on $b$ by $\iota$. Hence the \rcol\ monochromatic $2$-components at $t$ form a segment in the cyclic structure induced on the $2$-components at $t$. Similarly, the \bcol\ monochromatic $2$-components at $t$ form a segment. Let $\iota'$ be the embedding obtained from $\iota$ by performing the following operations at each tuplet $t$ of $G$: If we are in Case 1 of the construction, do nothing. Else, we first subdivide each edge incident with $x$ in $G$ so that the subdivision vertices lie on the boundary of a small disk around $x$ and then delete $x$ and its incident edges.
        
        Now we will obtain an embedding $\iota''$ of $H$ from $\iota'$ by performing some contractions of the plane, identifying the subdivision vertices into the vertices $x_0,x_R,x_B$. Since the other cases are analogous, suppose we are in Case 2 of the construction.
        Observe that the subdivision vertices appear in the same cyclic structure in the boundary of the small disk as their corresponding $2$-components at $t$ in $G$. Now we know that all the subdivision vertices from monochromatic \rcol\ $2$-components form a segment and all of the subdivision vertices from monochromatic \bcol\ $2$-components form a segment. Hence by contracting appropriate parts of the small disk, we can make the necessary identifications. Let this new embedding of $H$ be denoted $\iota''$.
        By construction $\iota''$ induces $\iota$ on $G$, and this completes the proof of our claim.
    \end{cproof}

    \begin{claim}\label{claim:wobbly_nice}
        The graph $H$ is nice.
    \end{claim}
    \begin{cproof}
        By \autoref{claim:wobbly_tuplet_struct}, each $2$-component at a tuplet of $G$ is attached to each vertex of the tuplet by a unique edge. Along with \autoref{claim:wobbly_2-con}, this completes the proof of our claim.
    \end{cproof}
    
    \begin{claim}\label{claim:wobbly_accurate}
        The graph $H$ is accurate for $G$.
    \end{claim}
    \begin{cproof}
        Let $c$ be a $3$-block of $H$ and $b$ and $b'$ two $3$-blocks of $G$ that are minors of $c$. We are required to show that $b$ and $b'$ are torsos of the same bag of the Tutte decomposition of $G$. Let $B$ and $B'$ denote the bags of the Tutte decomposition of $G$ corresponding to $b$ and~$b'$, respectively. By \autoref{claim:wobbly_tuplet_struct}, if $B$ and $B'$ can be distinguished by some $2$-separator of $G$, they can be separated by a pair of vertices in $H$. But this cannot happen since both $B$ and $B'$ (after possibly relabeling $x$ to some $x_i$ in each $B,B'$) are contained in the bag of the Tutte decomposition of $H$ corresponding to~$c$. Hence $B=B'$ and this completes our proof.
    \end{cproof}

    \begin{claim}\label{claim:wobbly_notbad}
        The graph bond graph $(H,b)$ is not badly-coloured.
    \end{claim}
    \begin{cproof}
        This claim follows from \autoref{claim:wobbly_tuplet_struct} and the way we constructed $H$. Indeed, each of the $2$-components at a tuplet $\{x_Q,y_Q\}$ in $H$ for some $Q\in\{0,R,B\}$ are coloured the same as their counterparts contained in them as minors at $\{x,y\}$ in $G$.
    \end{cproof}
    
    Combining \autoref{claim:wobbly_model}, \autoref{claim:wobbly_nice}, \autoref{claim:wobbly_accurate} and \autoref{claim:wobbly_notbad} completes the proof.
\end{pf}

\begin{pf}[Proof of \autoref{lem:b-model}]
    Let $(G,b)$ be a nice bond graph along with a partial \rcol/\bcol-colouring and suppose that $G$ is not badly-coloured. We are required to show that there exists a nice bond-fitting accurate model $(H,b)$ of $(G,b)$ which is not badly-coloured and furthermore, $\Cbb(H)\leq\Cbb(G)$ and $|E(H)|\leq|E(G)|+12\,\Cbb(G)$. 
    
    Firstly, we apply \autoref{lem:make_almost_fit} to obtain a nice accurate model $(L,b)$ of $(G,b)$ which is not badly-coloured, all whose tuplets which cross $b$ are almost fit. The application of \autoref{lem:make_almost_fit} requires replacing some tuplets with pairs of triangles or pairs of $4$-cycles, which meet three or four new $2$-separators with only two $2$-components, respectively. Every tuplet contributes $3=3\cdot (3-2)$ or $6=3\cdot (4-2)$ to $3\beta(G)$, respectively. Each totally-nested $2$-separator with only two $2$-components contributes $1$ to $\alpha(H)$. Hence $\alpha(G)+3\beta(G)\leq\alpha(H)+3\beta(G)$. Observing that $\gamma(H)=\gamma(G)$ proves $\Cbb(H)\leq\Cbb(G)$. Furthermore, we add at most $8$ edges per tuplet so $E(L)\leq E(G)+8\,\Cbb(G)$.
    
    From $(L,b)$, we obtain our required model $(H,b)$ by applying \autoref{lem:almost_to_fit}. Here we replace some tuplets with pairs or triples of tuplets, so that the total number of $2$-components at tuplets increases by $2$ or $4$ and the total number of tuplets increases by $1$ or $2$, respectively. Hence $\beta(H)\leq\beta(G)$. Observe that $\alpha(H)=\alpha(G)$ and $\gamma(H)=\gamma(G)$ proves that $\Cbb(H)\leq\Cbb(L)$. We perform the coaddition of at most four edges at each tuplet of $L$ to form $H$ and so $E(H)\leq |E(L)|+4\,\Cbb(L)$. In total, both $\Cbb(H)\leq \Cbb(G)$ and $|E(H)|\leq |E(G)|+12\, \Cbb(G)$, as required.
\end{pf}

\subsection{Delightful reductions}
\label{sec:del_red}

In this subsection we prove \autoref{list_for_bonds}, which is the last piece of our proof of \autoref{thm:main}.

Let $(G,b)$ be a $2$-connected bond graph and $p$ a bag of the Tutte decomposition of $G$. Let $t$ be a totally-nested $2$-separator of $G$ that lies in $p$. The \emph{bough} at $t$ for $p$ is the union of the $2$-components at $t$ which do not contain $p$. The \emph{bough segment} of $b$ at $t$ for $p$ consists of those edges of the bough at $t$ for $p$ that are in $b$. Let $P$ be the bond-torso of $p$. For an edge $e$ of $P$ that is not a torso edge, we call the singleton $\{e\}$ the \emph{bough segment} of $b$ at $e$ for $p$. In this way, the bough segments for $p$ partition $b$. By our definition of bond-torso, for each totally-nested $2$-separator $t$ of $G$ that lies in $p$ so that its bough segment is nonempty, one of the edges of its bough segment labels the torso edge of $P$ for $t$. 
This way each nonempty bough segment at a totally nested $2$-separator is mapped to an edge of $P$, and other bough segments are already edges of $P$. This gives a bijective map between the nonempty bough segments and the  restriction of $b$ to the bond-torso $P$; we refer to this bijection as the \emph{bough bijection}. 

\begin{lem}\label{segments_are_segments}
For every embedding $\iota$ of a $2$-connected bond graph $(G,b)$ and a bag $p$ of the Tutte decomposition of $G$, every bough segment for $p$ is a segment in the cyclic ordering $\iota$ induces on $b$.
\end{lem}
\begin{pf}
    For each totally-nested $2$-separator $t$ of $G$ in $p$, we apply \autoref{lem:bond_segment} to the $2$-component at $t$ which contains $p$ and observe that the complement of a segment is again a segment.
\end{pf}

In the context of \autoref{segments_are_segments}, we make the following definition: take the cyclic ordering $\iota$ induces on $b$, and collapse each nonempty bough segment to a singleton. We refer to the induced order of these singletons (resp.\ the bough segments) as the \emph{bough ordering induced by $\iota$ for $p$}; if the bough ordering contains at least two elements, then it is a cyclic ordering. 

\begin{lem}\label{compare_via_bough_ordering}
Let $p$ be a bag of the Tutte decomposition of a $2$-connected bond graph $(G,b)$.
If two embeddings of $(G,b)$ induce different bough orderings for $p$, then they induce different embedding on the torso of $p$. 
\end{lem}
\begin{pf}
Let $X$ be the graph with two vertices that has one edge between these two vertices for every bough at $p$. 
Denote the torso of $p$ by $Y$. 
Then $X$ can be obtained as a minor of $Y$ by contracting all edges outside $b$. 
The torso $Y$ in turn is a minor of $G$. Now if  two embeddings of $(G,b)$ induce different bough orderings for $p$, then they induce these same bough orderings on $X$ when considered as a minor of $G$ -- by contracting the edges outside $b$ and deleting all but one edge from every nonempty bough segment. Now this minoring goes via $Y$, and hence the two embeddings cannot agree on $Y$ since they differ on the minor $X$. Hence the two embeddings differ on $Y$, as desired. 
\end{pf}

Given a $2$-connected bond graph $(G,b)$ and two bags $p$ and $q$ of its Tutte decomposition with at least three vertices, the \emph{bough segment for $p$ containing $q$} is the bough segment for the unique bough for $p$ that contains $q$.

\begin{lem}\label{bough-segments-nested}
Given a $2$-connected graph $G$, and bags $p$ and $q$ of its Tutte decomposition with at least three vertices, the bough segment for $p$ containing $q$ includes all bough segments at $q$ except for the bough segment for $q$ containing $p$. 
\end{lem}
\begin{pf}
By properties of tree-decompositions, the bough for $p$ containing $q$ includes all boughs for $q$ except for the bough for $q$ containing $p$. Hence the lemma follows by applying the definition of bough segment. 
\end{pf}

Let $\iota$ be an embedding of a $2$-connected graph $G$. By \autoref{determine_embed}, the embedding $\iota$ is determined by the cyclic orderings of the $2$-components at tuplets of $G$ and by the orientations of the $3$-blocks of $G$. Given a $3$-block $p$ of $G$, we say that an embedding $\kappa$ is the \emph{embedding obtained from $\iota$ by reorienting $p$} if the embedding data for $\iota$ and $\kappa$ are identical except for the orientation of $p$, where they disagree.

\begin{lem}\label{lem:flip_segments}
    Let $(G,b)$ be a $2$-connected bond graph and $p$ a $3$-block of $G$. Let $\iota$ be an embedding of $(G,b)$ and let $\iota(b)=I_1I_2\ldots I_n$ be the cyclic order induced on $b$ by $\iota$, where $I_1,\ldots,I_n$ are the nonempty bough segments for $p$. Let $\kappa$ be the embedding obtained from $\iota$ be reorienting $p$. Then $\kappa(b)=I_nI_{n-1}\ldots I_1$.
\end{lem}
\begin{pf}
    Let $P$ be the bond-torso of the bag $p$. By the bough bijection, the edges in the restriction $b\restric_P$ of $b$ to $P$ correspond to the nonempty bough segments for $p$, and furthermore for any embedding $\phi$ of $G$, the order induced by $\phi$ on the restriction $b\restric_P$ corresponds to the bough ordering for $p$. Reorienting an embedding of $P$ reverses the cyclic order induced on $b\restric_P$. Hence reorienting $p$ in an embedding of $G$ reverses the bough order for $p$.

    It remains to show that the internal order of each bough segment is the same between $\iota$ and $\kappa$. For each bough at a $2$-separator of $G$ in $p$, let $G_t$ be the graph obtained from $G$ by minoring the $2$-component at $t$ containing $p$ down to a single edge $e_t$ (as in \autoref{2-sum-inverse}). For an embedding $\phi$ of $G$, let $\phi_t$ denote the embedding induced on $G_t$ by $\phi$.
    \begin{claim}
        We have $\iota_t=\kappa_t$.
    \end{claim}
    \begin{cproof}
        Since $\iota$ and $\kappa$ have the same embedding data except at $p$, we get that the embedding data for $\iota_t$ and $\kappa_t$ are identical and hence by \autoref{determine_embed}, $\iota_t=\kappa_t$.
    \end{cproof}
    
    Let $b_t$ be the restriction of $b$ to $G_t$ along with the edge $e_t$. Since $G$ can be obtained as a $2$-sum of $G_t$ with the bond-branch at $t$ containing $W$, and both sides of the $2$-sum contain edges from $b$, the set $b_t$ is a bond of $G_t$.
    \begin{claim}
        Let $\phi$ be an arbitrary embedding of $G$, let $J_1,\ldots,J_n$ be the orders induced by $\phi$ on each of the bough segments for $w$. Then the cyclic order induced by $\phi_t$ on $b_t$ is $e_t \cdot J_i$.
    \end{claim}
    \begin{cproof}
        This claim follows from the construction of $G_t$.
    \end{cproof}
    
    Combining the above two claims complete our proof.
\end{pf}

\begin{dfn}[Colour-locked]
    Let $G$ be a $2$-connected graph and $b$ a bond of $G$ along with a \rcol/\bcol-colouring of some of its edges. We say that a a $3$-block $p$ of $G$ is \emph{colour-locked} if its $b$-torso $P$ satisfies the following: the $b$-torso $P$ contains three coloured edges so that two of them receive the colour \rcol\ and two of them receive the colour $\bcol$.
\end{dfn}

\begin{lem}\label{lem:colour-locking}
    Let $(G,b)$ be a planar bond graph along with a partial \rcol/\bcol-colouring. Suppose that this partial colouring is not gaudy. Let $p$ and $q$ be colour-locked $3$-blocks of $G$. Then there exist embeddings $\iota$ and $\kappa$ for $p$ and $q$, respectively, so that every bond-respectful embedding of $G$ induces both $\iota$ and $\kappa$ or their reorientations, $(\iota)^-$ and $(\kappa)^-$.
\end{lem}
\begin{pf}
Since $p$ is colour-locked, there are bough segments $I_1$ and $I_2$ at $p$ that are different from the bough segment at $p$ containing $q$  so that they each contain a coloured edge and $I_1\cup I_2$ contains both \rcol\ and \bcol\ edges. 
If one of them contains edges of both colours, then one of them contains a \rcol\ edge and the other contains  a \bcol\ edge. Otherwise, we directly see that one of them contains a \rcol\ edge and the other contains a \bcol\ edge.
So by symmetry assume that $I_1$ contains a \rcol\ edge and $I_2$ contains a \bcol\ edge. 

Similarly, there are bough segments $I_3$ and $I_4$ at $q$ that are different from the bough segment at $q$ containing $p$ so that $I_3$ contains a \bcol\ edge and $I_4$ contains a \rcol\ edge. 
By \autoref{bough-segments-nested} the disjoint bough segments $I_1$ and $I_2$ are disjoint from the disjoint bough segments $I_3$ and $I_4$. Thus all the sets $I_1$, $I_2$, $I_3$ and $I_4$ are disjoint.

Since the partial colouring for $G$ is not gaudy, by \autoref{lem:existence_b-respect} the graph $G$ admits a $b$-respectful embedding; denote it by $\phi$. Let $\iota$ be the embedding that $\phi$ induces on the torso $p$ of $G$, which is a minor of $G$; and let $\kappa$ be the embedding that $\phi$ induces on the torso of $q$. 
Abbreviate $I=I_1\cup I_2\cup I_3\cup I_4$. 

\begin{claim}\label{proof_of structure}
The cyclic order on $b$ induced by $\phi$ restricted to $I$ has the form\footnote{In this lemma, we treat $I_i$ as a set without a given order. So unlike in \autoref{lem:flip_segments}, writing $I_1I_2I_3I_4$ does not specify the order completely.} $I_1I_2I_3I_4$ or its reverse $I_4I_3I_2I_1$.
\end{claim}

\begin{cproof}
Denote by $\sigma$ the cyclic order on $b$ induced by $\phi$ restricted to $I$. 
Since each set $I_j$ is a bough segment, the sets $I_j$ are segments of $\sigma$. And as shown above, the sets $I_j$ are mutually disjoint. 
So it remains to evaluate which orderings of these segments are possible. 

Since $\phi$ is $b$-respectful, the sets containing edges of the same colour must not be separated by sets containing edges of the opposite colour. Thus $I_2$ and $I_3$ must be adjacent, and $I_4$ and $I_1$ must be adjacent.
Since $I_3$ and $I_4$ are bough segments at $q$, and $I_1$ and $I_2$ are included in a single bough segment at $q$ (the bough segment containing $p$) they must be adjacent. 
Similarly,  $I_3$ and $I_4$ must be adjacent. Now we have gathered enough information, and we will figure out the possibilities in the following. 

If $\sigma$ contains $I_1$ followed by $I_2$, then $I_3$ must follow after that, and then $I_4$.
If $I_1$ is not followed by $I_2$, since they are adjacent, it must be that $I_2$ is followed by $I_1$, then $I_4$ comes after that and then $I_3$. Thus we get $\sigma=I_1I_2I_3I_4$ or $\sigma=I_2I_1I_4I_3$. 
This are the two outcomes of the claim for $\sigma$, after applying a cyclic permutation to the second term.
\end{cproof}

By \autoref{proof_of structure} the cyclic order on $b$ induced by $\phi$ restricted to $I$ has the form $I_1I_2I_3I_4$ or its reverse $I_4I_3I_2I_1$, by symmetry assume that it is $I_1I_2I_3I_4$.
Then for the reorientation of $\phi$,  the cyclic order on $b$ induced by it restricted to $I$ is $I_4I_3I_2I_1$.

Now let $\psi$ be an arbitrary $b$-respectful embedding of $G$. 
We distinguish two cases. 

{\bf Case 1:} the cyclic order on $b$ induced by $\psi$ restricted to $I$ has the form $I_1I_2I_3I_4$. 
So in the bough order of $\psi$ at $p$ we have that $I_2$ comes after $I_1$ and before the bough segment containing $q$, while in the bough order of the reorientation of $\phi$ at $p$ it is the other way round. 
Thus the bough order of $\psi$ at $p$ is not equal to the bough order of the reorientation of $\phi$ at $p$. Hence by \autoref{compare_via_bough_ordering} $\psi$ does not induce the same embedding on the torso of $p$ as the reorientation of $\phi$. Since the reorientation of $\phi$ induces the reorientation of $\iota$, the embedding $\psi$ induces $\iota$. An analogue argument for \lq $q$\rq\ in place of \lq $p$\rq\ gives that $\psi$ induces $\kappa$ on the torso of $q$. 

{\bf Case 2:} not Case 1. Then by \autoref{proof_of structure} applied to \lq $\psi$\rq\ in place of \lq $\phi$\rq\ the cyclic order on $b$ induced by $\psi$ restricted to $I$ has the form $I_4I_3I_2I_1$. 
Now we argue as in Case 1 with \lq $\phi$\rq\ in place of the \lq reorientation of $\phi$\rq. 
This completes the proof.
\end{pf}

\begin{dfn}
     Let $(G,b)$ be a bond graph along with a partial \rcol/\bcol-colouring. Let $\iota$ be a bond-respectful embedding for $(G,b)$. Then the \emph{list induced by $\iota$} is the list with one nontrivial part, consisting of all of the colour-locked $3$-blocks of $G$, equipped with the embeddings induced on the $3$-blocks of $G$ by $\iota$.
\end{dfn}

Recall that two lists are equivalent if their partitions are the same and one can be obtained from the other by reorienting all of the embeddings in some of the parts.
\begin{cor}\label{cor:respect_list_equivalence}
     Let $(G,b)$ be a bond graph along with a partial \rcol/\bcol-colouring. If this colouring is not gaudy, then any two bond-respectful embeddings induce equivalent lists.
\end{cor}
\begin{pf}[Proof:]
 apply \autoref{lem:colour-locking}.
\end{pf}

\begin{lem}\label{lem:induced_list_yehaw}
    Let $(G,b)$ be a bond graph along with a partial \rcol/\bcol-colouring. Let $\iota$ be a bond-respectful embedding of $G$ and $L$ be the list induced by $\iota$. If $G$ is bond-fitting, then every embedding of $G$ which satisfies $L$ is bond-respectful.
\end{lem}
\begin{pf}
Let $\kappa$ be an embedding of $G$ that satisfies $L$. 
If necessary, we replace $\kappa$ with its reorientation so that the embeddings it induces on the colour-locked $3$-blocks of $G$ coincide with those specified by $L$. By \autoref{determine_embed}, the embedding $\kappa$ is determined by the cyclic orderings of the $2$-components at tuplets and by the orientations of its $3$-blocks. We prove by induction on the number of disagreements between $\iota$ and $\kappa$—in terms of these cyclic orderings and orientations—that $\kappa$ is $b$-respectful. To do so, it suffices to show that a single such change preserves $b$-respectfulness. Thus, we may assume that $\kappa$ is obtained from $\iota$ by performing exactly one such modification.

    \textbf{Operation 1.} Suppose that $\kappa$ is obtained from $\iota$ by changing the cyclic order in which the $2$-components appear at a tuplet $t$. By \autoref{lem:bond_segment}, this reordering corresponds to permuting the segments determined by their intersections with the bond $b$, according to the cyclic order induced on $b$. Since $G$ is bond-fitting, all but at most one of the $2$-components at $t$ are monochromatic in a shared colour. It follows that $\kappa$ is bond-respectful.

    \textbf{Operation 2.} Suppose that $\kappa$ is obtained from $\iota$ by reorienting some $3$-block $p$ of $G$ which is not colour-locked. We distinguish two cases.

       \textbf{Case 1:} There is a colour, say \rcol, such that the bond-torso of $p$ contains at most one edge of colour \rcol. Hence applying the bough bijections yields that all but at most one bough segment for $p$ contain only \bcol\ edges. Thus by \autoref{lem:flip_segments} reorienting $p$ preserves $b$-respectfulness.

       \textbf{Case 2:} not Case 1; that is, the bond-torso of $p$ contains at least two edges of each colour. Since $p$ is not colour-locked, the bond-torso contains at most two coloured edges in total. 
       Applying the bough bijection gives that there are most two bough segments for $p$. Thus by \autoref{lem:flip_segments} reorienting $p$ preserves $b$-respectfulness. This completes the last case, and thus this proof.      
       \end{pf}

\begin{pf}[Proof of \autoref{list_for_bonds}]
    Let $(G,b)$ be a bond graph along with a partial \rcol/\bcol-colouring. Suppose that this colouring is not gaudy and that $G$ is bond-fitting. We are required to show that there exists a list $L$ on the $3$-blocks of $G$ so that an embedding of $G$ is bond-respectful if and only if it satisfies $L$.
    
    By \autoref{lem:existence_b-respect}, there exists a bond-respectful embedding $\iota$ of $G$. Now let $L$ be the list induced by $\iota$. By \autoref{lem:induced_list_yehaw}, every embedding of $G$ which satisfies $L$ is bond-respectful. Now let $\kappa$ be an arbitrary bond-respectful embedding of $G$. Let $M$ be the list induced by $\kappa$. By \autoref{cor:respect_list_equivalence}, the lists $M$ and $L$ are equivalent. Hence the embedding $\kappa$ satisfies $L$. Taking $L$ to be our desired list completes the proof.
\end{pf}

\begin{pf}[Proof of \autoref{thm:main}]
In \autoref{sec:main} we proved that Lemmas \myref{lem:A}, \myref{lem:B_2}, \ref{lem:list-combining_machine} and \ref{lem:niceness_main} together imply \autoref{thm:main}. These four facts are proved in sections \ref{sec:amazing_analogues}, \ref{sec:delightful_analogues}, \ref{sec:obs} and \ref{sec:nice}, respectively, concluding the proof. 
\end{pf}

\begin{pf}[Proof of \autoref{thm:algorithmic_main}]
This proof follows the structure of \autoref{thm:main}, augmented by a running time analysis.
Let a nonempty temporal sequence $(\Gcal,\{L,M\})$ of $2$-connected graphs with lists be given. Our goal is to either find an obstruction or construct a simultaneous embedding in polynomial time.
The first paragraph of the proof \autoref{thm:main} applies verbatim, and shows that we may without loss of generality assume that $\Gcal$ is nice and that we only have a set of lists consisting of a single list, which is called $N$. 
All this is done in a pre-processing step, which increases $\Gcal$ only by a polynomial amount in polynomial time. 
Next we find in polynomial time a last upleaf $A$, rooted at $i$ for some odd $i\in[n]$ as explained in that proof.

    We do not construct a simultaneous embedding explicitly, but instead use a recursive approach based on a parameter we call the flexiblity of $\Gcal$. By definition, the initial flexibility is polynomial in the size of $\Gcal$. At each step, we reduce the flexibility by at least one, but may increase the size of $\Gcal$ by an additive term that is polynomial in the current flexibility.
    More precisely, if we move from a temporal sequence $\Gcal$ to a new sequence $\Gcal'$ of lower flexibility, we ensure that $\text{size}(\Gcal') \leq \text{size}(\Gcal) + 13\cdot(\Cbb(\Gcal))^2$, 
    as detailed in the definition of an improvement (\autoref{dfn:improv}).
    Since the flexibility decreases in each step, the additive increase in size is bounded by a constant depending only on the initial size of $\Gcal$. It therefore suffices to verify that each such improvement step can be performed in polynomial time in the size of the current temporal sequence $\Gcal$.

Each improvement step falls into one of four cases, corresponding to applications of \autoref{lem:B_2} (\nameref{lem:B_2}), \autoref{claimX}, \autoref{unamazing_lemma}, or \myref{lem:A}.
All four cases can be handled in polynomial time. For \autoref{unamazing_lemma}, this is immediate. For \autoref{claimX}, it is also immediate, except for a single invocation of \myref{lem:A} within its proof.
It therefore suffices to verify that \myref{lem:A} and \myref{lem:B_2} can each be executed in polynomial time. These are proved in \autoref{sec:amazing_analogues} and \autoref{sec:delightful_analogues}, respectively. Although these sections are fairly long, it is straightforward to check that all the individual steps and lemmas they contain can be carried out in polynomial time.
\end{pf}

\section{Concluding remarks}
\label{sec:conc}

We conclude this paper by first explaining how the results of this paper can be applied in order to answer graph embedding problems that arise in different contexts, and secondly discuss new directions opened by this work. 

\subsection{Rooted-tree SEFE}
The \emph{Sunflower SEFE} problem has as input a family $G_1$, $G_2$, \ldots ,$G_k$ of graphs that all have a common subgraph $H$. The aim is to decide whether for all $i\in [k]$ the graphs $G_i$ admit embeddings that all induce the same embedding on $H$. We shall talk about the \emph{$k$-Sunflower SEFE} problem if the number $k$ of graphs is fixed. 
Note that $(G_1,G_2,\ldots, G_k,H)$ is an instance for $k$-Sunflower SEFE if and only if 
$(G_1,H,G_2,H,\ldots, H, G_k)$ is a simultaneously embeddable temporal sequence. 
By a result of \cite{gassner2006simultaneous}, 3-Sunflower SEFE is NP-hard to decide.
However, the corresponding construction makes use of the fact that the graph $H$ is not $2$-connected 
(while the first proof of this fact makes use of a disconnected graph $H$ \cite{gassner2006simultaneous}, later the constructions were improved relying only on a graph $H$ that is a tree \cite{ANGELINI201571}). By \cite{haeupler2010testing}, if all involved graphs are $2$-connected, then Sunflower SEFE is linear-time solvable. 

The Sunflower SEFE problem, may be generalised by replacing the star structure of the graphs $G_1,\ldots,G_k,H$ by the structure of a genuine rooted tree; that is, we are given a rooted tree $T$ and a graph $G_t$ assigned to each node $t$ of $T$ so that if for an edge $vw$ the endvertex $v$ is nearer than $w$ to the root, then $G_v\se G_w$. We refer to this problem as the rooted-tree SEFE problem. 

\begin{lem}\label{rooted_tree_sefe}
Rooted-tree SEFE problems $(T, (G_t|t\in V(T)))$ can be reduced to simultaneous embedding problems all whose graphs are in $(G_t|t\in V(T))$.
\end{lem}

\begin{pf}
Denote the root by $r$. Assume by induction that for each component $S$ of $T-r$ the rooted-tree SEFE problem $(S, (G_s|s\in V(S)))$ can be reduced to simultaneous embedding problems all whose graphs are in $(G_s|s\in V(S)))$.
Let $T_1,\ldots,T_k$ be the corresponding temporal sequences.
Now consider the temporal sequence $(T_1,G_r,T_2,G_r,\ldots,G_r,T_k)$. Clearly, it is simultaneous embeddable if and only if the rooted-tree SEFE problem $(T, (G_t|t\in V(T)))$ is satisfiable. 
\end{pf}

\subsection{Open problems and conjectures}

\autoref{thm:algorithmic_main} shows that the graph minors approach provides a systematic toolkit for studying temporal graphs, and opens the door to a graph minors theory for temporal graphs. A natural next step is to determine the 2-connected excluded minors for simultaneous embeddability of 2-connected temporal sequences. The most promising approach appears to be identifying the minimal members of each of our five classes of obstructions. We expect this set of excluded minors to be well-structured, in sharp contrast to the case without connectivity assumptions where we believe the list to be infinite and \lq wild\rq. Both of these directions are part of our ongoing work. A partial result in this direction is the complete characterisation of the (finitely many) excluded minors for simultaneous embeddability of 2-connected temporal contraction-only sequences in the upcoming work \cite{contractionsequences2025}.

A central question seems to be the following. 

\begin{oque}
Is there a polynomial time algorithm that checks whether a given temporal sequences is a temporal minor of another temporal sequence?
\end{oque}

We have studied planarity of temporal sequences. It is natural to ask whether similar results hold for other minor-closed properties. In particular, investigating an analogue of bounded tree-width might be the next step. Admitting a tree decomposition of width bounded by some given constant $w$ is a property like planarity, in the sense that such a tree decomposition enriches the structure of the graph and it yields a tree decompositions for all the minors of the graph with no increase in width or adhesion. Hence we can define a `simultaneous tree decomposition' and by extension `simultaneous treewidth'.

\begin{oque}\label{sim-tw}
Is there a polynomial-time algorithm that checks whether a temporal sequence has simultaneous bounded tree-width?
\end{oque}

We expect the answer to be as rich and intricate as in the planar case.

It is worth mentioning that much work has already done on the behaviour of tree decompositions on temporal graphs. For example, Korhonen, Majewski, Nadara, Pilipczuk, and Sokołowski devised a fully dynamic algorithm that maintains near-optimal tree decompositions and supports CMSO${}_2$ model-checking, yielding a dynamic analogue of Courcelle’s Theorem \cite{korhonen2023dynamic}. Complementary results include a dynamic MSO algorithm for graphs of bounded tree-depth by Dvořák, Kupec, and Tůma \cite{dvorak2014dynamic}, and dynamic rank-decompositions for LinCMSO${}_1$ problems by Korhonen and Sokołowski \cite{korhonen2024rankwidth}. Further developments include dynamic subgraph counting in sparse graphs by Dvořák and Tůma \cite{dvorak2013counting}, and approximation schemes for optimisation problems in planar and minor-free graphs by Korhonen, Nadara, Pilipczuk, and Sokołowski \cite{korhonen2024soda}. These works demonstrate the algorithmic reach of structural methods in dynamic settings, however they do not yet constitute a truly `temporal' version of treewidth that accounts for the whole sequence of graphs simultaneously. The study of true structural temporal graph width measures are still in its infancy, although there are rapid developments in this area \cite{Fluschnik2020,enright2025families}.

Beyond tree-width, many other structural properties arise naturally in the context of temporal sequences. These include bounded path-width, bounded adhesion, series-parallel graphs, outerplanar graphs, apex graphs, membership in a fixed minor-closed class, being a minor of a fixed (possibly infinite) graph \( H \) (such as the two-dimensional grid, a tessellation of the hyperbolic plane, or the Farey graph), and admitting a common structure as described by the Graph Minor Structure Theorem. In the latter case, it may be most natural to restrict attention to the $3$-connected setting. The most obvious question in this context is the following. 

\begin{oque}
Can you extend the results of this work from the plane to simultaneous embeddability in other surfaces, such as the torus?
\end{oque}

We view these problems as part of a broader research programme: extending the ideas and intuition of graph minor theory to temporal graphs. While the structural questions can often be phrased in a uniform way, we expect that their resolution will require a diverse range of new methods. 

Separately, in the course of this work, we developed the methods of \lq Stretching\rq\ and \lq befitting models\rq. We expect that both of them should be more widely applicable, and we believe that they are interesting enough to be studied in their own right. Another natural direction into which one can extend these ideas is to generalise the notion of temporal sequences even further. Indeed, a temporal sequence is a linearly ordered family of graphs, and one could more generally study tree-ordered families of graphs; this approach is introduced in the upcoming work
\cite{contractionsequences2025,fpttreesofgraphs2025}.
\begin{rem}\label{rem:weak_remark}
    A \emph{weak minor} of $G$ is a graph isomorphic to a minor of $G$ (alternatively, one can define \lq weak minors\rq\ directly using the definition of \lq minors\rq\ where one additionally allows relabelling of edges).
    In the paper \cite{weaktemporal2025}, we show that the results of this paper do not extend to the weak minor relation by proving that the $2$-connected simultaneous embedding problem for the weak minor relation is NP hard, even when we only allow edge-deletion throughout the sequence. Similarly, we show that the problem of determining simultaneous embeddability for $2$-connected (non-weak) temporal sequences is intractible when we relax the assumption that the type of minor operations (contraction or deletion) used in the temporal sequence are not fixed in advance. Since we believe that a minor theory should go hand in hand with algorithmic applications, this makes us believe that the minor relation we define in this paper is the correct one to develop this theory. 
\end{rem}
From a modelling perspective, planarity of temporal graphs is particularly relevant in applications such as brain networks, where wiring cost and physical embedding constrain long-range connections. Cortical connectivity has been observed to approximate nearly planar structures, reflecting a trade-off between cost-efficiency and functional integration \cite{bullmore2009complex, giusti2015clique}.

We conclude by stepping back and placing our work in a broader structural context. Consider two classes of objects equipped with a minor relation, where one class enriches the other, and the minor relation on the richer class induces the minor relation on the simpler one via a natural forgetting operation. In this paper, for instance, the richer structure is the class of graphs embedded in the plane, and the simpler structure is the class of graphs. Another example arises when graphs play the role of the richer structure and matroids form the simpler one.

This perspective suggests a possible algebraic reformulation of our theory. In particular, while the 2-connectivity assumption in \autoref{thm:algorithmic_main} is necessary in our setting, we speculate that it might be circumvented in a more algebraic context. For example:

\begin{oque}
Is there a polynomial-time algorithm that checks for temporal sequences of matroids whether they are simultaneously graphic?
\end{oque}

It is straightforward to see that \autoref{cor:poly} implies an affirmative answer for the subclass of 2-connected co-graphic matroids. 
Rutger Campbell suggested to us the following related question, which also falls under this algebraic umbrella:

\begin{oque}
Given a fixed finite field $k$, is there a polynomial-time algorithm that checks for temporal sequences of matroids whether they are simultaneously $k$-representable?
\end{oque}

We hope that this work lays the foundation for a temporal graph minor theory that unites topological intuition with algebraic structure, offering a systematic lens through which to explore evolving networks and their underlying structure.

\appendix
\begin{center}
{\LARGE\textbf{Appendix}}
\end{center}
\section{Algorithmic appendix}
\label{appendix:algorithms}

\begin{pf}[Proof of \autoref{lem:polynomial_disagreeable} for $2$-connected temporal sequences]
Let a $2$-connected temporal triple\\ $(G_1,G_2,G_3)$ be given. We are to find a polynomial-time algorithm that checks whether there is a pair $(H_1,H_3)$ so that $H_i$ is a $3$-connected minor of $G_i$ for $i=1,3$, and there exists a pair of temporal minors $(H_1\minorm\triple_1\minorp H_3)$ and $(H_1\minorm\triple_2\minorp H_3)$ of $\tseq{G}$, with $\triple_1$ and $\triple_2$ theta graphs, so that no choice of embeddings $\alpha_1$ for $H_1$ and $\alpha_3$ for $H_3$ satisfy both $\alpha_1\rightarrow\triple_1\leftarrow\alpha_3$ and $\alpha_1\rightarrow\triple_2\leftarrow\alpha_3$.

First, we notice that for this, we may assume that $H_1$ and $H_3$ are maximal $3$-connected minors of $G_1$ and $G_3$, respectively. Since the maximal $3$-connected minors are the $3$-blocks of the Tutte-decomposition, there are only linearly many choices for each of $H_1$ and $H_3$. A technical aspect is the labelling of the torso edges: 
given a $3$-block $b$, we label the torso edge $e$ representing the 2-separation $(A,B)$ of $G_i$ with $b\se A$ by the set of all edges in $B$; shall refer to this set as the \emph{reservoir} for $e$. This definition ensures that  every maximal $3$-connected minor is equal to a torso, where for each torso edge we pick a suitable label from its reservoir. 

For $j=1,2$, for each of the graphs $\triple_j$, we need to select three edges, so there are at most $n^3$ choices.
For each we can check in polynomial time whether it is a minor of $G_2$. 

Now we have at most quadratically many pairs of graphs $(H_1,H_3)$, so at most polynomially many pairs consisting of $(H_1,\triple_1,H_3)$ and $(H_1,\triple_2, H_3)$. For each such pair we check whether there is a choice of labels from the reservoirs that makes both these triples into temporal triples. Since the reservoirs are disjoint for distinct torso edges, this can be done in linear time. Since $H_1$ and $H_3$ each have at most two embeddings, we can check in constant time whether there are maps $\alpha_1$ and $\alpha_3$. 
This completes the construction of some polynomial-time algorithm for this problem.
\end{pf}

\begin{rem}
Making use of the block-cutvertex decomposition and the decomposition into connected components, it is straightforward to extend this proof to arbitrary temporal triples that need not be $2$-connected, providing a full proof of \autoref{lem:polynomial_disagreeable}. Since the focus of this paper are $2$-connected temporal sequences, we do not burden the reader with this technical extension. It remains open whether there is a linear time algorithm for  \autoref{lem:polynomial_disagreeable}. 
\end{rem}

\begin{rem}[Computing with lists]\label{rem:computing_with_lists}
A list for one graph $G$ assigns relative embeddings to pairs of $3$-blocks of $G$. In this remark, we outline the reasons why computing with lists can be done efficiently. We will consider only lists for one graph, since the case for sequences is analogous.
Each $3$-block of $G$ corresponds to a single bag of the Tutte decomposition, which has linearly-many bags with respect to the number of edges of $G$. In the other direction, each $3$-block can be obtained as a torso of a bag of the Tutte decomposition by choosing a torso edge for each edge in the auxiliary tree that is incident with the vertex that represents that bag. For each of these torso edges, we have linearly many options, and by definition each results in a Tutte-equivalent block.
The definition of list-embeddings is invariant under Tutte-equivalence, so we only ever need to store embeddings of one member of each equivalent class. That is, at most linearly many. We can keep track of which equivalence class each belongs to by always storing the label of the node of the Tutte decomposition tree which represents it. Storing a rotation system for a connected graph also only requires data of size linear in the number of edges. In conclusion, a list requires only data of linear size with respect to the number of edges in the graph.

It is also easy to see how one can check whether two lists are combinable in polynomial-time. Indeed, given a graph $G$ with $3$-blocks $b_1,\ldots,b_n$, and a pair of lists $L$ and $M$, we proceed greedily. For each $i=1,\ldots,n$ in order, we try to make the embedding assigned to $b_i$ in $L$ and $M$ the same. To do so, we are allowed to reorient all of the embeddings of the $3$-blocks in the part of $L$ containing $b_i$ or all of the embeddings in the part of $M$ containing $b_i$. After the embeddings assigned to $b_i$ are equal, we consider $b_i$ `fixed'. That is, we are not allowed to reorient its embeddings in $L$ or $M$ in any further steps. We proceed until either we obtain a combination or we are not able to continue, in which case we conclude that the lists are uncombinable.

\end{rem}

\begin{rem}
In the main result, we believe that the quadratic bound for the increase in the size of the temporal seequence at each step is not best possible. Indeed, in the hypotheses of Lemmas \ref{lem:b-model} and \ref{list_for_bonds} we require a \rcol/\bcol-colouring. However, we claim that it is possible to prove a version of each of these which works for arbitrarily many colours while only increasing the bounds by a constant factor. In this case, the constant $k$ in the bound in \autoref{lem:model_list_workhorse} is replaced by some constant that does not depend on the number of colours. Then we can improve the bound in \autoref{lem:B_1}, \autoref{lem:B_2} and finally the main result \autoref{thm:main}.
Since this construction would in the end also just yield a polynomial algorithm (though with a much better running time), we prioritized readability over a better polynomial running time, and proved this weaker version of \autoref{lem:model_list_workhorse} in the paper. 
\end{rem}

\bibliographystyle{plain}
\bibliography{wliterature}

\end{document}